\documentclass [11pt]{amsart}
\usepackage {amsmath, amssymb, a4wide, epsfig, enumerate, psfrag}
\usepackage[latin1]{inputenc}
\usepackage{verbatim}
\usepackage[all]{xy}

\title[Distribution of critically preperiodic parameters]{Distribution
  of rational maps  with a preperiodic critical point}

\date{\today}

\author{Romain Dujardin \and Charles Favre}

\address{Paris 7 et Institut de Math{\'e}matiques de Jussieu\\
         {\'E}quipe G{\'e}om{\'e}trie et Dynamique\\
         Case 7012, 2  place Jussieu\\
         F-75251 Paris Cedex 05\\
         France}
\email{dujardin@math.jussieu.fr}

\address{CNRS et Institut de Math{\'e}matiques de Jussieu\\
        {\'E}quipe G{\'e}om{\'e}trie et Dynamique\\
         Case 7012, 2  place Jussieu\\
         F-75251 Paris Cedex 05\\
         France}
\email{favre@math.jussieu.fr}

\dedicatory{Dedicated to the memory of Adrien Douady}

\subjclass[2000]{Primary: 37F45, Secondary: 32H50}

\newcommand{\cc}{\mathbb{C}}

\newcommand{\re}{\mathbb{R}}

\newcommand{\zz}{\mathbb{Z}}
\newcommand{\nn}{\mathbb{N}}
\newcommand{\PP}{\mathbb{P}}
\newcommand{\e}{\varepsilon}
\newcommand{\cv}{\rightarrow}

\newcommand{\fr}{\partial}
\newcommand{\om}{\Omega}
\newcommand{\set}[1]{\left\{#1\right\}}

\newcommand{\abs}[1]{\left\vert#1\right\vert}

\newcommand{\rest}[1]{ \arrowvert_{#1}}

\newcommand{\unsur}[1]{\frac{1}{#1}}
\newcommand{\cA}{\mathcal{A}}
\newcommand{\el}{\mathcal{L}}

\newcommand{\C}{\mathcal{C}}
\newcommand{\Per}{\mathrm{Per}}
\newcommand{\preper}{\mathrm{Preper}}

\newcommand{\hL}{\widehat{\Lambda}}
\newcommand{\hT}{\widehat{T}}
\newcommand{\ho}{\hat{\omega}}
\newcommand{\hf}{\hat{f}}
\newcommand{\hg}{\hat{g}}
\newcommand{\rond}{\hspace{-.1em}\circ\hspace{-.1em}}

\renewcommand{\=}{:=}

\newcommand{\itm}{\item[-]}

\DeclareMathOperator{\supp}{Supp}

\DeclareMathOperator{\ord}{ord}

\newtheorem{prop} {Proposition} [section]
\newtheorem{thm}[prop] {Theorem} 
\newtheorem{defi}[prop] {Definition}
\newtheorem{lem}[prop] {Lemma}
\newtheorem{cor}[prop]{Corollary}
\newtheorem{prop-def}[prop]{Proposition-Definition}
\newtheorem{theo}{Theorem} 
\newtheorem{coro}[theo]{Corollary}
\newtheorem{propo}[theo]{Proposition} 
\theoremstyle{remark}
\newtheorem{exam}[prop]{Example}
\newtheorem{rmk}[prop]{Remark}

\begin{document}

\begin{abstract}
  Let $\{f_\lambda\}$ be any algebraic family of rational maps
  of a fixed degree, with a marked critical point $c(\lambda)$.  We
  first prove that the hypersurfaces of parameters for which $c(\lambda)$ is
  periodic converge as a sequence of positive closed $(1,1)$ currents
  to the bifurcation current attached to $c$ and defined by
  DeMarco~\cite{DeM1}. We then turn our attention to the parameter
  space of polynomials of a fixed degree $d$. By intersecting the
  $d-1$ currents attached to each critical point of a polynomial,
  Bassaneli and Berteloot~\cite{bas-ber} obtained a positive measure
  $\mu_\mathrm{bif}$ of finite mass which is supported on the
  connectedness locus. They showed that its support is included in the
  closure of the set
  of parameters admitting $d-1$ neutral cycles. We show that the support
  of this measure is precisely the closure of the set of strictly
  critically finite polynomials (i.e. of Misiurewicz points).
\end{abstract}

\maketitle
\setcounter{tocdepth}{1}

\tableofcontents

\newpage


 \section*{Introduction} \label{sec:intro}
 
 It is a central problem in dynamics to understand how the dynamics of
 a map can change under perturbation. In the context of rational maps
 of the Riemann sphere, this question has received a lot of attention,
 and is still a very active area of research. The seminal paper of
 Ma{\~n}{\'e}-Sad-Sullivan~\cite{MSS} has first paved the way to understand
 structural stability of general holomorphic dynamical systems, and
 was completed a few years later by the construction of a Teichm{\"u}ller
 theory for rational maps in~\cite{mcmsul}. Besides these general
 results, many particular families have been studied in detail. Among
 them are the family of quadratic polynomials~\cite{DH,mandel}, the
 set of cubic polynomials~\cite{BrH1,BrH2,Kiwi}, or the space of
 quadratic rational maps~\cite{rees1,rees2,milnor1}.

\medskip

\begin{center} $\diamond$ \end{center}
\medskip

Let $(\Lambda, f)=\set{f_\lambda, \lambda\in \Lambda}$ be a
holomorphic family of rational maps of degree $d\ge2$, parameterized
by a smooth connected complex manifold $\Lambda$ (of any dimension).
We say that a critical point is {\em marked} if it can be followed
holomorphically along $\Lambda$ (see Section~\ref{subs:active} for
precise definitions). Following the terminology of McMullen
\cite{McM1}, a critical point is {\em passive at $\lambda_0\in
  \Lambda$} if the family $\{f^n_\lambda\, c(\lambda)\}_n$ is normal
in some neighborhood of $\lambda_0$.  Otherwise, $c$ is {\em active}.
It follows from~\cite{MSS} that a rational map is unstable if and only
if at least one of its critical points is active.

Following L.DeMarco~\cite{DeM1} it is possible to associate a natural
positive closed (1,1) current $T$ to a marked critical point
$\lambda\mapsto c(\lambda)$. The support of $T$ is the {\em activity
  locus} of $c(\lambda)$.  Her method requires  lifting  the
situation to $\mathbb{C}^2\setminus\{0\}$. Our first aim is to present
a coordinate free presentation of her results --see
Proposition-Definition~\ref{def:bifurcation} and
Theorem~\ref{thm:support} below. Our presentation is quite classical
in higher dimensional holomorphic dynamics and owes much to the work
of V.Guedj~\cite{guedj1,guedj2}, see also~\cite{DS}.

We then study the distribution of parameters for which $c(\lambda)$ is
(pre-)periodic. A classical normal families argument using Montel's
Theorem implies that if $c(\lambda)$ is active at $\lambda_0$, then
arbitrarily close to $\lambda_0$ there are parameters for which
$c(\lambda)$ is preperiodic.  We prove a {\em quantitative} version of this
fact under the form of an equidistribution theorem, which is a parameter space
analog of the following classical theorem: the periodic points of a
fixed rational map equidistribute towards the maximal entropy measure,
see~\cite{Tor, Lyu}.

For any $n>k\ge0$, let us introduce the set $\mathrm{Per}(n,k)$ of parameters
for which $f^n_\lambda\, c(\lambda) = f^k_\lambda \, c(\lambda)$. It
is either equal to the whole variety $\Lambda$ or it is a
hypersurface. In the former case, we call the family \emph{trivial}:
the activity locus is then empty, and the current $T$ equals $0$.  In
the latter case, we consider $[\mathrm{Per}(n,k)]$ the current of
integration over the divisor $\mathrm{Per}(n,k)$ (each component is
counted with its multiplicity as a solution of 
$f^n_\lambda\, c(\lambda) = f^k_\lambda \, c(\lambda)$).  
\begin{theo}\label{t:cvg}
  Let $(\Lambda, f, c)$ be a non-trivial holomorphic family of
  rational maps of $\PP^1$ of degree $d\ge 2$ with a marked critical
  point. If $\Lambda$ is quasi-projective, then for any
  sequence of integers $0\le k(n) < n$, the following convergence
  statement holds:
 $$\lim_{n\to\infty} \frac{[\Per(n,k(n))]}{d^n+d^{(1-e)k(n)}}  = T$$
 where $e\in\{0,1\}$ is the cardinality of the exceptional set of a
 generic map $f_\lambda$.
\end{theo}
In other words, $e=0$ if, for some (hence for a generic) $\lambda$,
$f_\lambda$ is not M{\"o}bius-conjugate to a polynomial. On the other
hand, $e= 1$ means that all $f_\lambda$ are polynomials. We refer to
Section~\ref{subs:exc} for a brief discussion on exceptional points in
families of rational maps. 

We also prove (Theorem \ref{t:quasiproj}) that, except for some
trivial cases, in a family of rational maps parameterized by a
quasiprojective manifold, the current $T$ is always non trivial. See
the comments following Theorem~\ref{p:motion} for more details.

A natural method for proving the theorem is the following. We already
mentioned that for any rational map $f$, periodic points
equidistribute towards the measure of maximal entropy. These results
can be put together to yield a convergence theorem in $\Lambda\times
\PP^1$. There exists a positive closed $(1,1)$ current $\hT$ on
$\Lambda\times \PP^1$ such that $d^{-n}[\{ (\lambda,z),\, f^n_\lambda
(z) = z] \to \hT$.  We can now look at the hypersurface $\Gamma \= \{
(\lambda, c(\lambda))\} \subset \Lambda \times \PP^1$. The projection
$\pi$ onto the first factor gives an analytic isomorphism from
$\Gamma$ to $\Lambda$, and the constructions are done in such a way
that $T= \pi_* ( \hT|_\Gamma)$.  Observe that $d^{-n}[\{
(\lambda,z),\, f^n_\lambda (z) = z]|_\Gamma \to \hT|_\Gamma$ is
equivalent to the statement of our theorem.  However, it is important
to note that the convergence of currents is not strong enough to
ensure the convergence on $\Gamma$.  Our result thus needs a separate
argument.

Our proof relies on potential theory, and bears some similarity with
the argument designed by Brolin for the proof of his celebrated
equidistribution theorem \cite{Bro}.  This results in the appearance
of the assumption of quasiprojectivity, because we use the maximum principle
for plurisubharmonic functions in the course of the proof.

The theorem applies immediately in a number
of interesting examples, including the family of polynomials of degree
$d$ (Section~\ref{sec:poly}), and the family of rational maps of degree
$d$ (Section~\ref{sec:rational}). In the latter family, Bassanelli-Berteloot
have recently used the description of the bifurcation current in terms of  Lyapunov exponents
to prove equidistribution results of parameters admitting a periodic
cycle with a fixed multiplier, see~\cite{BB2}.

By applying our convergence theorem to
the one-parameter family of unicritical polynomials of degree $d$, we
get the following corollary which is due to  Levin~\cite{levin}.
\begin{coro}[\cite{levin}] \label{cor:unicritical}
Let $\mathcal{M}_d$ be the set of complex numbers $c$ for which $f_c(z)=z^d+c$
has connected Julia set (the Mandelbrot set). Then
$$
\lim_{n\to\infty} \frac1{d^n}\sum_{f_c^n(0)=0} \delta_c = \mu
$$
where $\mu$ is the harmonic measure on $\mathcal{M}_d$.
\end{coro}
It is remarkable that this statement can be obtained by arithmetic
methods based on height theory.  It is a consequence of~\cite{Aut},
and it was explicitely stated (when $d=2$)
in~\cite[Theorem~8.13]{baker-hsia}.  Using height theory presents the
advantage to yield a precise estimate on the speed of convergence.
 For instance, a proof of the following
result is given in~\cite{FRL}.
\begin{theo}[\cite{FRL}]\label{cor:unic} With notation as in
  Corollary~\ref{cor:unicritical}, let $F_n\subset \mathbb{C}$ be a
  sequence of disjoint finite sets, invariant under the absolute
  Galois group of $\mathbb{Q}$, and included in the union
  $\bigcup_{n\neq k} \{ f_c^n(0) = f_c^k(0)\}$. Then for any compactly
  supported $C^1$ function $\varphi$, we have
  $$
  \left| \frac1{|F_n|} \sum_{c\in F_n} \varphi(c) - \int \varphi\,
    d\mu \right| \le C\, \left(\frac{\log |F_n|}{|F_n|}\right)^2 \times \sup\{
  |\varphi|, |\varphi'| \} ~,$$
  where $|F|$ denotes the cardinality of
  $F$.
\end{theo}
Unfortunately, it is not clear how to extend this method to higher
dimensional parameter spaces.

\medskip

Another issue in the proof of the convergence Theorem \ref{t:cvg} is
whether preperiodic critical points are active or not. Of course if a
critical point is preperiodic to a repelling cycle it is active, and
if it is preperiodic to an attracting cycle it is passive. In the
neutral case, if for instance $\Lambda$ is the space of all rational 
maps  (or all polynomials), such a critical point can 
directly be perturbed to become prerepelling, so it is active.
In the case of a general family $\Lambda$,
the situation is more delicate, and it seems that this question, although natural, 
has not been previously
addressed.

\begin{theo}\label{p:motion}
Let $(\Lambda, f, c)$ be any holomorphic family of
  rational maps of $\PP^1$ of degree $d\ge 2$ with a marked critical
  point.  Assume $U\subset \Lambda$ is a connected open subset
  where $c$ is passive. Then exactly one of the following cases hold.
\begin{enumerate}
\item $c$ is never preperiodic in $U$. In this case the closure of the
  orbit of $c$ can be followed by a holomorphic motion.
\item $c$ is persistently preperiodic in $U$. 
\item The set of parameters for which $c$ is preperiodic is a closed
  subvariety in $U$. Moreover, either there exists a persistently
  attracting (possibly superattracting) cycle attracting $c$
  throughout $U$, or $c$ lies in the interior of the linearization
  domain associated to a persistent irrationally neutral periodic
  point. 
\end{enumerate}
\end{theo}
This result is in fact a consequence of a purely local statement,
Theorem \ref{t:local}.  As a consequence of these techniques, we give in
Section \ref{subs:quasiproj} the following generalization of a theorem of McMullen
\cite{McM2}: for any {\em algebraic} family of rational maps with a
marked critical point $(\Lambda, f, c)$, either $f$ is constant, or
$c$ is persistently preperiodic, or the activity locus of $c$ is non
empty (see Theorem \ref{t:quasiproj} for a precise statement).

\begin{center} $\diamond$ \end{center}
\medskip

When all critical points are marked, the bifurcation
locus is the union of the various activity loci of the critical
points. This can be refined by introducing successive bifurcation
loci, indexed by the number of critical points being active.

A natural expectation is that if $k$ marked critical points are active
at $\lambda_0$, there should be a nearby parameter where the $k$
critical points are preperiodic. This is of particular interest when
$k$ is maximal, so that the perturbed map becomes critically finite.
However, there is no reasonable analogue of Montel's Theorem for
sequences of holomorphic mappings in higher dimension (due to the
Fatou-Bieberbach phenomenon), and as it turns out, this expectation is
wrong. Indeed there exist cubic polynomials with both critical points
active such that for every perturbation, one of the two critical
points is attracted by an attracting cycle (see
Example~\ref{ex:douady}).  We thank A. Douady for kindly communicating
this example to us.

It is one of the main ideas in higher dimensional holomorphic dynamics
that pluripotential theory and the use of positive closed currents can
serve as a natural substitute to Montel's Theorem.  Later on we shall see that
in order to obtain a correct statement, we need to replace the locus
where $k$ critical points are active by the support of the wedge
product of $k$ suitable bifurcation currents.

\medskip

In section \ref{sec:poly}, we implement this strategy in the space
$\mathcal{P}_d$ of polynomials of degree $d$ with all critical points
marked, which, up to finite branched cover, is biholomorphic to
$\cc^{d-1}$. In this space we can consider $d-1$ positive closed
currents of bidegree (1,1): $T_0,\ldots ,T_{d-2}$, associated to the
marked critical points, as well as the so called {\em bifurcation
  current} $T_{\rm bif}=\sum T_i/(d-1)$ introduced by
DeMarco~\cite{DeM1} and its successive powers $(T_{\rm bif})^{\wedge
  k}$, $k=1, \cdots, d-1$ introduced by Bassanelli and
Berteloot~\cite{bas-ber}. Although both papers are concerned with the
more difficult situation of rational maps of degree $d$, their results
apply equally in the context of polynomial maps (see also \cite{mai}
for related definitions in a wider context).

The support of $(T_{\rm bif})^{\wedge k}$ is contained in the set of
parameters where $k$ critical points are active. This inclusion can be
strict, as shown by the already mentioned Example~\ref{ex:douady}.
Bassanelli and Berteloot proved that near any point in $\supp(T_{\rm
  bif}^{\wedge k})$, there is a parameter with $k$ neutral periodic
orbits. Here, we obtain the following result:
\begin{theo}\label{thm:approx-all}
  There exists a sequence of codimension $k$ algebraic subvarieties
  $W_n$ (not necessarily reduced) and positive real numbers $\alpha_n$
  such that $\alpha_n [W_n]$ converges to $(T_{\rm bif})^{\wedge k}$,
  and $W_n$ is supported in the set of parameters for which $k$
  critical points are preperiodic.
\end{theo}
In particular near every $\lambda\in \supp(T_{\rm bif}^{\wedge k})$
there is a parameter for which $k$ critical points are preperiodic.
Theorem~\ref{thm:approx-all} is proved by successive applications of
Theorem~\ref{t:cvg} on the subvarieties where a fixed number of
critical points is periodic with a fixed period.  

\medskip

The case where $k=d-1$ deserves special attention.  In this case the
suitably normalized positive measure $\mu_{\rm bif} \= c\, (T_{\rm
  bif})^{\wedge d-1}$ is a probability measure, supported on the
boundary of the connectedness locus, which is compact in
$\mathcal{P}_d$.  Recall that a polynomial is said to be {\em of
  Misiurewicz type} if all its critical points are preperiodic to
repelling cycles. It is classical that if all critical points are
strictly preperiodic, then the Misiurewicz property holds. From
Theorem \ref{thm:approx-all} we immediately get the following
corollary.

\begin{coro}\label{cor:mis}
  The support of the bifurcation measure $\mu_{\rm bif}$ is contained
  in the closure of Misiurewicz parameters.
\end{coro}

\smallskip

We obtain several interesting characterizations of the bifurcation
measure, which are well known when $d=2$ (and more generally for
unicritical polynomials).  First, when $d=2$, $\mu_{\rm bif}$ is the
harmonic measure of the Mandelbrot set. To give a precise statement in
higher dimension, we need to use a parameterization of $\mathcal{P}_d$
by the affine space $\cc^{d-1}$ (see Section \ref{sec:param} for more
details).

\begin{propo}\label{prop:equilibrium}
The bifurcation measure is the pluricomplex equilibrium
measure of the connectedness  locus $\mathcal{C}\subset \cc^{d-1}$. As a 
consequence, $\supp(\mu_{\rm bif})$ is the Shilov boundary of $\mathcal{C}$.
\end{propo}
Being the pluriharmonic measure of the connectedness locus shows that
this measure is natural from the point of view of complex analysis.

\medskip

Our next result is a characterization of the measure $\mu_{\rm bif}$
in $\mathcal{P}_d$ as the landing measure of a family of
external rays. Let us explain how these rays are defined.  Let $P$ be
a polynomial whose Green function takes the same value $r>0$ at all
critical points.  The set $\Theta$ of external angles of rays landing
at the critical points gives a natural way to describe the
combinatorics of $P$, see \cite{Go1,BFH,kiwi-portrait} where $\Theta$
is referred to as the \emph{critical portrait} of $P$.  Now, we may
deform $P$ in the shift locus by leaving $\Theta$ unchanged and
letting $r$ vary in $\mathbb{R}_+^*$ --this is the operation of {\em
  stretching}, as defined in \cite{BrH1}.  This defines a ray in
parameter space, corresponding to $\Theta$.

We may now consider the set
$\mathsf{Cb}$ of all possible combinatorics/critical portraits.  
This space 
is a compact finite
dimensional ``manifold'', endowed with a natural measure
$\mu_{\mathsf{Cb}}$ arising from the translation structure on the
angle space $\re/\zz$ (see Proposition~\ref{p:trans} and
Definition~\ref{d:meas}).  As observed by \cite{BMS}, Fatou's Theorem
implies that almost every ray lands when $r\cv0$
(Proposition~\ref{p:landing}).  We may thus define a measurable {\em
  landing map} $e: \mathsf{Cb} \cv \mathcal{C}\subset \mathcal{P}_d$.
The basic link between critical portraits and the bifurcation measure
is given by the following
\begin{theo}\label{thm:landing} The image of $\mu_{\mathsf{Cb}}$ under
  landing is $\mu_{\rm bif}$, that is, 
$e_*\mu_{\mathsf{Cb}} = \mu_{\rm bif}$.
\end{theo} 

The landing of external rays was extensively studied by J.~Kiwi. He
obtained a fundamental continuity result, which we now briefly
describe.  A combinatorics $\Theta$ is of {\em Misiurewicz type} if
all angles in $\Theta$ are strictly preperiodic under multiplication
by $d$. Kiwi's result \cite[Corollary 5.3]{kiwi-portrait} is that the
landing map $e$ is continuous at Misiurewicz parameters. This
generalizes to higher degrees the well known theorem by Douady and
Hubbard that rational external rays of the Mandelbrot set land.

This, combined  with our description of $\mu$ in
terms of external rays, allows us to give more properties of $\supp(\mu_{\rm
  bif})$. 
\begin{theo}\label{thm:supp}
Every Misiurewicz parameter lies inside $\supp(\mu_{\rm
  bif})$.
\end{theo}

In particular the inclusion in Corollary~\ref{cor:mis} is an equality.

\medskip

Finally, the landing Theorem~\ref{thm:landing} allows us to give some
dynamical properties of $\mu_{\rm bif}$-almost every polynomial. 
\begin{theo}\label{t:CE}
  The Topological Collet-Eckmann property holds for  $\mu_{\rm
    bif}$-almost every polynomial $P$.
\end{theo}
This result in turn  implies
\begin{coro}\label{cor:CE}
For a $\mu_{\rm bif}$-generic polynomial $P$, we have that:
\begin{itemize}
  \itm all cycles are repelling; \itm the orbit of each critical
  point is dense in the Julia set; \itm $K_P=J_P$ is locally connected
  and has Hausdorff dimension strictly less than 2.
\end{itemize}
\end{coro}
The Topological Collet-Eckmann property (TCE for short) is a way to
estimate quantitatively the recurrence of critical points. We refer
to~\cite{prs} for (many) other characterizations of the TCE condition and
references.  In the case of unicritical
polynomials, the TCE condition is equivalent to the more standard
Collet-Eckmann condition. In this case, Theorem~\ref{t:CE} is due to
Graczyk-\'Swi\c{a}tek \cite{gs} and Smirnov \cite{smirnov}.

\begin{center} $\diamond$ \end{center}
\medskip

Let us close this introduction by indicating the structure of
this article.  In Section~\ref{sec:local}, we prove a local result,
Theorem~\ref{t:local}, describing the behaviour of a passive point
which lands on a periodic cycle. This result is the key to the proof
of Theorem~\ref{p:motion}.  Sections~\ref{sec:gen} to~\ref{subs:cv}
deal with general families of rational maps with a marked critical
point.  We begin in Section~\ref{sec:gen} with some generalities on
passive and active points in families of rational maps, and briefly
discuss on exceptional points. We also include the proof of
Theorem~\ref{p:motion} and a description of algebraic families of
rational maps with a marked and passive critical point in the spirit
of~\cite{McM2}.  Section~\ref{sec:bif} is devoted to the construction
of the bifurcation current. We show that its support is the activity
locus.  Section~\ref{subs:cv} contains the proof of a slightly more
general version of our convergence result. 

Sections~\ref{sec:param} to \ref{sec:external} are devoted to the
parameter space of polynomials. Its basic properties, as well as the
parameterization by $\cc^{d-1}$, are described in
Section~\ref{sec:param}. In Section~\ref{sec:poly}, we describe the
structure of the higher bifurcation currents and prove
Theorem~\ref{thm:approx-all} and its corollary.
Section~\ref{sec:external} is devoted to the description of $\mu_{\rm
  bif}$ in terms of external rays, leading to
Theorems~\ref{thm:landing}, \ref{thm:supp} and \ref{t:CE} and
Corollary~\ref{cor:CE}.

We conclude the paper with Section~\ref{sec:rational}, where we show
how to extend to the space of rational maps some of the results of 
Section~\ref{sec:poly}.
\medskip

\noindent {\bf Acknowledgements.} 
It is a pleasure to thank the ACI-jeunes chercheurs ``Dynamique des
applications polynomiales'' for providing a very stimulating working
atmosphere. We thank Mattias Jonsson for several
discussions concerning the construction of the bifurcation current,
and Adrien Douady and Tan Lei for their careful attention and their
explanations.

\newpage


\section{Families of holomorphic germs and periodic points}\label{sec:local}

The main results of the present article deal with rational or
polynomial mapppings on the Riemann sphere. In the course of the proof
of Theorem~\ref{t:cvg}, we shall however need a result of more local
nature. Our aim in this section is to describe completely the set of 
possible  situations where a holomorphically varying point falls into a periodic
cycle but does not present any bifurcation (see Theorem~\ref{t:local}
below). We believe that this may have independent interest.

Before stating the theorem, let us introduce some terminology.  A
holomorphic family of holomorphic maps defined on the unit disk
$\Delta$ and parameterized by a complex manifold $\Lambda$ is a
holomorphic map $f: \Lambda\times \Delta \to \mathbb{C}$.  In general,
we write $f(\lambda, z) = f_\lambda(z)$. In this section, we always assume that
$f_\lambda (0) =0$ for all $\lambda$. Notice that for fixed $\lambda$
and $n$, the set of points $U_n(\lambda)\subset \Delta$ for which
$f_\lambda, \cdots, f^n_\lambda$ are well-defined forms a decreasing
sequence of open neighborhoods of $0$, whose intersection might be
reduced to the origin.

Recall that a map $f: \Delta \to \mathbb{C}$ fixing the origin is
linearizable if there exists a holomorphic and locally invertible germ $\phi$ such that $\phi
\circ f (z) = \mu\, \phi(z)$ where $\mu = f'(0)$. Any germ with
 $|\mu|\neq 1$ is linearizable, see~\cite{CG,mibook}. When $\mu^l =1$, then $f$
is linearizable iff $f^l$ is the identity map. When $|\mu|=1$ and
$\mu$ is not a root of unity, the situation is much more
complicated, but any germ is at least \emph{formally} linearizable,
see~\cite[Problem 8.4]{mibook}. The domain of linearization is the maximal
\emph{open} set $U$ for which there exists a biholomorphism $\phi$
into $\mathbb{C}$ conjugating $f$ to $w \mapsto \mu w$.
\begin{thm}\label{t:local}
  Let $f_\lambda$ be any holomorphic family of holomorphic maps
  parameterized by a connected complex manifold $\Lambda$. Suppose that
  each $f_\lambda$ is defined on the unit disk with values in
  $\mathbb{C}$, and leaves the origin fixed i.e.  $f_\lambda(0) =0$.
  Let $\lambda \mapsto p(\lambda)$ be any holomorphic map such that
  $p(\lambda_0) =0$ for some parameter $\lambda_0$.
 
  Assume that for all $n\in \mathbb{N}$, the function $f^n_\lambda\,
  p(\lambda)$ is well-defined and takes its values in the unit disk.
  Then one of the following three cases holds.
\begin{enumerate}
\item For every $\lambda\in \Lambda$, the point $0$ is attracting or
  superattracting, and $p(\lambda)$ lies in the (immediate) basin of
  attraction of $0$.
\item The point $p$ is periodic for all parameters, i.e.
  $f^l_\lambda\, p(\lambda) = p(\lambda)$ for some $l$ and all
  $\lambda\in \Lambda$.
\item The multiplier of $f_\lambda$ at $0$ is constant and equals
  $\exp(2i\pi \theta)$, with $\theta \in \mathbb{R}\setminus
  \mathbb{Q}$. For all $\lambda\in \Lambda$, the map $f_\lambda$ is
  linearizable and $p(\lambda)$ lies in the interior of the domain of
  linearization of $f_\lambda$.
\end{enumerate}
\end{thm}

\begin{proof}
  Suppose first that $0$ is an attracting fixed point of
  $f_{\lambda_0}$.  Our aim is to show that for all parameters, the
  fixed point $0$ remains attracting, and $p(\lambda)$ is attracted
  towards $0$.  Note first that there exists a neighborhood $U$ of
  $0\in \Lambda$, and a fixed disk $D$ containing $0$ in the dynamical
  plane such that $f_\lambda (D)$ is relatively compact in $D$ for all
  $\lambda \in U$. In particular, $f^n_\lambda (z) \to 0$ for all $z
  \in D$.
  
  Now by assumption the sequence $\{ f^n_\lambda\, p(\lambda)\}_{n\in
    \mathbb{N}}$ forms a normal family. Any cluster value of this
  sequence vanishes identically on the open set of parameters in $U$
  for which $p(\lambda) \in D$. So $ f^n_\lambda\, p(\lambda)$
  actually converges to zero uniformly on compact sets on $\Lambda$.
  We infer that $|f'_\lambda(0)| \le 1$ for all $\lambda$. As the
  multiplier of $f_{\lambda_0}$ at $0$ has modulus $<1$ and $\lambda
  \mapsto f'_\lambda(0)$ is holomorphic, the Maximum Principle implies
  that $0$ is an attracting fixed point for $f_\lambda$ for every
  $\lambda$. This shows that Case~(1) holds.

\medskip

From now on, we assume that the multiplier of $f_{\lambda_0}$ at
$0$ has modulus $\ge1$. For the sake of simplicity, we also assume
that $\Lambda$ has dimension one.  We explain at the end of the proof
how to deal with the general case.  The key computation is contained
in the next lemma. It is a classical result in the case $\mu=1$, and
serves as a basis for the definition of the iterative logarithm,
see~\cite{ecalle}.
\begin{lem}\label{l:comp}
  Write $f_\lambda(z) = \mu z + \sum_{k+l \ge 2} a_{kl} \lambda^k z^l$
  with $\mu \in \mathbb{C}^*$ and $a_{kl} \in \mathbb{C}$.  Then for
  any integer $n$, we have $f^n_\lambda(z) = \mu^n z + \sum_{k+l \ge
    2} a^n_{kl} \lambda^k z^l$, with
\begin{equation}\label{e:poly}
a^n_{kl} = \sum_{r=0}^l \mu^{rn} P_{rkl}(n)
\end{equation}
where $P_{rkl}(n)$ is a \emph{polynomial} in $n$.
\end{lem}
The proof of this lemma will be given at the end of this section. We now
continue with the proof of the theorem.

For each $n$, write the expansion of $\lambda$ of
$f^n_\lambda\,p(\lambda)$ into increasing powers as
$f^n_\lambda\,p(\lambda)= \sum_{j\ge 1} q_j^n \lambda^j$ with
$q_j^n\in\mathbb{C}$.  By the preceding lemma, $f_\lambda^n\,
p(\lambda)$ equals $ \mu^n p(\lambda) + \sum_{k+l \ge 2} a^n_{kl}
\lambda^k (p(\lambda))^l$.  Identifying the $\lambda^j$ terms of both
expressions, we infer that $q_j^n = c\cdot \mu^n + \sum_{k,l\le j
}c_{kl} a^n_{kl}$ for some constants $c, c_{kl}$. The order of
vanishing $\ord_0 (p)$ of $p$ at $0$ is greater than $1$, hence
$\ord_0 (p^l) \ge l$. This explains why we can take $l\le j$ in the
above sum.  In particular, for any $j$ , we can write
\begin{equation}\label{e:pol}
q_j^n = \sum_{r=0}^j \mu^{rn} Q_{rj}(n)~,
\end{equation}
where the $Q_{rj}$ are polynomials.

We now translate our main assumption on $p$ into estimates on the
coefficients of the power series expansion of $f^n_\lambda\,
p(\lambda)$. By assumption $f^n_\lambda\, p(\lambda)$ is a family of
holomorphic functions with values in the unit disk. So the Cauchy estimates
imply that  for each $j$
\begin{equation}\label{e:bdd}
\sup_n |q_j^n| \le C(j)<+\infty~.
\end{equation}
The proof of the theorem is now based on the comparison
between~\eqref{e:pol} and the estimates~\eqref{e:bdd}.  We proceed by
a case by case analysis, assuming first that $|\mu|>1$; then that $\mu$ is a
root of unity; and finally dealing with the case of $\mu = \exp(2i\pi
\theta)$, with $\theta$ irrational.

\bigskip

Suppose first $|\mu| >1$ (this case is classical). Fix any integer
$j\ge 1$.  Then $q_j^n \sim c\cdot n^d \mu^{rn}$ for some
$c\in\mathbb{C}^*$ and for $r\in\mathbb{N}$ maximal such that $Q_{rj}$
is non-zero and $d = \deg (Q_{rj})$.  The estimates~\eqref{e:bdd}
imply $d=0$ and $r=0$, hence $q_j^n$ is constant independent on $n$.
We conclude that $f^n_\lambda\, p(\lambda)=p(\lambda)$ for all $n$.
In particular, $f_\lambda\, p(\lambda) =p(\lambda)$, and Case~(2) of
the theorem holds.  Notice that since $0$ is a simple root of the
equation $f_{\lambda_0}(z)-z = 0$, in fact $p(\lambda)\equiv 0$.

\bigskip

For the remaining part of the proof, we assume $|\mu|=1$, and write $\mu =
\exp(2i\pi\theta)$ for some real number $\theta$.

\medskip

When $\theta$ is a rational number, $\mu$ is a root of unity and
$\mu^l =1$ for some integer $l$. Equation~\eqref{e:pol} implies that
$q_{j}^{ln}$ is actually \emph{polynomial} in $n$. The
estimates~\eqref{e:bdd} now show that for each $j$
this polynomial is constant. We
thus conclude that $f^{nl}_\lambda\, p(\lambda) = \sum q_j \lambda^j$ for
all $n\ge0$.  In particular, $f^lp(\lambda) = p(\lambda)$. As before
we are in Case~(2) of the theorem. Notice however that in general $p$
need not be fixed, and that even if it is fixed, it may be different
from $0$.

For later reference, let us summarize this discussion in a lemma.
\begin{lem}\label{l:multi-first}
If the multiplier of $f_{\lambda_0}$ at $0$ has modulus $>1$ or is a
root of unity, then $f^l_\lambda p(\lambda) = p(\lambda)$ for some
integer $l$ and all $\lambda$.
\end{lem}

Finally, suppose that $\theta$ is an irrational number, and fix some
integer $j\ge1$. Let $d= \max_r \deg (Q_{rj})$. Assume that $d>0$ and
write $Q_{rj}(T) = Q^d_{rj}T^d+ O(T^{d-1})$, so at least one
of the $Q^d_{rj}$ is non zero when $r$ ranges from 0 to $j$. 
Now for any complex number $\zeta$ of
norm $1$, there exists an increasing sequence of integers $n_k$ such
that $\mu^{n_k}$ converges towards $\zeta$. Fix $\zeta$ such that
$\sum_0^j\zeta^r Q^d_{rj}\neq 0$. We infer that when $k\cv\infty$,
$q^{n_k}_j \sim n_k^d \times \sum_0^j\zeta^r Q^d_{rj}$, which
contradicts \eqref{e:bdd}. We thus conclude that $d=0$. In other
words, for every $j$ there exists a polynomial $Q_j$ such that
$$q_j^n = Q_j(\mu^n)~.$$

When the point $p$ is fixed for all parameters, Case~(2) of the
theorem holds so there is nothing to prove.  We thus assume that
outside some discrete subset of $\Lambda$, $f_\lambda\, p(\lambda)
\neq p(\lambda)$ (recall that $\dim \Lambda =1$). Our aim is to
prove that we are in Case~(3) of the theorem. We first prove the
result in the neighborhood of $\lambda_0$.
\begin{lem}\label{l:local}
  There exists a neighborhood $U$ of $\lambda_0$ such that for all
  $\lambda\in U$, $f_\lambda$ is linearizable at 0, with multiplier
  $\mu=\exp(2i\pi\theta)$ independent of $\lambda$ and $p$ belongs to
  the interior of the domain of linearization. Moreover, $f^n_\lambda
  \, p(\lambda) \neq f^m_\lambda \, p(\lambda)$ for all $n\neq m$ and
  all $ U\ni\lambda \neq \lambda_0$.
\end{lem}
We then globalize this result. First, it is  clear that the
multiplier of $f_\lambda$ at $0$ is constant equal to $\mu$ for all
$\lambda \in \Lambda$. Next we have:
\begin{lem}\label{l:discrete}
  The set of points for which $f^n_\lambda p(\lambda) = f^m_\lambda
  p(\lambda)$ for some $n >m$ is \emph{discrete} in $\Lambda$.
\end{lem}
Let $\mathcal{F}$ be the set of points for which $p(\lambda)$ is
preperiodic. The previous lemma shows it is discrete. On $\Lambda
\setminus \mathcal{F}$, the points $f^n_\lambda p(\lambda)$, $n\in
\nn$, move holomorphically and without collision, thus there is a
holomorphic motion of the orbit of $p(\lambda)$ parameterized by
$\Lambda \setminus \mathcal{F}$. By~\cite{MSS}, it automatically
extends to the closure of the orbit, that we denote by
$\gamma_\lambda$.  For $\lambda \notin \mathcal{F}$ close enough to
$\lambda_0$, $\gamma_\lambda$ is a circle surrounding $0$ by
Lemma~\ref{l:local}. As $\Lambda\setminus \mathcal{F}$ is connected,
we infer that for all $\lambda \in\Lambda\setminus \mathcal{F} $, the
map $f_\lambda$ admits an invariant quasicircle surrounding $0$, on
which the dynamics is conjugate to an irrational rotation. Because of
the latter property, we further deduce that this quasicircle does not
contain $0$.

We are now in position to prove that the fixed point $0$ is
linearizable for any $\lambda\in \Lambda$. If $\lambda \in
\Lambda\setminus \mathcal{F}$, the connected component of
$\Delta\setminus\gamma_\lambda$ containing $0$ is a neighborhood of
$0$ that does not escape under iteration, hence $0$ is a linearizable
fixed point. If on the other hand $\lambda\in \mathcal{F}$, take a
small loop $\ell$ around $\lambda$, and avoiding $\mathcal{F}$. For
the parameters $\lambda'\in \ell$, we deduce by continuity the
existence of a disk $D(0,r)$ of a fixed size in dynamical space that
does not escape under iteration. By the maximum principle, this also
holds true for $\lambda$, and we conclude that 0 is also linearizable
for $f_\lambda$.

We summarize this discussion in the following lemma.
\begin{lem}\label{l:closure}
  The map $f_\lambda$ is linearizable for all parameters, and the point
  $p(\lambda)$ lies in the \emph{closure} of the domain of linearization of
  $f_\lambda$. In particular, either $p(\lambda) =0$ or its orbit is infinite.
\end{lem}
Finally we have:
\begin{lem}\label{l:interior}
For every $\lambda$ the point $p(\lambda)$ lies in the \emph{interior} of the
domain of linearization of $f_\lambda$.
\end{lem}
A proof is given below.  This concludes the proof of
Theorem~\ref{t:local} in  case the parameter space has dimension
one. The same proof gives the result in any dimension. Namely if the
multiplier $\mu$ of $f_{\lambda_0}$ at $0$ has norm $<1$ then
$f^n_\lambda \,p(\lambda)$ converges to an attracting fixed point for
all parameters. We are thus in Case~(1). When $|\mu|>1$, $p(\lambda)
=0$ for each disk passing through $\lambda_0$ hence everywhere. When
$\mu^l =1$, then $f^l_\lambda p(\lambda) = p(\lambda)$ for all disks
containing $\lambda_0$ hence everywhere. In both situations, we are in
Case~(2). Finally when $|\mu|=1$ is not a root of unity and we are not
in Case~(2), then for any $\lambda\in \Lambda$ 
there exists an immersed  holomorphic disk
$D\subset \Lambda$ containing $\lambda_0$ and $\lambda$  (see for
instance~\cite{fornaess-stout}).  We can now apply the theorem to the
family $f_\lambda$ restricted to $D$ to conclude that $f_\lambda$ is
linearizable and $p(\lambda)$ is in the interior of the domain of
linearization of $f_\lambda$.
\end{proof}
\begin{proof}[Proof of Lemma~\ref{l:local}]
  We use an argument of ``iteration in complex time'' in the spirit of
  \cite{ecalle}.  Pick a small disk $U$ around $\lambda_0$ and a small
  neighborhood $D$ containing $0$ in the dynamical plane such that the
  following hold:
\begin{itemize}
  \itm $f_\lambda\, p(\lambda) \neq p(\lambda)$ for all $\lambda \in
  U\setminus \set{\lambda_0}$; 

\itm the series $f^n_\lambda
  p(\lambda)= \sum_{j\ge 1} q_j^n \lambda^j$ converges uniformly in
  $U$ for all $n$; 

\itm the only solution to the equation $f_\lambda
  (z) = z$ for $z \in D$ is $z=0$.
\end{itemize}
The Cauchy estimates imply that $|q_j^n| \le C^j$ for some constant $C>0$ and
all $j,n$. As $q^n_j = Q_j (\mu^n)$, with  $Q_j$ a polynomial, by
continuity we obtain that $|Q_j(\zeta)|\le C^j$ for all $|\zeta|=1$. Then
  $|Q_j(\zeta)|\le C^j$ for  $|\zeta|\leq1$ by the Maximum Principle.
The
map $\phi_\lambda(\zeta) \mapsto \sum Q_j(\zeta) \lambda^j$ is thus
analytic and continuous in $\overline{\Delta}\times U$, and by construction
$f_\lambda \phi_\lambda (\zeta) = \phi_\lambda (\mu \zeta)$.

\smallskip

For $\lambda\in U$, the map $\phi_\lambda$ thus semiconjugates
$f_\lambda$ to $\zeta \mapsto \mu \zeta$ in the neighborhood of $0$,
but it needs not {\em a priori} be a conjugacy .  Nevertheless we claim that for
$\lambda\in U$, $f_\lambda$ is conjugate to $\mu \zeta$ in a disk
$D(0,r)$ of fixed size. The proof of Lemma \ref{l:interior} below will actually show that 
the conjugacy is global. 

Indeed, suppose that $\phi_\lambda (\zeta) = 0$ for some $\zeta \in
\fr \Delta$.
Then $\phi_\lambda (\mu^n \zeta)=0$ for all $n$, hence $\phi_\lambda
\equiv 0$ on $\fr \Delta$, hence on $\Delta$. In particular
$\phi_\lambda(1)= p(\lambda)= \phi_\lambda(\mu)=
f_\lambda p(\lambda)=0$ so by assumption $\lambda = \lambda_0$.  In
particular we obtain that for $\lambda\in \fr U$, $\phi_\lambda(S^1)$
remains at definite distance $r$ from the origin, and from this we
easily deduce that for $\lambda\in \fr U$, $\phi_\lambda(\Delta)\supset
D(0, r)$.

Thanks to the semiconjugacy, for every $\lambda\in \fr U$, points in
$D(0,r)$ never escape under iteration. By the Maximum Principle, the
same holds for $\lambda\in U$.  From this we deduce that
$|f_\lambda'(0)|\le 1$. The map $\lambda \to f_\lambda'(0)$ being
holomorphic, with values of modulus $1$ at $\lambda =\lambda_0$, we
conclude that it is constant equal to $\mu$.  Any indifferent point in
the Fatou set is linearizable, see~\cite[Theorem~II.6.2]{CG}. The
sequence $\{f^n_\lambda\}$ forms a normal family on $D(0,r)$ for each
$\lambda\in U$, so $D(0,r)$ is contained in the domain of
linearization. As $p(\lambda) \cv 0$ when $\lambda \cv \lambda_0$,
this implies that $p$ is in the interior of the domain of
linearization of $f_\lambda$ for small enough $\lambda$.

Finally note that by construction $f^n_\lambda \, p(\lambda) \neq
f^m_\lambda \, p(\lambda)$ for all $n\neq m$ and $\lambda\in
U\setminus{\lambda_0}$.
\end{proof}

\begin{proof}[Proof of Lemma~\ref{l:discrete}]
  Take a parameter $\lambda_1$ such that $ f^{n_1+k_1}_{\lambda_1}\,
  p(\lambda_1) = f^{k_1}_{\lambda_1}\, p(\lambda_1)$ for some
  intergers $n_1, k_1 \ge1$. Define $g_\lambda = f^{n_1}_\lambda$ and
  $q(\lambda) = f^{k_1}_\lambda\, p(\lambda)$. The point
  $q(\lambda_1)$ is now fixed by $g_{\lambda_1}$. Let $\mu_1$ be the
  multiplier of $g_{\lambda_1}$ at $q(\lambda_1)$. If $|\mu_1|<1$,
  then $q(\lambda)$ is attracted to an attracting fixed point for all
  parameters which is impossible as $0$ is indifferent for
  $f_{\lambda_0}$. If $|\mu_1|>1$ or if $\mu_1$ is a root of unity,
  Lemma~\ref{l:multi-first} implies that $q(\lambda)$ is periodic in a
  neighborhood of $\lambda_1$, whence $p(\lambda)$ is preperiodic for
  all parameters. Again this is impossible. We may thus apply
  Lemma~\ref{l:local} to $g$ and $q$. We conclude that in a punctured
  neighborhood of $\lambda_1$, the point $q(\lambda)$ is not
  $g$-preperiodic. This implies that $p$ is not $f$-preperiodic in
  this neighborhood.
\end{proof}

\begin{proof}[Proof of Lemma~\ref{l:interior}]
  Recall that in the proof of lemma \ref{l:local} we have constructed
  a natural holomorphically varying semiconjugacy
  $\phi_\lambda(\zeta)$ between $f_\lambda$ and $\zeta\mapsto
  \mu\zeta$. We now prove that it is a conjugacy. Indeed there exists a
  linearizing coordinate around zero in dynamical space, and for small
  $\lambda$, $p(\lambda)$ is inside the linearization domain. We
  denote by $\psi_\lambda$ the unique mapping such that
  $f_\lambda\rond\psi_\lambda= \psi_\lambda(\mu\,\cdot)$ and
  $\psi_\lambda'(0)=1$. Let also $s=\psi_\lambda^{-1}(z)$, so that
  $\widetilde{f_\lambda}(s)=\mu s$, where
  $\widetilde{f_\lambda}=\psi^{-1}_\lambda\rond
  f_\lambda\rond\psi_\lambda$. For small $\lambda$, in the linearizing
  coordinate,
  $\widetilde{f_\lambda}^n(\psi_\lambda^{-1}p(\lambda))=\mu^n
  \psi_\lambda^{-1}p(\lambda)$, so that by definition of
  $\phi_\lambda$ we have $\psi_\lambda^{-1}\rond \phi_\lambda(\zeta)=
  \zeta \, \psi_\lambda^{-1}p(\lambda)$. In particular, $\phi_\lambda$
  is a conjugacy for small $\lambda$, and
\begin{equation}\label{e:lambda}
\phi_\lambda(\zeta)=\psi_\lambda(\zeta\psi_\lambda^{-1}p(\lambda)) 
\end{equation}

Let us now consider the radius of convergence $R(\lambda)$ of the
power series defining $\zeta\mapsto \phi_\lambda(\zeta)$. We know that
this function is defined on the unit disk so $R(\lambda) \geq 1$.
Since $\phi_\lambda(1)=p(\lambda)$, to get the desired conclusion it
is enough to prove that $R>1$ everywhere. This will be a consequence
of a subharmonicity property of $R$.

For $\lambda$ close to zero, $\phi_\lambda$ is $p(\lambda)$ multiplied
by the linearizing coordinate so $R(\lambda)>1$ since $p(\lambda)$ is
inside the linearization domain. Write $\phi_\lambda(\zeta)=\sum_{k\ge
  0} a_k(\lambda) \zeta^k$.  The coefficient $a_k(\lambda)$ is defined
by the formula $a_k(\lambda) = \int_{S^1} \phi_\lambda (\zeta)
\zeta^{-k}$ so it depends holomorphically on $\lambda$. Since
$\phi_\lambda$ takes its values in the unit disk, $|a_k| \le 1$.  The
radius of convergence equals $R^{-1}(\lambda) = \limsup_k
|a_k(\lambda)|^{1/k}$. As $|a_k| \le 1$, the function $-\log R$ is the
supremum of a sequence of non-positive subharmonic function on
$\Lambda$. Its upper-semicontinuous regularization, $\rho^*$ thus
defines a non-positive subharmonic function. As $\lambda\cv0$,
$p(\lambda)\cv 0$ and $f_\lambda$ is linearizable on a fixed disk
$D(0,r)$, so $R(\lambda) \cv \infty$ (see \eqref{e:lambda}), hence
$\rho^* <0$ there.  By the Maximum Principle, we conclude that $\rho^*
<0$ everywhere. As $-\log R \le \rho^*$, we get that $R(\lambda) >1$
for all $\lambda \in \Lambda$.
\end{proof}

\begin{proof}[Proof of Lemma~\ref{l:comp}]
  Write $f^n_\lambda(z)= \mu^n z + \sum_{k+l\ge 2} a_{kl}^n \lambda^k
  z^l$.  Then $f^{n+1} = f^n \circ f$ from which we deduce that 
  $$
  f^{n+1}_\lambda(z)= \mu^n \left( \mu z+ \sum_{k+l\ge2} a_{kl}
    \lambda^k z^l\right) + \sum_{i+j\ge 2} a^n_{ij}\lambda^i \left( \mu z+
    \sum_{p+q\ge2} a_{pq} \lambda^p z^q\right)^j~.
$$
For $k+l\ge2$, we thus infer that 
$$
a^{n+1}_{kl} = \mu^n a_{kl} + \sum_{i+j\ge 2} a^n_{ij} \times \text{ Term
  in } \lambda^{k-i} z^l \text{ of } \left[  \left(\mu z+
    \sum_{p+q\ge2} a_{pq} \lambda^p z^q\right)^j\right]
$$
The sum over $(i,j)$ in the right hand side is finite, as we
necessarily have $i \le k$ and $j\le l$. Note also that for $(i,j) =
(k,l)$, the contribution of the sum is exactly $\mu^l a^n_{kl}$.  We
conclude that there exist constants $c, c_{ij}$, independent of $n$,
such that
$$
a^{n+1}_{kl} = \mu^l a^n_{kl} + c \cdot \mu^n + \sum_{(i,j)<(k,l),\,
  i+j\ge 2}
c_{ij} a^n_{ij}~,
$$
where $\le$ (resp. $<$) here denotes the partial order on
$\mathbb{N}^2$ given by $(i,j) \le (k,l)$ iff $i\le k$ and $j\le l$
(resp. at least one inequality is strict).

We complete the proof by an induction on $(k,l)$ with
respect to the previously defined partial order. Indeed if for every
$(i,j)<(k,l)$, $a_{ij}^n$ is of the form $\sum_{r=0}^j
\mu^{rn}P_{rij}(n)$ we get that
$$a^{n+1}_{kl} = \mu^l a^n_{kl} + \sum_{r=0}^l \mu^{rn}Q_{rkl}(n)~.$$
Since for every complex number $\nu$, the sum $\sum_{k=0}^n \nu^k k^d$ is of the
form $\nu^n P(n)$ where $P$ is a polynomial, 
the sequence $a^n_{kl}$ is of the required
form~\eqref{e:poly}.
\end{proof}

\section{Generalities on families of rational maps} \label{sec:gen}


\subsection{Active and passive points}\label{subs:active}
We start with a definition.
\begin{defi}\label{def:family}
  A holomorphic family of rational maps of degree $d$ with a marked critical
  point is a triple $(\Lambda, f, c)$ such that:
\begin{itemize}
\itm $\Lambda$ is a smooth connected complex manifold;
\itm $f: \Lambda \times \PP^1 \to \PP^1$ is a holomorphic map;
\itm for any $\lambda \in \Lambda$, the map $z \mapsto f_\lambda(z)
  \= f (\lambda, z)$ is a rational map of $\PP^1$ of degree $d$;
\itm $c: \Lambda \to \PP^1$ is a holomorphic map such that
  $f'_\lambda c(\lambda) =0$ for all $\lambda$.
\end{itemize}
We call $\Lambda$ the parameter space of the family. 
\end{defi}
Throughout the paper we will focus our interest
 on the locus of points in parameter space
for which $c(\lambda)$ is dynamically unstable. This is a classical
 notion, which goes back to Levin \cite{leactive} and Lyubich
 \cite{lyubactive}.  We use the terminology of McMullen \cite{McM1}.  
\begin{defi}\label{def:active}
  The marked critical point $c$ is \emph{passive} at $\lambda_0\in\Lambda$ if
  $\set{\lambda\mapsto f_\lambda^n c(\lambda)}_{n\in\nn}$ forms a
  normal family of holomorphic functions in the neighborhood of
  $\lambda_0$. Otherwise $c$ is said to be \emph{active} at
  $\lambda_0$.  The set of points $\cA \subset \Lambda$ where $c$ is
  active is called the \emph{activity locus}.
\end{defi}
Notice that the family $(f_\lambda)$ being structurally stable on its Julia
set is actually equivalent to the fact that all marked critical
points are passive, see e.g. \cite{McM2}. In
particular the results of \cite{MSS} imply that the complement of the
activity locus is always an open dense set in $\Lambda$.
\begin{lem}\label{lem:active_prerep} Let $(\Lambda, f, c)$ be a
  holomorphic family of rational maps with a marked critical point.
If $c$ is active at $\lambda_0$, then there exists a nearby
  parameter $\lambda$ such that $c(\lambda)$ is prerepelling
  (i.e. preperiodic to a repelling cycle). 
\end{lem}
\begin{proof} 
  The lemma follows from a classical normal family
  argument.  Fix three repelling periodic points at the parameter
  $\lambda_0$.  They persist in a neighborhood of $\lambda_0$. By
  conjugating with a holomorphically varying M{\"o}bius transformation, we
  may further assume they persistently  equal  $\{ 0, 1 ,
  \infty\}$.  Since the family
  $\set{f_\lambda^nc(\lambda)}_{n\in\nn}$ is not normal in any
  neighborhood of $\lambda_0$, it cannot avoid these three points.
\end{proof}
The following  result is a little bit  more delicate (see also \cite{leactive,McM1}).
\begin{prop}\label{prop:mcmullen}
Suppose $c$ is active at $\lambda_0$. In any neighborhood of
$\lambda_0$, there exists a parameter $\lambda$ for which $f^n_\lambda
c(\lambda)$ converges to an \emph{attracting cycle} when $n\to\infty$.
\end{prop}
\begin{proof} We will actually prove that there exists a parameter close to $\lambda_0$ where 
$c$ is periodic, hence superattracting.
Replacing $\Lambda$ by a suitable finite ramified cover, we may follow holomorphically a
preimage of $c(\lambda)$. Observe that this operation preserves active and passive critical points..
We thus get a holomorphic map
$c_{-1}(\lambda)$ such that $ f_{\lambda}\,(c_{-1}( \lambda))=
c(\lambda)$.  Similarly we
get a holomorphic map $c_{-2}$ with
$f_\lambda c_{-2}(\lambda)=c_{-1}(\lambda)$. There are two cases: either $c(\lambda)$,    
$c_{-1}(\lambda)$ and $c_{-2}(\lambda)$ are disjoint in the neighborhood of $\lambda_0$ or there is a parameter 
near $\lambda_0$ where $c$ is fixed. In the latter case we are done. In the former, 
we can  choose a holomorphically varying family of M{\"o}bius maps
$\phi_\lambda$ such that $\phi_\lambda\, c(\lambda) =0$,
$\phi_\lambda\, c_{-1}(\lambda) =1$ and $\phi_\lambda\,
c_{-2}(\lambda) = \infty$. Replacing $f_\lambda$ by the family
$\phi_\lambda \circ f_\lambda \circ \phi_\lambda^{-1}$ (still denoted by $f_\lambda$), 
we get that $\infty \mapsto 1 \mapsto 0 \equiv c(\lambda)$. 
Now  if $c$ is active at $\lambda_0$, the family $(f^n_\lambda c(\lambda))$ 
cannot avoid $\set{0,1,\infty}$ and we conclude that there are parameters close to $\lambda_0$ 
where $c$ becomes periodic. 
\end{proof}
We now prove Theorem~\ref{p:motion} stated in the introduction, which 
classifies the dynamics  of passive critical points. This is essentially a 
 consequence of Theorem~\ref{t:local}. 
\begin{proof}[Proof of Theorem~\ref{p:motion}] 
  Suppose $c$ is not stably preperiodic, and
  $f^{n_0+k_0}_{\lambda_0}\,c(\lambda_0)=
  f^{n_0}_{\lambda_0}\,c(\lambda_0)$ for some $\lambda_0\in U$, and
  minimal $n_0\ge0$, $k_0\ge1$.  By Theorem~\ref{t:local}, there
  exists a periodic cycle $\set{p(\lambda), \ldots,
    f^{k_0-1}(p(\lambda))}$ which is either attracting or with
  multiplier $\exp(2i\pi\theta)$ with $\theta$ irrational such that
  $p(\lambda_0) = f^{n_0}_{\lambda_0}\,c(\lambda_0)$; and either $c$
  is attracted towards the orbit of $p$ for all parameters, or $c$
  eventually lies in the interior of the domain of linearization of
  $p$ (hence $f$ has a Siegel disk at $p$).

\smallskip
  
  Suppose we are in the first case. Pick a small neighborhood $U$ of
  $\lambda_0$, and small disjoint disks $U_i$ such that 
\begin{itemize}
\itm
$U_0=U_{k_0}
  \ni p(\lambda), U_1\ni f_\lambda\,p(\lambda), \cdots , U_{k_0-1} \ni
  f^{k_0-1}_\lambda \, p(\lambda)$ for all $\lambda \in U$; 
\itm
$f_\lambda$ is injective on $U_i$ for all $i$ and all $\lambda$; 
\itm $f_\lambda U_i $ is relatively compact in $U_{i+1}$ for all $i$
  and all $\lambda$.
\end{itemize}
Replacing $U$ by an even smaller open set containing $\lambda_0$, we
may assume that $f^n_\lambda\, c(\lambda) \in \cup_0^{k_0-1}U_i$ for
all $\lambda \in U$ and all $n$. It is then clear that $c(\lambda)$ is
preperiodic iff $f^n_\lambda\, c(\lambda) = p(\lambda)$ for some $n$,
and that it is equivalent to $f^{n_0}_\lambda\, c(\lambda) =
p(\lambda)$.  The latter condition defines a closed analytic subspace
in $U$.

\medskip

When $p$ is neutral, the proof goes through with the following
additional remarks. Replace the last condition in the construction of
$U_i$ above by taking $U_i$ to be an open set containing the closure
of the orbit of $f^{n_0}_\lambda\,c(\lambda)$ for all parameters.
This is possible even though the Siegel disk need
  not move continuously with $\lambda$.
Indeed, by the proof of lemma \ref{l:local}, for
$\lambda$ close to $\lambda_0$, there is a disk of fixed size
$D(p(\lambda),r)$ contained in the linearization domain.

As before, we
conclude that $c(\lambda)$ is preperiodic iff
$f^{n_0+k_0}(c(\lambda))= f^{n_0}(c(\lambda))$.
\end{proof}
\subsection{Quasiprojective parameter space} \label{subs:quasiproj}

Following McMullen, call \emph{isotrivial} a family in which any two
members are conjugate by a M{\"o}bius transformation; and
\emph{algebraic} a family parameterized by a quasiprojective
variety.  

McMullen proved in~\cite{McM2} that any non isotrivial
algebraic family of rational maps admits bifurcations, with the only
exception of families of {\em flexible Latt{\`e}s examples}.  A
consequence of this result is that any non isotrivial algebraic family
is critically finite or has bifurcations. Here, building
on~\cite{McM2}, we prove the more precise result that in a non
isotrivial algebraic family of rational maps {\em each} critical point
either presents bifurcations or is stably preperiodic.

\begin{thm}\label{t:quasiproj}
Let $(\Lambda, f, c)$ be a holomorphic family of rational maps of
degree $d\geq 2$ with a marked critical point, where $\Lambda$ is an 
irreducible quasiprojective complex variety. Assume that the
activity locus of $c$ is empty. Then: 
\begin{itemize}
\itm either all maps $f_\lambda$ are conjugate to each other by M{\"o}bius
transformations;
\itm or there exist integers $(n,k)$ such that for every $\lambda\in
\Lambda$, $f_\lambda^n(c_\lambda)=f_\lambda^k(c_\lambda)$.
\end{itemize}
\end{thm}  

\begin{proof} Taking successive hyperplane sections, we reduce the
  proof to the case where $\Lambda$ is a Riemann surface of finite
  type, that is a compact Riemann surface with finitely many points
  deleted. The crucial fact is that there are only finitely many non
  constant holomorphic maps $\Lambda \cv
  \PP^1\setminus\set{0,1,\infty}$ (see \cite[p.478]{McM2} for a
  proof).

We start with an easy lemma.
\begin{lem}\label{l:constant}
  Let $(\Lambda, f, c)$ be as above, and assume that $c$ is not stably
  preperiodic.  Assume that for large $n$, the function $f^n_\lambda
  \,c(\lambda)$ is constant and does not depend on $\lambda$. Then the
  family is constant, that is, for all $\lambda, \lambda'$, we have
  $f_\lambda= f_{\lambda'}$.
\end{lem}

\begin{proof} For some integer $n_0$ and all $n\geq n_0$, 
  $u_n=f^n_\lambda\, c(\lambda)$ does not depend on $\lambda$.
  Moreover, since there exists a parameter for which $c$ has infinite
  orbit, the set $\set{u_n, n\geq n_0}$ is infinite. Take $\lambda
  \neq\lambda'$. Then $f_\lambda(u_n)= f_{\lambda'}(u_n)$ for all
  $n\ge n_0$. Any two rational maps of constant degree, agreeing on an
  infinite set are equal. Hence $f_\lambda= f_{\lambda'}$.
\end{proof}

We continue with the proof of the theorem. We thus assume that $c$ is passive
throughout $\Lambda$, and  not stably preperiodic.
We shall prove that any two rational maps in this family are
conjugated by a M{\"o}bius transformation.

\medskip

Suppose first that $c$ is \emph{never} periodic. As in Proposition~\ref{prop:mcmullen}, 
replacing $\Lambda$
by a suitable finite ramified cover, we may follow holomorphically two 
preimages of $c(\lambda)$.  We thus get  holomorphic maps
$c_{-1}(\lambda)$ and $c_{-2}(\lambda)$ such that $ f_{\lambda}\,(c_{-1}( \lambda))=
c(\lambda)$ and $f_\lambda(c_{-2}(\lambda))=c_{-1}(\lambda)$.   Note that any finite
ramified cover of a Riemann surface of finite type is of finite type
too. Also by assumption, $c$ is never periodic so $c(\lambda)$,
$c_{-1}(\lambda)$ and $c_{-2}(\lambda)$ are always distinct. 
As before we can thus choose a holomorphically varying family of M{\"o}bius maps
$\phi_\lambda$ such that $\phi_\lambda\, c(\lambda) =0$,
$\phi_\lambda\, c_{-1}(\lambda) =1$ and $\phi_\lambda\,
c_{-2}(\lambda) = \infty$ and we replace $f_\lambda$ by the family
$\phi_\lambda \circ f_\lambda \circ \phi_\lambda^{-1}$.

As $c$ is never periodic, for any $n>0$ and any $\lambda\in \Lambda$,
$f^n_\lambda(c(\lambda))$ avoids $\{ 0, 1, \infty \}$. So the family of
holomorphic maps $f^n_\lambda\,c(\lambda):
\Lambda\cv\PP^1\setminus\set{0,1,\infty}$ takes only finitely many
nonconstant terms. Since we have moreover assumed that $f^n_\lambda(c(\lambda))$ 
is not persistently preperiodic, $f^n_\lambda(c(\lambda))$  is constant for large $n$. 
By Lemma~\ref{l:constant}, we conclude that the family is constant.

\medskip

Suppose then that $c$ is \emph{periodic} for some parameter.  We may
assume that for some $\lambda_0$,
$f^{k_0}_{\lambda_0}\,c(\lambda_0)=c(\lambda_0)$, with minimal $k_0$.
Using Theorem~\ref{p:motion}, we infer that $c$ is attracted to an
attracting periodic orbit throughout $\Lambda$. The multiplier of this
periodic orbit defines a holomorphic function on a Riemann surface of
finite type with values in the unit disk. It is hence constant equal
to $0$.

We conclude that there is a persistent superattracting cycle
$\set{p(\lambda), \ldots, f^{k_0}\,p(\lambda)}$ attracting $c$, and
$c(\lambda_0)=p(\lambda_0)$. The cycle can be followed holomorphically
because the multiplier never equals 1.  Replacing the family
$\{f_\lambda\}$ by $\{f^{k_0}_\lambda\}$, we may assume that
$p(\lambda)$ is a fixed point.  Let $m-1$ be the largest integer $q$
such that $\frac{d^qf}{d\lambda^q}\,p(\lambda)$ is identically zero.
Thanks to the B{\"o}ttcher theorem, outside some  discrete subset
$\mathcal{D}\subset \Lambda$, $f_\lambda (\zeta) = \zeta^m$ in a
suitable coordinate $\zeta$ centered at $p(\lambda)$. In this case, we
can even choose $\zeta$ to depend holomorphically on $\lambda$.  On
the complement $\mathcal{D}$, the multiplicity is larger so 
$f_\lambda (\zeta) = \zeta^{m_\lambda}$
for some integer $m_\lambda > m$.

Let $d_{\PP^1}$ be the spherical distance on $\PP^1$ ($d_{\PP^1}\leq
1$). If $z$ lies in the basin of attraction of $p(\lambda)$, let
$$g_\lambda(z)= \lim_{n\cv\infty} \unsur{m^n} \log d_{\PP^1}
(f_\lambda^{n}(z), p(\lambda)).$$
Notice that the limit is the same
when one replaces $d_{\PP^1}$ by any equivalent distance.  In the
local chart $\zeta$ mentioned above, we thus get $g_\lambda(\zeta) =
\log |\zeta|$ for any $\lambda \notin \mathcal{D}$.  It follows that
$(z,\lambda)\mapsto g_\lambda(z)$ is plurisubharmonic and continuous
in the domain consisting of couples $(z,\lambda)$ such that $\lambda
\notin\mathcal{D}$ and $z$ in the basin of $p$.

In particular the function $\lambda\mapsto g_\lambda\,c(\lambda)$ is
subharmonic, non positive, and continuous outside the discrete subset
of $\mathcal{D}\subset\Lambda$. It may thus be extended across
$\mathcal{D}$ as a non positive plurisubharmonic function on
$\Lambda$.  Now $\Lambda$ is a Riemann surface of finite type, hence
$g_\lambda\,c(\lambda)$ extends across the punctures of $\Lambda$ and
defines a non positive subharmonic function a \emph{compact} Riemann
surface.  It is thus constant.  Since
$g_{\lambda_0}\,c(\lambda_0)=-\infty$, we get that
$g_\lambda\,c(\lambda)\equiv -\infty$, that is $c(\lambda)\equiv
p(\lambda)$ is persistently periodic.

This concludes the proof of the theorem.
\end{proof}


\subsection{Exceptional points} \label{subs:exc}

To get a precise statement of our convergence theorem, we need a short
discussion on exceptional points. Let $f : \PP^1 \to \PP^1$ be any
rational map of degree $d\ge2$. A point $z$ in the Riemann sphere is
called \emph{exceptional}, if $z$ is totally invariant by $f$ or by
$f^2$.  The set of exceptional points of $f$ is denoted by ${\mathcal
  E}(f)$, and  contains  at most two points. When
$\mathrm{card}({\mathcal E}(f)) =1$, then $f$ is conjugate to a
polynomial for which ${\mathcal E}(f) = \{ \infty \}$.  When
$\mathrm{card}({\mathcal E}(f)) =2$, then $f$ is conjugate to $z
\mapsto z^d$ for some integer $d \in {\mathbb Z} \setminus\{\pm 1,
0\}$.

Let now $(\Lambda, f)$ be a holomorphic family of rational maps.
Define $\mathcal{E}\= \{ (\lambda, z),\, z \in
\mathcal{E}(f_\lambda)\}$. This defines a closed analytic subset in
$\Lambda \times \PP^1$. It can be empty, or reducible, or singular.
For the sake of simplicity, we write $e$ for the
cardinality of $\mathcal{E}(f_\lambda)$ for  \emph{generic}
$\lambda$. As $\Lambda$ is connected, this is precisely $\min_\Lambda
\mathrm{card}({\mathcal E}(f_\lambda))$.

\smallskip

\noindent \underline{$e=2$.}
Then $\mathrm{card}({\mathcal E}(f_\lambda)) = 2$ for \emph{all}
$\lambda$. In this case, $\pi_1$ induces a non-ramified $2$ to $1$
cover of $\mathcal{E}$ onto $\Lambda$, and $\mathcal{E}$ is a smooth
hypersurface. Locally at any parameter $\lambda_0$, the family is
conjugate to the trivial family $f_\lambda(z) = z^d$,
$d\in\mathbb{Z}\setminus \{ \pm 1, 0\}$.  In general, the family need
not be globally trivial (e.g. $f_\lambda(z) = \lambda z^d$
over $\mathbb{C}^*$).

\smallskip

\noindent \underline{$e=1$.}
At a generic point, $\mathrm{card}({\mathcal E}(f_\lambda)) = 1$. Let
$\mathcal{E}'$ be the irreducible component of $\mathcal{E}$
containing ${\mathcal E}(f_\lambda)$ for generic $\lambda$. Then
$\pi_1:\mathcal{E}'\to\Lambda$ is an isomorphism, and after passing to
the universal cover of $\Lambda$, we may assume that
$\pi_2(\mathcal{E}') = \infty \in \PP^1(\mathbb{C})$. In particular,
$f_\lambda$ is a \emph{polynomial for all $\lambda$}.

\smallskip

\noindent \underline{$e=0$.}
This is the generic case. Outside a  proper closed Zariski subset
of $\Lambda$, the set $\mathcal{E}(f_\lambda)$ is \emph{empty}.


\section{The bifurcation current}\label{sec:bif} 

We now associate a natural positive closed $(1,1)$ current to the data
$(\Lambda, f,c)$. This will be the main object of interest in
the paper.  We refer the reader to Demailly \cite{de} for basics on
positive closed currents.

We first fix some notation.  Write $\hL \= \Lambda \times
\PP^1$.  The family of maps $f_\lambda$ lifts to a holomorphic map
$\hf : \hL \to \hL$ sending $(\lambda, z) $ to $(\lambda,
f_\lambda(z))$.  We denote by $\pi_1: \hL \to \Lambda$ and $\pi_2: \hL
\to \PP^1$ the two natural projections. The map
$\hat{p} : \Lambda \to \hL$ defined by $\hat{p}(\lambda) = (\lambda,
p(\lambda))$ induces an isomorphism from $\Lambda$ onto its image
$\Gamma$ which is a smooth submanifold of $\hL$. It defines a section
of $\pi_1$, that is, $\pi_1\rond \hat{p}$ is the identity map.

Let now $\omega$ be the unique smooth positive closed $(1,1)$ form on
$\PP^1$ which is invariant under the unitary group, and normalized by
$\int_{\PP^1} \omega =1$ (the Fubini-Study form).  It induces the
spherical distance on $\PP^1$, which we denote by $d_{\PP^1}$. We set
$\ho \= \pi_2^*\omega$.
\begin{prop-def}\label{def:bifurcation}
  For any holomorphic family $(\Lambda, f, c)$ of rational maps of
  degree $d\ge2$ with a marked critical point, the sequence of
  positive closed $(1,1)$ currents $d^{-n}\hf^{n*} \ho$ converges to a
  positive closed $(1,1)$ current $\hT$ on $\hL$.
  
  This current admits a continuous potential locally at any point.
  The image under $\pi_1$ of the restriction of $\hT$ to the
  submanifold $\Gamma \= \{ (\lambda, c(\lambda)) \in \hL\}$ is thus
  well-defined as a positive closed $(1,1)$ current on $\Lambda$. It
  is called the \emph{bifurcation current} of the family, and we
  denote it by $T$.
\end{prop-def}
\begin{proof}
  All assertions are \emph{local} in $\Lambda$. We may thus assume
  that $\Lambda$ is the unit open polydisk in $\cc^n$, and that the
  family is defined over a neighborhood of the closed unit polydisk.
  
  We first claim the existence of a continuous and bounded function
  $\hg$ on $\hL$ such that 
\begin{equation}\label{e1}
d^{-1} \hf^* \ho - \ho = dd^c \hg~.
\end{equation}
To see this, choose homogeneous coordinates $[z:w]$ on $\PP^1$, and
let $\pi : \mathbb{C}^2\setminus \{ 0 \} \to \PP^1$ be the natural projection.
Take a lift of the family of rational maps $f$ to $\cc^2$. Write
$f_\lambda [z:w] = [ P_\lambda (z,w): Q_\lambda(z,w)]$ where
$P_\lambda$, $Q_\lambda$ are two holomorphic families of homogeneous
polynomials of degree $d$ with no common factors (for each $\lambda$).
This is always possible by restricting $\Lambda$ if necessary. The function
$$
\hg (\lambda, z) \= \frac1{2d} \log\frac{|P_\lambda(z,w)|^2+
  |Q_\lambda(z,w)|^2}{(|z|^2+|w|^2)^d}
$$
clearly satisfies  our requirements.

Applying $d^{-k} f^{k*}$ to~\eqref{e1}, and summing from $k=1$ to
$n-1$, we get that $d^{-n} \hf^{n*} \ho - \ho = dd^c \sum_0^{n-1}
d^{-k}\, \hg\circ\hf^k$. As $\sup_{\hL} |d^{-k} \, \hg \circ \hf^k |
\le d^{-k} \sup_{\hL}|\hg|$, the sequence of functions $\sum_0^{n-1}
d^{-k}\, \hg\circ\hf^k$ converges uniformly on $\hL$ to a continuous
function $\hg_\infty$. In particular, $d^{-n} \hf^{n*} \ho \cv \hT \=
\ho + dd^c \hg_\infty$.

The restriction to a subvariety of a continuous $(1,1)$ form and of a
continuous function is always well-defined. We may thus set
$\hT|_\Gamma \= \ho|_\Gamma + dd^c (\hg_\infty|_\Gamma)$. Finally,
$\pi_1$ is an isomorphism from $\Gamma$ onto $\Lambda$, so we may 
define $(\pi_1)_* ( \hT|_\Gamma)$ as the bifurcation current. 
\end{proof}
\begin{thm}\label{thm:support}
The support of the bifurcation current coincides with the activity locus.
\end{thm}
From this theorem and the convergence statement in Definition
\ref{def:bifurcation}, we get the following interesting
corollary.
\begin{cor}\label{cor:normal-sub}
  Assume that there exists an increasing sequence of integers
  $(k_n)$ such that $\set{f_\lambda^{k_n} c(\lambda)}$ is a normal
  family in some neighborhood of $\lambda_0$. Then $c$ is passive at
  $\lambda_0$. 
\end{cor}
\begin{proof}[Proof of Theorem~\ref{thm:support}]
  Suppose first that $c$ is passive at $\lambda_0$. We may thus assume
  that $f_\lambda^{k(n)} c(\lambda) \cv h(\lambda)$ in a neighborhood
  $U$ of $\lambda_0$ for some sequence $k(n)$ increasing to infinity.
  Notice that the map $\lambda \to f^n_\lambda c(\lambda)$ is
  identical to the composite map $\pi_1\circ \hf^n \circ \hat{p}$.  On
  $U$, we infer that $(\pi_1\circ \hf^{k(n)} \circ \hat{p})^*\, \omega
  - h^* \omega$ tends to zero in the supremum norm of forms (hence in
  the weak topology of currents).  By construction, we have $T|_U =
  \lim_{n\to\infty} d^{-n} (\pi_1\circ \hf^{n} \circ \hat{p})^*
  \omega$, whence $T|_U \equiv 0$.

\medskip

Conversely, assume that $T|_U \equiv 0$ in a neighborhood $U$ of
$\lambda_0$. We need to show that the sequence $\{ f^n_\lambda
c(\lambda) \}$ is normal in a neighborhood of
$\lambda_0$.  We may assume without loss of generality that $\dim
(\Lambda) =1$.

On $U$, we have 
\begin{equation*}
  (\pi_1 \circ \hf^n \circ \hat{p}) ^* \omega = d^n \left( d^{-n}
    (\pi_1 \circ \hf^n \circ \hat{p}) ^* \omega - T \right)= dd^c
  \left[ d^n \, \sum_n^\infty \frac{\hg\circ \hf^k\circ \hat{p}}{d^k}\right]~,
\end{equation*}
where $dd^c \hg = d^{-1} \hf^* \ho - \ho$ as in the previous proof.
The function $\hg$ being continuous, the sequence $ d^n \,
\sum_n^\infty d^{-k} \hg\circ \hf^k\circ \hat{p}$ is uniformly bounded
in $n$.  In particular, $(\pi_1 \circ \hf^n \circ \hat{p}) ^* \omega =
dd^c \phi_n$ for a sequence of continuous functions $\phi_n$ with
$|\phi_n|\le C$. Note that $dd^c\phi_n$ is positive, thus $\phi_n$
is \emph{psh on $U$}.

Let $\omega_\Lambda$ be any smooth positive $(1,1)$ form on $U$.  As
$\dim(\Lambda) =1$, we have
\begin{equation*}
(\pi_1 \circ \hf^n \circ \hat{p})^* \, \omega 
 = \|d \left( f^n_\lambda c(\lambda)\right)
\|^2\,  \omega_\Lambda~,
\end{equation*}
where $\|d \left( f^n_\lambda c(\lambda)\right) \|$ is the norm of the
differential of the map $\lambda \mapsto f^n_\lambda c(\lambda)$
computed in the metrics induced by $\omega_\Lambda$ in $\Lambda$ and
$\omega$ in $\PP^1$. Fix any relatively compact open set
$V\subset\subset U$. The relative capacity of the compact set
$\overline{V}$ with respect to $U$ is by definition the non-negative
real number $\mathrm{cap}_U(\overline{V})\= \sup \{ \int_{\overline{V}}
dd^c u,\, 0\le u\le 1, u \text{ psh}\}$, see~\cite{BT}.  In particular,
we obtain
\begin{equation*}
\int_V \|d \left( f^n_\lambda c(\lambda)\right)
\|^2\,  \omega_\Lambda 
\le \int_{\overline{V}}  dd^c (\phi_n+ C) \le 2\, C\, 
\mathrm{cap}_U(\overline{V})~. 
\end{equation*}
When $V$ decreases to a point, its capacity tends to zero.  In
particular, by choosing $V$ sufficiently small, we may assume that
$\sup_n \int_V \|d \left( f^n_\lambda c(\lambda)\right)\|^2
\omega_\Lambda < \e_0 <1$.

Now look at the sequence of (closed) analytic curves $\Gamma_n = \{
(\lambda, \hf^n c(\lambda))\} \subset V \times \PP^1$. The area of
$\Gamma_n$ with respect to the metric $\pi_1^* \omega_\Lambda +
\pi_2^* \omega$ is precisely $\int_V (1+ \|d \left( f^n_\lambda
  c(\lambda)\right) \|^2)\, \omega_\Lambda$. Reducing $V$ again if
necessary, we may assume that $\sup_n \mathrm{Area}(\Gamma_n) < \e_0
<1$.  By Bishop's theorem, from any subsequence, one can extract a
(sub-)subsequence converging in the Hausdorff topology to an analytic
curve $\Gamma_\infty$. The area of $\Gamma_\infty$ is also less than
$\e_0<1$, thus $\Gamma_\infty$ cannot contain any fiber of $\pi_2$
(the area of any such fiber is $\mathrm{Area}_\omega(\PP^1)=1$).  All
curves $\Gamma_n$ being graphs over $\Lambda$, we conclude that the
limit $\Gamma_\infty$ itself is a graph over $\Lambda$.

We have thus shown that on a fixed small neighborhood $V$ of
$\lambda_0$, any subsequence of $\{ f^n_\lambda c(\lambda)\}$ admits
a converging subsequence, since the associated sequence of graphs
does. We conclude that $\{ f^n_\lambda
c(\lambda)\}$ is a normal family on $V$. Hence $c$ is passive at
$\lambda_0$.
\end{proof}

\section{The convergence theorem} \label{subs:cv}
We now come to the proof of Theorem~\ref{t:cvg}.
\begin{defi}\label{def:per}
  Let $(\Lambda, f, c)$ be a holomorphic family of rational maps of
  $\PP^1$ with one marked critical point. For any $n\neq m\ge0$, we
  denote by $\Per(n,m)$ the (not necessarily reduced) analytic subset
  of $\Lambda$ defined by the equation $f^n_\lambda c(\lambda)
  =f^m_\lambda c(\lambda)$. For the sake of simplicity, we write
  $\Per(n)$ for $\Per(n,0)$. As a set, it consists of the parameters for
  which $c(\lambda)$ is periodic of period dividing $n$.

  A family is called \emph{trivial} when $\Per(n,m) = \Lambda$ for
  some $n\neq m$.
\end{defi}
We stress that we consider the analytic sets $\Per(n,m)$ as defined by 
the equations $f^n_\lambda c(\lambda)  =f^m_\lambda c(\lambda)$, so they can come 
along with multiplicities.
Note that the activity locus is empty in the case where $\Lambda$ is a
trivial family. Note also that when the cardinality of the exceptional
set of a generic map is $2$, then the family is trivial.

In the sequel, we shall work exclusively with non-trivial families of
rational maps of $\PP^1$. For any such family and any $n\neq m$,
$\Per(n,m)$ is a hypersurface in $\Lambda$, which might be singular,
or reducible, or even non-reduced. For convenience, we state again the
convergence theorem.
\begin{thm}\label{thm:cvg}
  Let $(\Lambda, f, c)$ be a non-trivial holomorphic family of
  rational maps of $\PP^1$ of degree $d\ge 2$ with a marked critical
  point.  We moreover assume one of the following two assumptions holds.
\begin{itemize}
\item[(H1)]{\it The manifold $\Lambda$ is quasi-projective.}
\item[(H2)]{\it For any $\lambda \in \Lambda$, there exists an
    immersed analytic curve $\Gamma \subset \Lambda$ containing
    $\lambda$ such that the closure of the complement of $\{ \lambda,
    \, c(\lambda) \text{ is attracted to a periodic cycle}\}$ is a
    proper compact subset of $\Gamma$.}
\end{itemize}
Let also $e\in\{0,1\}$ be the cardinality of the exceptional set of a
generic $f_\lambda$.  Then for any sequence of integers $0\le k(n) <
n$, the following convergence statement holds:
 $$\lim_{n\to\infty} \frac{[\Per(n,k(n))]}{d^n+d^{(1-e)k(n)}}  = T$$
\end{thm}
Notice that in the assumption (H2),
 we allow the curve $\Gamma$ to be singular, or 
  non properly embedded in $\Lambda$. Recall also that we have 
classified in Theorem \ref{t:quasiproj} the algebraic families of 
rational maps for which 
$T\equiv 0$. We refer to the next
sections for applications of this result to concrete families.

\smallskip

An important fact here is that the subvarieties $\Per(n)$,
$\Per(n,k)$, have several irreducible components: e.g. if $d\vert n$,
$\Per(d) \subset \Per(n)$. Also periodic points are preperiodic so
$\Per({n-k})\subset \Per(n,k)$. We  a priori have no control on the
multiplicities of the various irreducible components of $\Per(n,k)$. 
It could be possible that the multiplicity of $\Per({n-k})$ as a
component of $\Per(n,k)$ is so large  that most of the mass of
$[\Per(n,k)]$ would actually be concentrated on $\Per({n-k})$.

Let $\preper({n,k})\subset \Per(n,k)$ be the closure of the subset of
parameters where $c(\lambda)$ is {\it strictly} preperiodic.  More
precisely, $\preper({n,k})$ is locally defined in the open Zariski
subset $\bigcup_0^n\{ f^j_\lambda \, c(\lambda) \neq c(\lambda)\}$ by $f^n
_\lambda\, c(\lambda) = f^k_\lambda \, c(\lambda)$.  It extends
uniquely to an analytic subset of $\Lambda$.  It is a hypersurface,
which consists of a union of irreducible components of $\Per(n,k)$,
and endowed with the same multiplicities as $\Per(n,k)$.

Because we allow arbitrary sequences $0\leq k(n)\leq n$,
we can strengthen the previous result as follows. 
\begin{cor}\label{cor:strpreper}
  Assume that $(\Lambda, f, c)$ is a family of rational maps satisfying (H1) or 
  (H2)  as in the preceding theorem. Let $e\in\{0,1\}$ be the generic
  cardinality of $\mathcal{E}(f_\lambda)$.  Then for every fixed
  $k\geq 1$, we have
$$\unsur{d^n+d^{(1-e)(n-k)}}[\preper({n,n-k})] \cv T~.$$
\end{cor}
\begin{proof} We assume $e=0$, the other case is similar.
  If $\lambda\in\Per({n,n-k})$ and $c(\lambda)$ is periodic, then its
  period divides $k$. So we may write
  $$[\Per({n,n-k})] = [D_n]+[\preper({n,n-k})],$$ where $D_n$ is a divisor supported on 
$\Per(k)$. We only
  need to prove that $[D_n]/(d^n+d^{n-k})\cv 0$. But if not, $T$ would
  give some mass to the subvariety $\Per(k)$. This is impossible as
  $T$ has continuous potential.
\end{proof}

\begin{proof}[Proof of Theorem~\ref{thm:cvg} in the case where $e =0$.]
  Let $F_n:\Lambda\cv \PP^1$ be the map defined by
  $F_n(\lambda)=f_\lambda^n(c(\lambda))$. Recall from Proposition
  \ref{def:bifurcation} that the bifurcation current is defined as the
  limit of the sequence $F_n^*\omega$, where $\omega$ is the
  Fubini-Study (1,1) form on $\PP^1$. We first prove that the
  potentials of the currents $ [\Per(n, k(n))] - (F_n^*\omega +
  F_{k(n)}^{*}\omega)$ converge to zero on the activity locus.  A
  classical argument based on potential theory and which goes back to
  Brolin~\cite{Bro} allows us to derive the convergence of the
  sequence of currents everywhere on $\Lambda$.

\medskip 

We first note that $(d^n+ d^{k(n)})^{-1}\left( [\Per(n, k(n))] - 
(F_n^{*}\omega + F_{k(n)}^*\omega)\right) = dd^c h_n$, with
\begin{equation}\label{eq:potential}
h_n \= \frac1{d^n+ d^{k(n)}}\, \log\,
d_{\PP^1} \left( f_\lambda^n c(\lambda), f_\lambda^{k(n)} c(\lambda)
\right)~.
\end{equation}
Here $d_{\PP^1}$ denotes the spherical distance on $\PP^1$, given in
homogeneous coordinates by $$d_{\PP^1}^2 ([z_0:z_1], [w_0:w_1]) =
\frac{|z_0w_1-z_1w_0|^2}{(|z_0|^2+|z_1|^2)(|w_0|^2+|w_1|^2)}~.$$
Our aim is to show that $\lim_n h_n =0$.  Observe
 that $d_{\PP^1}\le 1$, whence $\limsup_n h_n \le 0$.
 
 The function $h_n$ is not plurisubharmonic, but standard compactness
 results for families of psh functions still apply in this context.
 More precisely,  Proposition~\ref{def:bifurcation}
 implies that $(d^n+ d^{k(n)})^{-1}(F_n^{*}\omega + F_{k(n)}^*\omega)
 = \omega + dd^c \widetilde{g}_n$, where $\widetilde{g}_n$ is a
 sequence of continuous functions converging uniformly on compact
 subsets to a continuous function $g_\infty$.  Locally near any
 parameter we may write $\omega$ as the $dd^c$ of a smooth function,
 so that $ (d^n+ d^{k(n)})^{-1}(F_n^{*}\omega + F_{k(n)}^*\omega) =
 dd^c v_n$, where $v_n$ is a sequence of continuous psh functions
 converging uniformly to a continuous psh function $v_\infty$.
 
 Let $u_n \= h_n + v_n$. It is a upper-semicontinuous function, and $dd^c u_n$
 is a positive closed current, hence $u_n$ is a psh function. The
 sequence  $(v_n)$  converges uniformly, and $\limsup_n h_n \le 0$ so
 $(u_n)$ is a sequence of psh functions which is locally
 uniformly bounded from above.  In particular we have the following
 dichotomy, see~\cite[p.94]{hormander}.  Either $u_n$ (hence $h_n$) converges
 uniformly on compact subsets to $-\infty$; or there exists a
 convergent subsequence $u_{n_j} \cv u$, and $u$ is again a psh
 function. In the latter case, we infer that $h_{n_j} \cv h\=
 u-v_\infty$ and $dd^c h = \lim_j dd^c h_{n_j}$.
 
 In any case, we may pick a subsequence, which we still denote by
 $h_n$ converging to a function $h$ which is either the sum of a
 continuous function and a psh function, or identically $-\infty$.
 The proof will be complete if we show that $h=0$.

\smallskip

\noindent \emph{Claim~1}: If 
$f^n_{\lambda_0}c(\lambda_0)$ converges towards an attracting (or
super-attracting) cycle, then $h(\lambda_0) =0$.  In particular
$h\not\equiv -\infty$ in this case.

Suppose by contradiction that $h(\lambda_0) <0$. The  Hartogs'
Lemma~\cite[p.94]{hormander} applied to the sequence of psh functions
$u_n$ and the function $v_\infty$, implies the existence of $\e_0>0$
and a neighborhood $U$ of $\lambda_0$, such that $h_n|_U < -\e_0<0$
for infinitely many $n$.  In particular,
\begin{equation}\label{eq:claim}
d_{\PP^1} \left( f_\lambda^n c(\lambda), f_\lambda^{k(n)} c(\lambda)
\right) \le \exp( -\e_0\, d^n) \text{ for all } \lambda \in U~.
\end{equation}
Let  $V\subset U$ be a  small open set
 such that $c(\lambda)$ converges to an attracting (or
super-attracting) cycle for  $\lambda \in V$. As $d_{\PP^1} (
f_\lambda^n c(\lambda), f_\lambda^{k(n)} c(\lambda))\le \exp( -\e_0
d^n)$, the cycle is in fact superattracting with multiplicity $d$,
hence totally invariant. For
any $\lambda\in V$, we thus have
$\mathrm{card}\,{\mathcal{E}}(f_\lambda)\ge1$, a contradiction.

\smallskip

\noindent \emph{Claim~2}: $h\equiv 0$ on the activity locus. In
particular, $h\not\equiv -\infty$ when the activity locus is
non-empty.

Indeed by Proposition~\ref{prop:mcmullen} and Claim~1, if $\lambda_0$
belongs to the activity locus, there exists a sequence of parameters
$\lambda_n\cv\lambda_0$ with $h(\lambda_n)=0$. Since $h$ is upper
semicontinuous and non positive, we conclude that $h(\lambda_0)=0$.

\smallskip

\noindent \emph{Claim~3}: $(d^n+ d^{k(n)})^{-1}
[\mathrm{Per}(n,k(n))]\to T\equiv 0 $ on the passivity locus. 

We use the classification of Theorem \ref{p:motion}. 
Let $U$ be a connected open set in the passivity locus.  If
$\mathrm{Per}\, (n,k)\cap U = \emptyset$ for all $n\neq k$ then there
is nothing to prove. Otherwise, $\mathrm{Per}\, (n_0+k_0,k_0)\cap U
\neq \emptyset$ for some $n_0, k_0\ge 1$. By assumption, $c(\lambda)$
is not stably preperiodic, so case (2) of this proposition is excluded.
If now for some $\lambda\in U$, the point $c(\lambda)$ is attracted by a
(super-)attracting cycle, then this behaviour persists throughout $U$
and the first claim implies that $h=0$, hence $T=0$, in $U$. The
remaining case is handled by the following lemma.
  \begin{lem}\label{l:passive}
  Let $U$ be a connected  open set where $c(\lambda)$ 
 stays in the domain of linearization of an irrational
  neutral cycle. Then
  $\bigcup_{n>k}\mathrm{Per}\,(n,k)$ is a (closed) subvariety in
  $U$. Furthermore
  for any compact subset $K\subset U$, there exists a constant
   $C = C(K) < + \infty$ such that
  $$
  \sup_{n>k} \mathrm{Mass}\, [\mathrm{Per}\,(n,k)](K) \le C~.
  $$ 
  \end{lem}
Assuming the lemma for the moment, we proceed with the proof of the
theorem. 

\smallskip 

\noindent \emph{Claim~4}: the restriction of $h$ to any germ of analytic
curve is continuous (if not identically $-\infty$).

The statement is local so let $v$ be a local potential of $T$ near
a parameter $\lambda$. Because $v$ is continuous, $h$ is
continuous iff $\tilde{h} \= h+ v$ is. As $dd^c (h+ v) $ is a positive
current, $\tilde{h}$ is psh.  On the support of $T$, $h$ is
identically zero by the second claim, so $\tilde{h}|_{\mathrm{Supp}\,
  (T)}$ is continuous. The third claim shows that $h$ is pluriharmonic
on the complement of $\mathrm{Supp}\, (T)$. We infer that $\tilde{h}$
is a psh function which is continuous on the support of $dd^c
\tilde{h}$. This is known to imply continuity only in dimension 1. So
let $\Gamma$ be any germ of analytic curve parameterized by a function
$\phi: \Delta\to \Lambda$. The function $\tilde{h}\circ \phi$ is
subharmonic, and continuous on the support of its Laplacian. The
continuity principle (see \cite[p.54]{tsuji}) then implies that
$\tilde{h}\circ \phi$ is continuous.

\bigskip

We are now in position to conclude the proof of the theorem.  When the
activity locus is empty, the theorem follows from Claim~(3).
Otherwise, we pick a convergent subsequence $h_n\cv h$. Claim~(2)
implies that $h \not\equiv -\infty$.  We need to show that $h \equiv
0$. 

Pick $\lambda_0\in\Lambda$ and choose an analytic curve $\Gamma$
containing $\Lambda_0$ which is either quasiprojective if (H1) is
satisfied, or given by the condition (H2). Let $\Omega$ be a
connected component of the intersection of $\Gamma$ with the passivity
locus. The second, third and fourth claims imply that $h$ is harmonic
on $\Omega$, continuous on $\overline{\Omega}$ and zero on the
boundary.

Under the condition (H1), the maximum principle proved
in~\cite[Theorem~III.28]{tsuji} and applied to $h$ gives immediately
$h\equiv 0$.  Under (H2), either the component $\Omega$ is
unbounded so by assumption it is contained in the set where
$f^n_\lambda c(\lambda)$ converges to an attracting periodic cycle,
and $h\equiv 0$. Or $\Omega$ is relatively compact in $\Gamma$ and the
minimum principle applied to $h$ forces $h$ to be
identically zero on $\Omega$.

This shows that $h|_\Gamma =0$, and so $h=0$ everywhere on $\Lambda$.
This concludes the proof of the theorem in the case where 
$e =0$.
\end{proof}
\begin{proof}[Proof of Lemma~\ref{l:passive}]
  The fact that $V\=\bigcup_{n>k} \mathrm{Per}\, (n,k)$ is a closed
  subvariety was already proved in Theorem \ref{p:motion}.  To
  conclude the proof we need to show that for any point $\lambda_0\in
  V$, the local multiplicity $m(n,k,\lambda_0)\in \mathbb{N}$ at
  $\lambda_0$ of the equation $f^n_\lambda\, c(\lambda) = f^k_\lambda
  \, c(\lambda)$ is bounded uniformly on $n$ and $k$. Our aim is to
  choose adequate coordinates both in the parameter space and the
  dynamical plane such that the computation of these multiplicities
  becomes simple. Note that this problem is local both in the parameter
  and dynamical space.

  First, we may assume that for some $n_0$, $f^{n_0}_{\lambda_0}\,
  c(\lambda_0)$ is periodic for period $k_0$, otherwise
  $m(n,k,\lambda_0)=0$ for all $n\neq k$.  Write $g_\lambda =
  f_\lambda^{n_0}$.  Since $c$ is passive, by choosing a suitable
  small neighborhood $U\subset \Lambda$ containing $\lambda_0$, we
  have $\Per(n,k)\cap U \neq \emptyset$ only if $n= n' n_0 + k_0 +r$
  and $k = k' n_0 + k_0 +r$ for some integers $ 0\le r < k_0$ and $n',
  k'$. In particular, in terms of divisors we have
  $$\left[ f_\lambda^{n}\,c(\lambda)=
    f_\lambda^{k}\,c(\lambda)\right]= \sum_{r=0}^{k_0-1}
  \left[ g_\lambda^{n'}\left( f^r_\lambda\,c(\lambda) \right)=
    g_\lambda^{k'}\left(f^r_\lambda\,c(\lambda)\right)\right]~,$$
  taking multiplicities into account.  We bound the
  multiplicities of the equations on the right hand side when $r=0$,
  the arguments being analoguous for $r>0$.

\smallskip

Define $p(\lambda ) \= f^{k_0}_\lambda\, c(\lambda)$.  Then
$g_{\lambda_0}$ has a fixed point at $p(\lambda_0)$.  By reducing $U$
if necessary, and changing coordinates in the parameter plane, we may
suppose that $p(\lambda_0) = 0$ and $g_\lambda(0) = 0$ for all
parameters.  By assumption, $g_\lambda$ is linearizable at $0$ with a
multiplier $\mu = \exp (2 i \pi \theta)$ with $\theta \in
\mathbb{R}\setminus \mathbb{Q}$, and $p(\lambda)$ lies in the domain
of linearization of $g_\lambda$ for every $\lambda\in U$. In other
words, we may assume that $g_\lambda(z) = \mu z$ for all $\lambda\in
U$ and all $z$ small enough.

  Now choose coordinates $\lambda=(\lambda_1, \cdots , \lambda_d)$
  around $\lambda_0 =0$.  Write the expansion into power series of $p$
  as $ p(\lambda) = P_l(\lambda) +  O(|\lambda|^{l+1})$,
  where $P_l$ is a non-zero homogeneous polynomial of degree $l$. Then
  $$g^{n'}_\lambda\, p(\lambda)-g^{k'}_\lambda\, p(\lambda) =
  (\mu^{n'} - \mu^{k'}) P_l(\lambda) + O(|\lambda|^{l+1}).$$
  Since $\mu$ is not a root of unity, we observe that the multiplicity
  of this expression is constant equal to $l$, independently of $n'$
  and $k'$. This concludes the proof of the lemma.
\end{proof}

\begin{proof}[Proof of Theorem~\ref{thm:cvg} in the 
case where $e =1$.]
  
  The problem is purely local.  By the discussion of
  Section \ref{subs:exc}, we may thus assume that $f_\lambda$ is a
  polynomial for all $\lambda$.  In this case, one easily checks that
  the current $\hT$ can be defined in $\Lambda\times \cc$ by the formula
  $$\hT = \lim_n \, dd^c_{\Lambda\times \mathbb{C}}\, \frac1{d^n}
  \log^+|f_\lambda(z)|,$$
  and that $ T =
  \lim_{n\to\infty}d^{-n}\log^+|f^n_\lambda\, c(\lambda)|$ in
  this case. We then write $$[\Per(n,k)] = dd^c \log | f^n_\lambda\,
  c(\lambda)-f^k_\lambda\, c(\lambda)|~,$$
  and the convergence theorem
  thus amounts to proving $\lim h_n =0$ with
\begin{equation}\label{eq:potential2}
h_n \= \frac1{d^n}\, \log 
\frac{|f^n_\lambda\, c(\lambda) -f^k_\lambda\, c(\lambda)|}
{\max\{ |f^n_\lambda\, c(\lambda)|, 1\}}~.
\end{equation}
For each $\lambda$, let $K_\lambda$ be the filled-in Julia set of
$f_\lambda$, i.e. the set of points of bounded orbit in $\mathbb{C}$.
It is a compact subset of $\mathbb{C}$. Moreover,  reducing
$\Lambda$ if necessary, it is included in a fixed disk of radius $C>0$.

Suppose $c(\lambda) \notin K_\lambda$, and fix $\e>0$. Then for $n$
sufficiently large, and any $k<n$, we have $|f^k_\lambda\,
c(\lambda)|\le \e\, |f^n_\lambda\, c(\lambda)|$. So $\lim h_n
(\lambda) = 0$.  Otherwise $c(\lambda) \in K_\lambda$, and it is clear
in this case that $\limsup h_n (\lambda) \le 0$.

The proof now proceeds exactly like in the former case, through the
 proof of Claims~1,~2,~3, and~4.  To adapt the proof  of Claim~1, 
 remark that since $e<2$, any
 superattractive point at finite distance has multiplicity not greater
 than $d-1$.
\end{proof}


\section{Parameterizing the space of critically marked
  polynomials}\label{sec:param}
Let $\mathcal{P}_d$ be the space of all polynomials of degree $d\ge 2$
with $d-1$ marked critical points up to conjugacy by an affine
transformation. A point in $\mathcal{P}_d$ is represented by a
$d$-tuple $(P,c_0, \cdots , c_{d-2})$ where $P$ is a polynomial of
degree $d$, and the $c_i$'s are complex numbers such that $\{ c_0,
\cdots , c_{d-2}\}$ is the set of all critical points of $P$. For
each $i$, $\mathrm{Card}\, \{ j,\, c_j = c_i \}$ is the order of
vanishing of $P'$ at $c_i$. Two points $(P,c_0, \cdots , c_{d-2})$ and
$(P',c'_0, \cdots , c'_{d-2})$ are identified when there exists an
affine map $\phi$ such that $P' = \phi P \phi^{-1}$, and $c'_i =
\phi(c_i)$ for all $0\le i\le d-2$.

The set $\mathcal{P}_d$ is a quasiprojective variety of dimension
$d-1$, and is isomorphic to the quotient of $\mathbb{C}^{d-1}$ by the
finite group of $(d-1)$-th roots of unity acting linearly and
diagonally on $\mathbb{C}^{d-1}$ (see~\cite{silverman} or Proposition
\ref{p:ramified} below).  For instance, $\mathcal{P}_3$ is isomorphic
to the quadric cone $\{z^2 = xy \} \subset \mathbb{C}^3$. When
$d\ge3$, this space admits a unique singularity at the point $(z^d,
0,\ldots, 0)$.  Depending on the problem we consider, we shall
work directly on $\mathcal{P}_d$, viewed as an abstract variety, or
with an ``orbifold parameterization'' of this space by $\cc^{d-1}$,
which we describe shortly. Working with the parameterization is better
suited for computing the masses of the bifurcation currents.
 On the other hand the description of the bifurcation
measure in terms of external rays is simpler in $\mathcal{P}_d$
(see Section \ref{sec:external}).

We now describe our orbifold parameterization of $\mathcal{P}_d$.

Let first  $\tilde{\pi}$ be the map $\tilde{\pi}: \mathbb{C}^{d-1} \to
\mathcal{P}_d$ which maps $(c_1, \cdots, c_{d-2},\alpha)
\in\mathbb{C}^{d-1}$ to the primitive of $z\prod_1^{d-2} (z-c_i)$
whose value at $0$ is $\alpha$.  For the sake of simplicity, we write
$c$ for $(c_1, \cdots, c_{d-2})$, and for $\alpha \in \mathbb{C}$, we have
\begin{equation}\label{e:111}
\widetilde{P}_{c,\alpha} (z) 
= \frac1d\, z^d+\sum_{j=2}^{d-1}(-1)^{d-j}\, \sigma_{d-j}(c)\,
\frac{z^j}{j} + \alpha~,  
\end{equation}
where $\sigma_i(c)$ is the symmetric polynomial in $\{c_j\}_1^{d-2}$
of degree $i$. The critical points of $\widetilde{P}_{c ,\alpha}$ are
$\{0, c_1, \cdots , c_{d-2}\}$, so there is a natural map $\tilde{\pi}:
\mathbb{C}^{d-1} \to \mathcal{P}_d$, which is a finite ramified cover. 

For our purpose, it is better to work with a slightly different
parameterization. For $(c , a) \in \mathbb{C}^{d-2} \times
\mathbb{C}$, we define
$$
P_{c , a} \= \widetilde{P}_{c , a^d}~.
$$
The critical set of $P_{c, a}$ is again given by $(0,c_1, \cdots,
c_{d-2})$, so that we get a natural map $\pi: \mathbb{C}^{d-1}
\to\mathcal{P}_d$, $(c , a) \mapsto (P_{c , a}, 0, c_1, \cdots ,
c_{d-2})$.  The advantage of this parameterization is that all
currents of bifurcation have the same mass, see
Proposition~\ref{p:defcur} below. 
Both orbifold parameterizations of $\mathcal{P}_d$ are very much 
inspired by the
one described in~\cite{BrH1}. Notice however that they chose to have
\emph{centered} polynomials, while we choose to have a  critical
point at $0$.

\begin{prop}\label{p:ramified}
The natural map $\pi: \mathbb{C}^{d-1} \to \mathcal{P}_d$ is a finite
ramified cover of degree $d(d-1)$.

Its critical set is precisely $\{ (c ,a),\, P_{c, a}(0) = 0 \}= \{ a =
0 \}$. The set of critical values of $\pi$ is the set of polynomials
with marked critical points $(P, c_1, \cdots , c_{d-1})$ such that
$P(c_1) = c_1$.

Finally, the number of preimages of a critical value under $\pi$
 is $1$ in the case
of $P(z) = z^d$ and $d-1$ otherwise.  
\end{prop}
\begin{proof}
  We need to understand when two parameters $(c, a)$ and $(c' , a')$
 map to the same point in $\mathcal{P}_d$, 
i.e. when there exists an affine map $\phi$
  such that $\phi \circ P_{c ,a} \circ \phi^{-1} = P_{c',a'}$,
  $\phi(c_i)= c'_i$ for all $i$, and $\phi(0)=0$. The latter fact
  implies that $\phi = \zeta z$ for some $\zeta$. But the polynomials
  $P_{c ,a}$, $P_{c' ,a'}$ have the same leading monomial, so
  $\zeta^{d-1} =1$. We conclude that $\pi (c ,a) = \pi( c', a')$ iff
  there exists a $(d-1)$-th root of unity $\zeta$ such that $\zeta c_i =
  c'_i$ for all $i$, and $\zeta a^d = (a')^d$.  The proposition
  follows easily from this fact.
\end{proof}


\section{Higher bifurcation currents and the bifurcation measure}
\label{sec:poly}

In this section, we work exclusively with the parameterization 
$(c , a) \mapsto P_{c ,  a}$ of the
space of polynomials with all critical points marked, as described before,
and write $c_0 \= 0$.  We introduce higher bifurcation currents
and the bifurcation measure.  This leads in particular to the proof of
Theorem~\ref{thm:approx-all}.

\subsection{Basics}
Recall that the filled-in Julia set $K(P)$ of a polynomial $P$ is the
set of points with bounded orbits. The filled-in Julia set is
compact, and its boundary is the Julia set of $P$.
\begin{defi}
  The \emph{connectedness locus} is the set of parameters in
  $\mathbb{C}^{d-1}$ for which the filled-in Julia set of $P_{c ,a}$
  is connected. Equivalently, it is the set of parameters for which
  all critical points have bounded orbit.  We denote it by $\C$.
\end{defi}
The following fundamental result comes from~\cite{BrH1}:
\begin{prop}\label{p:compact}
The connectedness locus is compact in $\mathbb{C}^{d-1}$.
\end{prop}

Recall that the Green function of a polynomial $P$ of degree $d$ is by
definition $g_P = \lim_n d^{-n} \log^+|P^n|$.  Let $z\mapsto g_{P{c
    ,a}}(z)$ be the Green function of $P_{c,a}$. We define a function
on the parameter space $\cc^{d-1}$ by putting $G(c , a)= \max \set{
  g_{P{c ,a}}(c_k),~0\leq k\leq d-2}$. We will give more properties of
$G$ later on. The previous proposition is an obvious consequence of
the following estimate, which is a result of~\cite{BrH1}. Because we
work in different coordinates, we provide a detailed proof.
\begin{prop}\label{l:growth}
$G(c ,a) = \log^+ \max \{|a|, |c_k| \} + {O}(1)$.
\end{prop} 
 The proof is based on the
  following two lemmas.
  \begin{lem}\label{l:lower}
 For all $z\in \mathbb{C}$, we have   $$
  g_{P_{c, a}}(z) \le \log \max \{ |z| , A \} + \frac{\log C}{d-1}~,
  $$
  for some constant $C$ depending only on $P$ and with $A = \max \{
  |c_k|, |a| \}$.
  \end{lem}
\begin{lem}\label{l:upper}
  For all $z\in \mathbb{C}$, we have $$
  \max\{ g_{P_{c, a}}(z), G(c
    ,a)\} \ge \log \left| z - \delta\right| - \log 4~,
  $$
where $\delta = \sum c_k /(d-1)$.
\end{lem}
\begin{proof}[Proof of Proposition~\ref{l:growth}]  Let 
  $A=\max\set{\abs{a}, \abs{c_k}}$.  Since $|c_k| \le A$ for all $k$,
  Lemma~\ref{l:lower} yields $g_P(c_k) \le \log A + \log C /(d-1)$.
  The same estimate holds for $g_P(0)$, whence $$G(c, a)= \max
  \{g_{P_{c , a}} (0), g_{P_{c , a }}(c_k) \} \le \log \max \{ |c_k|,
  |a| \} + \log C / (d-1).$$
  
  To get the estimate from below, we apply Lemma~\ref{l:upper}. Let
  $A\geq 2$. Assume first that all complex numbers $\delta, c_1 -
  \delta, \cdots , c_k - \delta$ have modulus $< A/2$, in which
  case $A=a$. Then by lemma \ref{l:upper} we get $$d\times G(c ,a)\geq
  \max\set{g_{P_{c , a}}(a^d), G(c , a)}\geq \log\left| a^d -
    \delta\right| - \log 4 \geq d \log\abs{A}- \log 8.$$
  Here the
  first inequality follows from the fact that $a^d= P_{c , a}(0)$ and
  the last one holds because $\abs{A^d-\delta}\geq A^d/2$.

In the opposite case, among the complex numbers $\delta, c_1 - \delta,
\cdots , c_{d-2} - \delta$, at least one has modulus $\ge A/2$.  So we
deduce from lemma \ref{l:upper} that $G(c, a)\ge \log \max \{|c_k|,
|a| \} -\log 8$.
\end{proof}

\begin{proof}[Proof of Lemma~\ref{l:lower}]
  By definition, $P_{c ,a} = \widetilde{P}_{c , a^d}$,
  and~\eqref{e:111} yields
$$
|P(z)| \le d^{-1} |z|^d \left( 1 + d \times \max \left\{
    \frac{\sigma_{d-j}(c)}{j|z|^{d-j}}, \frac{|a|^d}{|z|^d} \right\}
\right)~,$$
where $\sigma_{d-j}(c)$ is the symmetric polynomial of
degree $d-j$ in the $c_k$'s. When $|z|\ge A = \max \{ |c_k|,
|a| \}$, we infer that $|P(z)| \le C |z|^d$ for a real number $C$ depending
only on $d$. By the maximum principle, $|P(z)| \le CA^d$ when $|z|\le
A$.  These estimates classically imply the statement of the lemma.
\end{proof}
\begin{proof}[Proof of Lemma~\ref{l:upper}]
  Fix a polynomial $P$.  Recall that the B{\"o}ttcher coordinate is a
  \emph{univalent} function $\varphi_P: \{ g_P > G(P)\} \to \mathbb{C}$
  satisfying the functional equation
  $\varphi_P\circ P(z) = \phi_P(z) ^d$. In particular, $\log
  |\varphi_P| = g_P$ where the left hand side is defined.  It is actually
  possible to choose a branch of the $(d^n)^\mathrm{th}$-root such that 
 $\varphi_P(z)\= \lim_{n\to\infty} (P^n(z))^{1/d^n}$
  where the convergence is uniform, see~\cite{mibook}.  When $P=P_{c ,
    a}$, a direct induction shows that $P^n(z) = z^{d^n} -
  \frac{d^n}{d-1} \sigma_1(c) z^{d^n-1} + \mathrm{l.o.t}$, so that
  $\varphi(z) = z - \delta + O(1/z)$ with $\varphi \= \varphi_{P_{c ,a}}$ and
  $\delta \= \sigma_1(c)/(d-1) = \sum c_k/(d-1)$. Let $\psi \=
  \varphi^{-1}$. It is  a univalent function on $\mathbb{C}\setminus
  \overline{D}\left(0, \exp\,G(c, a)\right)$ and $\psi(z) = z
  +\delta + O(1/z)$ at infinity.

  We can thus apply~\cite[Corollary~3.3]{BrH1} --which is a version of
  the Koebe $1/4$-Theorem-- to the univalent
  function $\psi - \delta$.  This yields $\psi (\mathbb{C}\setminus
  \overline{D}(0,r)) \supset \mathbb{C} \setminus
  \overline{D}(\delta,2r)$ when $r> \exp\, G(c ,a)$. Now pick any
  $z\in \mathbb{C}$, and write $g_P(z) = \log r$. Then, since 
$\abs{\phi(z)}=r$,   $z \notin \psi
  \left( \mathbb{C}\setminus \overline{D}(0, 2\max\{ r, \exp G(c ,
  a)\})\right)$, thus $|z-\delta|\le 4 \max \{ r , \exp \, G( c , a) \}$. We
  conclude by taking  logarithms in both sides.
 \end{proof}

We close this paragraph by introducing some terminology.
\begin{defi}
  A polynomial is said to be \emph{critically finite} if all its critical
  points are preperiodic. It is of \emph{Misiurewicz type} if all its critical
  points are mapped to  repelling periodic orbits.
\end{defi}

It is a classical fact that if all critical points are strictly
preperiodic, then the polynomial is Misiurewicz (see \cite[p.92]{CG}).

\subsection{Currents and measure of bifurcation}

We now define $d-1$ positive closed $(1,1)$ currents on the parameter
space $\mathbb{C}^{d-1}$, each describing the bifurcation of one
critical point. This part is essentially a rephrasing of the results
of Section~\ref{sec:bif} in this specific situation.

As before, consider the Green function $z\mapsto g_{P_{c ,a}}(z)$ 
of the  polynomial $P_{c ,a}$.
 It is a continuous
positive subharmonic function, and the filled-in Julia set of $P_{c ,a}$
equals $\{ g_{P_{c ,a}} =0 \}$.  The function $((c , a), z) \mapsto g_{P_{c
    ,a}}(z)$ is psh in
$\cc^{d-1} \times \mathbb{C}$. In particular, for all $0\le k \le
d-2$, $g_k \=g_{P_{c , a}}(c_k)$ induces a psh function on
$\mathbb{C}^{d-1}$.  

\medskip

Let $\omega$ be the Fubini Study metric on $\cc^{d-1}$. If $T$ is a
positive closed current of bidegree $(k,k)$ on $\cc^{d-1}$, we let
$\mathrm{Mass}(T)$ be its {\em
  mass} relative to $\omega$, that is 
$$\mathrm{Mass}(T)=\int_{\cc^{d-1}} T\wedge \omega^{d-1-k}.$$
It is
convenient to write $d^c= (\fr-\overline\fr)/2i\pi$, so that if
$z=(z_1, \ldots, z_{d-1})\in \cc^{d-1}$, the mass of $dd^c
\log\abs{z_1}$ is 1. In terms of potentials, if
$u=\log^+\max \{ |z_i|\}+O(1)$, then the mass of $dd^cu$ is again 1.
 We also repeatedly use the following
B{\'e}zout-type estimate (see for instance \cite{FS}): if $T_1,\ldots,
T_q$ are positive closed currents with finite mass and locally bounded
potentials (so that their wedge product is well defined, see \cite{BT,
  de}), then $$\mathrm{Mass}(T_1\wedge \ldots \wedge T_q)\leq
\mathrm{Mass}(T_1)\cdots \mathrm{Mass}(T_q).$$
For $0\leq k\leq d-2$, we  set $T_k = dd^c g_k=dd^cg_{P_{c , a}}(c_k)$.

\begin{prop}\label{p:defcur}
  For each $0\le k \le d-2$, $T_k$ is a positive closed current of
  bidegree $(1,1)$, with continuous potential and mass 1 in
  $\cc^{d-1}$.  The support of $T_k$ is precisely the activity locus
  of the critical point $c_k$, and is equal to the boundary of the
  closed set $\{ g_{P_{c, a}}(c_k) =0 \}$.
\end{prop}
By summing these currents, we get a current whose support is the
bifurcation locus.  This is the {\em bifurcation current}, originally
 considered in~\cite{DeM1}. It was proved to equal $dd^c\,
\mathrm{Lyap} (P)$ in~\cite{DeM2} for families of  rational maps.
For convenience, we restate and reprove her results in our context of polynomial maps.
\begin{prop}[\cite{DeM1,DeM2}]\label{p:demarco}
  Define $T_\mathrm{bif}\= (d-1)^{-1} \sum_0^{d-2} T_k$. This is a
  positive closed $(1,1)$ current of mass $1$, whose support is
  precisely the complement of the set of 
  polynomials structurally stable on their Julia sets.  
Further,  $T_\mathrm{bif}=d^{-1} \, dd^c\,
  \mathrm{Lyap} (P_{c ,a})$, where $\mathrm{Lyap}(P)$ is the Lyapunov
  exponent of  $P$ relative to its measure of maximal entropy.
\end{prop}
We can now define higher dimensional bifurcation currents by intersecting
the currents $T_k$. The currents $T_k$ admit locally continuous
potentials, so  their intersections are well-defined.
We first note the following fact.
\begin{prop}\label{p:zero}
For each $k$, we have $T_k \wedge T_k =0$.
\end{prop}
\begin{defi}
  For each $1\le l \le d-1$, we define the \emph{$l$-bifurcation
    current} to be:
  $$
  T_\mathrm{bif}^l = \frac{(d-l-1)!}{(d-1)!}\,\left( \sum_{i=0}^{d-2} T_i
  \right)^{\bigwedge l}
$$
It is a positive closed current of bidegree $(l,l)$.
\end{defi}
These higher currents of bifurcation were first considered
by~\cite{bas-ber} in the more general context of rational maps.

\begin{prop}\label{p:bifcur}
  The current $T^l_\mathrm{bif}$ is a non zero positive closed $(l,l)$
  current, and its trace measure does not charge pluripolar (hence
  analytic) sets. It has finite mass $1$, and its support is included
  in the set of polynomial for which (at least) $l$ critical points
  are active.
\end{prop}
\begin{rmk}
  As we noticed in the introduction, for $l>1$ the support of
  $T^l_\mathrm{bif}$ is not equal to the set where $l$ critical points
  are active. The following example was indicated to us by A.~Douady,
  and was studied in great detail in the Ph. D. thesis of
  P.~Willumsen~\cite{willumsen}.
\end{rmk}
\begin{exam}\label{ex:douady}
  Take $P = z + 1/2\, z^2 + z^3$. Any cubic polynomial $Q$ close
  enough to $P$ has one or two fixed points near $0$. 
 In case there are two fxed points,  denote by
  $\mu_1$ and $\mu_2$ their multipliers. The residue of $1/(P-z)$ at
  $0$ equals $-4$, so $(\mu_1 -1)^{-1} + (\mu_2 -1)^{-1}$ is close to 
$-4$. It is not difficult to check that in this case,
  either $\mu_1$ or $\mu_2$ has modulus $<1$ (see \cite{mibook}), 
  so the associated fixed point attracts one critical 
  point. Therefore, by
  Proposition~\ref{p:bifcur} such a $Q$ cannot lie in the support of
  the measure $T^2_\mathrm{bif}$. 
  
  Hence near $P$, the measure $T^2_\mathrm{bif}$ has support in the subvariety of
  parameters with a fixed point of multiplier 1. Since this measure does not charge
  curves, $P$ does not lie in the support either. 
  
  On the other hand,
  the two critical points of $P$ are complex conjugate, so both of
  them are attracted by the parabolic fixed point $0$. At least one of them is active:
  indeed if not, for nearby parameters we would
  have a critical point attracted by a repelling fixed point. By symmetry, we conclude that 
  both are active. 
\end{exam}

\medskip

For $l = d-1$, we get a positive measure that we denote by
$\mu_\mathrm{bif}$. Its study will be our main focus in the sequel. 
The following proposition summarizes its first properties.
\begin{prop}\label{p:bifmeas}~
\begin{itemize}
  \itm $\mu_\mathrm{bif} = T_0 \wedge \cdots \wedge T_{d-1}$. 
  
  \itm The measure $\mu_\mathrm{bif}$ is a positive probability
  measure, supported on the boundary of the connectedness locus $\C$.
  
  \itm It does not charge pluripolar sets. In particular, it does not
  charge analytic subsets.
  
  \itm It is the pluricomplex equilibrium measure of the compact set
  $\C \subset \mathbb{C}^{d-1}$. In particu\-lar, its support is the
  Shilov boundary of $\C$.
\end{itemize}
\end{prop}
We refer to~\cite{klimek} for basic notions in pluripotential theory,
such as ``pluripolar sets'', or ``equilibrium measure''. See \cite{BT3}
for the notion of Shilov boundary and the fact that it equals
$\supp(\mu_{\mathrm{Bif}})$. The last statement of the proposition
says in particular that $\mu_\mathrm{Bif}$ is natural from the point
of view of complex analysis.
\begin{proof}[Proof of Proposition~\ref{p:defcur}]
  By definition, $T_k = dd^cg_k $ where $g_k= g_{P_{c ,a}}(c_k)$ is
  psh and continuous. The continuity of $g_k$ is classical in this
  context, or follows directly from 
  Proposition~\ref{def:bifurcation}. The support of
  $T_k$ is the activity locus of $c_k$ by Theorem~\ref{thm:support}.
  
  Assume now that  $g_k(c, a) >0$. Then $c_k$ is attracted towards the
  superattracting point $\infty$, hence $c_k$ is passive at $(c ,a)$,
  so $T_k=0$ near $(c ,a)$. It is also clear that $T_k=0$ in the
  interior of $\{g_k =0 \}$, so that $\mathrm{Supp}\, T_k \subset
  \partial \{ g_k =0 \}$. Conversely, by the maximum principle,
 $g_k$ cannot be pluriharmonic in an open set $U$ intersecting 
 $\partial \{ g_k =0 \}$. We
  conclude that $\mathrm{Supp}\, T_k = \partial \{ g_k =0 \}$.

  \smallskip It remains to compute the mass of the current $T_k$.
Proposition~\ref{l:growth} asserts that  $g_k(c, a)\leq \log\max\set{\abs{c_k},
  \abs{a}} +O(1)$ so $\mathrm{Mass}(T_k)\leq 1$.

 If now $L$ is a complex line in $\cc^{d-1}$ and $T$ is a positive current 
of bidegree (1,1) with continuous potential, we define the restriction
$T\rest{L}$
of $T$ to $L$ by first restricting the potential and then taking $dd^c$. 
The restriction $T\rest{L}$ is a positive measure on $L$, and 
$\mathrm{Mass}(T\rest{L})\leq \mathrm{Mass}(T)$ by Bezout ($T\rest{L}=
T\wedge [L]$).

To get the opposite inequality for $\mathrm{Mass}(T_k)$, we restrict $T_k$
to the line of unicritical polynomials $L:=\set{c_1=\ldots=
  c_{d-2}=0}$. It is enough to prove that
$\mathrm{Mass}(T_k\rest{L})=1$. 
With our parameterization, unicritical
polynomials are of the form $\unsur{d}z^d+a^d$. For every $k$, 
$g_k(0, a)$ equals $g_a(0)$, where $g_a$ is the Green function of 
$\unsur{d}z^d+a^d$. From proposition~\ref{l:growth} we infer that $g_a(0)=
\log\abs{a}+ O(1)$ at infinity, so the mass of the restriction is 1. \end{proof}

\begin{proof}[Proof of Proposition~\ref{p:demarco}]
  All currents $T_k$ have mass $1$, hence $T_\mathrm{bif}$  also has
  mass $1$.

  By~\cite{MSS}, a polynomial is structurally stable on $J$ iff all its
  critical points are passive. Using Proposition~\ref{p:defcur}, we conclude
  that a polynomial is structurally unstable iff it belongs to the union
  of the support of the $T_k$'s. The latter is precisely the support
  of $T_\mathrm{bif}$.
  
  The Lyapunov exponent of a polynomial $P$ can be computed in terms
  of the Green function, see~\cite{manning,prz}. We get $\mathrm{Lyap}(P_{c ,
    a}) = \log d + \sum_0^{d-2} g_{P_{c , a}}(c_k)$. We already know that $dd^c
  g_{P_{c , a}}(c_k) = T_k$, whence $dd^c\, \mathrm{Lyap}(P_{c , a})=
  \sum_0^{d-2} T_k = (d-1) \times T_\mathrm{bif}$.
\end{proof}

\begin{proof}[Proof of Proposition~\ref{p:zero}]
  As before write $g_k = g_{P_{c ,a}}(c_k)$, so that $dd^c g_k = T_k$.
  The sequence of continuous psh functions $\max \{ g_k, \e \}$ 
  decreases to $g_k$ when $\e$ decreases to zero. From~\cite{BT}, we
  infer that $dd^c \max \{ g_k ,\e\} \wedge dd^c g_k $ converges weakly
  to $dd^c g_k \wedge dd^c g_k = T_k \wedge T_k$. Since $g_k$ is
  pluriharmonic where it is positive, the support of $dd^c
  \max \{ g_k ,\e\} \wedge dd^c g_k $ is contained in the intersection
  $\{ g_k = \e \} \cap \{ g_k = 0 \}$ which is empty. Hence $T_k
  \wedge T_k =0$.
\end{proof}

\begin{proof}[Proof of Proposition~\ref{p:bifcur}]
  The fact that $T^l_\mathrm{bif}$ is a positive closed $(l,l)$
  current whose trace measure does not charge pluripolar set is a
  consequence of the definition and the fact that all currents $T_l$
  have continuous potentials, see~\cite{de}. 
  
  Take a polynomial $P$ for which at most $l-1$ critical points are
  active. Then $T_j =0$ for at least $d-l$ distinct integers $j$.
  Combining this with Proposition~\ref{p:zero}, we conclude that
  $T^l_\mathrm{bif}=0$ near $P$. The support of the $l$-th bifurcation
  current is thus contained in the set of polynomials with at least $l$
  active critical points.

 It remains to compute the mass of $T^l_\mathrm{bif}$ for every $2\leq
 l\leq d-1$. To this end,
 we rely on the following computation. 
\begin{lem}\label{l:wedge}
Let  $I=(i_1, \cdots, i_l)$ be a  multi-index with  $l$ distinct
entries  in $\{ 0, \cdots, d-2\}$. Then
\begin{equation}\label{e:wedge}
T_{i_1}\wedge \cdots \wedge T_{i_l} = (dd^c)^l G_I \text{ with } G_I =
\max \{ g_{i_1} , \cdots , g_{i_l} \}~.
\end{equation}
\end{lem}
Let us continue with the proof of the proposition. 
By Lemma~\ref{l:wedge}, we have that $(dd^c)^{d-1} G = T_0 \wedge \cdots \wedge
T_{d-2}$ where $G = \max \{ g_0, \cdots , g_{d-2} \}$.
Recall from Proposition~\ref{l:growth} that  $G(c ,a)= \log
\max\set{\abs{c_k}, \abs{a}}+O(1)$. 
Standard estimates in pluripotential theory (see
e.g. \cite[p. 212]{klimek}) imply that 
 the measure $(dd^c)^{d-1}G$
has unit mass.  Let now $i_1, \cdots, i_{d-1}$ be any reordering of the $d-1$
integers $\{ 0, \cdots, d-2\}$. Then 
\begin{eqnarray*}
1= \mathrm{Mass}\, (T_{i_1}\wedge \cdots \wedge T_{i_{d-1}})
&\le &
 \mathrm{Mass}\,\, T_{i_{d-1}} \times  
\mathrm{Mass}\, (T_{i_1}\wedge \cdots \wedge T_{i_{d-2}})
\\
& = & 
\mathrm{Mass}\, (T_{i_1}\wedge \cdots \wedge T_{i_{d-2}})
\le \prod_{j=1}^{d-1}  \mathrm{Mass}\, T_{i_j} =1~.
\end{eqnarray*}
Here the first inequality is the B{\'e}zout-type estimate for currents.
We conclude that $\mathrm{Mass}\, (T_{i_1}\wedge \cdots \wedge
T_{i_{d-2}})=1$.  Proceeding by descending induction, we infer that
for any choice of distinct integers $\{i_1, \cdots, i_l\} \subset \{
0, \cdots , d-2\}$, the mass of $(T_{i_1}\wedge \cdots \wedge
T_{i_l})$ equals $1$. Whence
\begin{equation*}
  \mathrm{Mass}\, \left(\sum_0^{d-2} T_i \right)^{\wedge l} = \sum
  _{i_1\neq i_2 \cdots \neq i_l} \mathrm{Mass}\, (T_{i_1} \wedge \cdots \wedge
  T_{i_l}) = (d-1) (d-2) \cdots (d-l) = \frac{(d-1)!}{(d-1-l)!}~.
\end{equation*}
This concludes the proof of the proposition.
\end{proof}

\begin{proof}[Proof of Lemma~\ref{l:wedge}]
We  prove that $T_0 \wedge \cdots \wedge T_l = (dd^c)^{l+1} G_l$
with $ G_l = \max \{ g_0, \cdots , g_l\}$ for all $l\le d-1$. The
same proof gives the general equality~\eqref{e:wedge}.

We proceed by induction on $l$.  The statement is clear for $l=0$.
Suppose that we know that $T_0 \wedge \cdots \wedge T_{l-1} = (dd^c)^l
G_{l-1}$. We want to compute $T_0 \wedge \cdots \wedge T_l$. We
proceed as follows:
\begin{eqnarray*}
  T_0 \wedge \cdots \wedge T_{l-1} \wedge T_l 
&=&
  dd^c \left( g_l\, (dd^c)^{l} G_{l-1} \right) 
=  dd^c \left( G_l\, (dd^c)^l G_{l-1} \right) 
  \\
&=&  dd^c \left( G_{l-1}\, (dd^c)^{l-1} G_{l-1}\wedge dd^c G_l \right)
\\
&=&  dd^c \left( G_l\, (dd^c)^{l-1} G_{l-1} \wedge dd^c G_l \right)
= (dd^c)^{l+1} G_l ~.
\end{eqnarray*}
Let us justify these sequence of equalities. The first one is by
definition. For the second one,  observe that the support of $(dd^c)^l
G_{l-1}$ is contained in the intersection of the supports of $T_0,
\ldots ,T_{l-1}$, and $g_0= \cdots = g_{l-1}=0$ on this set. So $g_l =
G_l$ there. The third one is just reordering
the wedge product. For the fourth, note that on $\{G_{l-1} < G_l \}$
we have $G_l = g_l > 0$. Thus $G_l$ is pluriharmonic and $dd^c G_l =
0$.  The last equality is obtained by repeating the same argument
$l-1$ times.
\end{proof}

\begin{proof}[Proof of Proposition~\ref{p:bifmeas}]
  The equation $\mu_\mathrm{bif} = T_0 \wedge \cdots \wedge
  T_{d-1}$ is a consequence of Proposition~\ref{p:zero}.  The first
  two facts are consequences of Proposition~\ref{p:bifcur}. It  remains
  to prove the last one.
  
  Consider the function $G = \max g_k$. This is a continuous psh
  function on $\mathbb{C}^{d-1}$, and $\{G =0 \} = \C$. We claim
  that
\begin{equation}\label{e:379}
G = \sup \left\{ u \text{ psh }, \, u - \log^+ \max
  \{|a|, |c_k| \} \le  {O}(1), \, u \le 0 \text{ on } \C \right\}~. 
\end{equation}
In the terminology of~\cite{klimek}, $G$ is the pluricomplex Green function
of $\C$ with pole at infinity. The equilibrium measure of $\C$ is by
definition the Monge-Amp{\`e}re measure of $G$. It  thus
equals  $(dd^c)^{d-1} G = \mu_\mathrm{bif}$ by~\eqref{e:wedge}.

The proof of~\eqref{e:379} is  standard, we include it for completeness. 
Take any psh function $u$ with $u\le 0$ on $\C$
and $u - \log^+ \max \{|a|, |c_k| \} \le O(1)$. Choose $\e>0$, and set
$u_\e\= (1-\e) u - \e$. Then $u_\e < G$ on a neighborhood of $\C$ as
$G\le0$; this also holds 
 in the complement of some ball $B(0,R)$ by
Proposition~\ref{l:growth}.  Now pick any parameter $(c, a)\notin
\mathcal{C}$. Consider an  open set $ \Omega\subset
\mathbb{C}^{d-1}\setminus \C$, $(c, a)\in\om$, 
with boundary contained in $U\cup \cc^{d-1}\setminus
B(0,R)$, where $U$ is some small neighborhood of
$\mathcal{C}$. On $\Omega$  we have $(dd^c)^{d-1} G =0$, so $G$ is maximal
as a psh function, see~\cite{klimek}. This implies that 
$u_\e \le G$ on $\Omega$. We have thus
proved that $u_\e \le G$ everywhere. By letting $\e\cv 0$, we conclude
that $u \le G$. 
\end{proof}


\subsection{Density of Misiurewicz points}
We now aim at proving Theorem~\ref{thm:approx-all} cited in the
introduction.  It will be a consequence of the following more general
result.
\begin{thm}\label{thm:allbidegree}
  Choose a reordering $\{ i_0, \cdots , i_{d-2}\}$ of $\{ 0 , \cdots ,
  d-2\}$ and fix $0\le l\le d-2$. For each $0\le j < l$, let
  $(n_j,m_j)$ be any couple of integers such that $n_j > m_j$. 
For $l=0$, set $W=\cc^{d-1}$, and for $l\geq 1$, let
  $W\subset\mathbb{C}^{d-1} $ be the analytic subvariety (not
  necessarily reduced) consisting of parameters $c ,a$ for which $P_{c
    , a}^{n_j} (c_{i_j}) = P_{c , a}^{m_j} (c_{i_j})$ for all $0\le j <l$.
  Then, the following hold. 
\begin{enumerate}
\item
All irreducible components of $W$ have codimension $l$.
\item
Pick any $n >k$, and denote by $\mathrm{Per}\, (n , k)$ the
set of polynomials such that $P^{n} (c_{i_l}) = P^{k} (c_{i_l})$.  Then the
variety $\mathrm{Per}\, (n , k)\cap W$ is a hypersurface in $W$.
\item
For any sequence of integers $k(n) <n$, one has 
$$
\lim_{n\to \infty} \frac1{d^n}\,[\mathrm{Per}\, (n , k(n))\cap W] =
T_{i_l} \wedge [W]~.
$$
in the weak topology of currents.
\end{enumerate}
\end{thm}
\begin{cor}\label{cor:cvgall}
  Fix a collection of distinct integers $i_0, \cdots , i_l \in \{ 0,
  \cdots , d-2\}$, and for any $(n_0, \cdots , n_l) \in
  (\mathbb{N}^*)^l$ choose a collection of integers $k(n_0) <n_0,
  \cdots , k(n_l) < n_l$. Then define the analytic subset $W_{n_0,
  \ldots, n_l} =
  \bigcap_{j=1}^l\, [ P^{n_j} (c_{i_j}) = P^{k(n_j)} (c_{i_j})]$.

Then $W_{n_0,
  \ldots, n_l}$ has pure codimension $l+1$ and  
\begin{equation}\label{e:induc}
\lim_{n_l \to \infty} \cdots \lim_{n_1\to\infty} \,\lim_{n_0\to
  \infty} \,
\frac1{d^{n_l+ \cdots +n_1 + n_0}}\, [W_{n_0,
  \ldots, n_l}] =T_{i_l} \wedge \cdots \wedge
T_{i_1}\wedge T_{i_0}~.
\end{equation}
\end{cor}
With this result in hand, 
Theorem~\ref{thm:approx-all} follows easily:
\begin{proof}[Proof of Theorem~\ref{thm:approx-all}]
  The current $T^l_\mathrm{bif}$ is a multiple of the positive closed
  $(l,l)$ current $\left(\sum_0^{d-2} T_i \right)^{\wedge l} = \sum
  _{i_1\neq i_2 \cdots \neq i_l} T_{i_1} \wedge \cdots \wedge
  T_{i_l}$. The theorem is then a direct consequence of
  Corollary~\ref{cor:cvgall}.
\end{proof}
From this we also deduce Corollary~\ref{cor:mis}.  Notice that a
critically finite polynomial is hyperbolic iff its critical points are
all {\em periodic}. Recall that a polynomial is Misiurewicz if all
critical points are preperiodic to repelling periodic orbits.
\begin{cor}\label{cor:meas}
  There exists a sequence of atomic measures supported on the set of
  Misiurewicz parameters (resp. on critically finite hyperbolic
  parameters) and converging to $\mu_\mathrm{bif}$. In particular, the
  support of $\mu_\mathrm{bif}$ is contained in the closure of the set
  of 
  Misiurewicz points (resp. of hyperbolic critically finite  parameters).
\end{cor}
\begin{proof}
Corollary~\ref{cor:cvgall} implies that
$$
\lim_{n_{d-2} \to \infty} \cdots \lim_{n_0\to \infty} \,
\frac1{d^{n_{d-2} + \cdots + n_0}}\, \left[ \bigcap_{j=0}^{d-2} \left\{
    P^{n_j}_{c ,a} (c_j) =c_j \right\} \right] =T_0 \wedge
\cdots \wedge T_{d-2}~.
$$
On the left hand side we have a sequence of atomic measures
supported on critically finite and hyperbolic parameters. The right
hand side is $\mu_\mathrm{bif}$ by Proposition~\ref{p:bifmeas}.

The same Corollary~\ref{cor:cvgall} again implies that 
$$
\lim_{n_{d-2} \to \infty} \cdots \lim_{n_0\to \infty} \,
\frac1{d^{n_{d-2} + \cdots + n_0}}\, \left[ \bigcap_{j=0}^{d-2} \left\{
    P^{n_j}_{c ,a} (c_j) = P_{c , a}^{n_j-1}(c_j) \right\} \right] =T_0 \wedge
\cdots \wedge T_{d-2}~.
$$
To conclude the proof we use the same argument as in
Corollary~\ref{cor:strpreper}.  Each divisor $H_j \=\left\{ P^{n_j}_{c
    ,a} (c_j) = P^{n_j-1}_{c , a}(c_j) \right\}$ can be decomposed as a sum of
two divisors $H_j = \mathrm{Preper}_j\, + \mathrm{Fix}_j$ where
$\mathrm{Preper}_j$ is the set of points in $H_j$ for which $c_j$ is
\emph{strictly} preperiodic, and $\mathrm{Fix}_j$ is supported on 
the set on which $c_j$ is \emph{fixed}.  As $T_j$ does not charge hypersurfaces, we
have $d^{-n_j} \mathrm{Preper}_j\to T_j$, therefore
$$
\lim_{n_{d-2} \to \infty} \cdots \lim_{n_0\to \infty} \,
\frac1{d^{n_{d-2} + \cdots+ n_0}}\, \left[\bigcap_{j=0}^{d-2}
  \mathrm{Preper}_j\,\right] =T_0 \wedge \cdots \wedge T_{d-2}~.
$$
The supports of the atomic  measures on the left hand side
are contained  in the set of Misiurewicz parameters.
\end{proof}

\begin{proof}[Proof of Corollary~\ref{cor:cvgall}]
  The proof is by induction on $l$.  For $l=0$, this is Theorem
  \ref{thm:cvg} for the critical point $c_{i_0}$ (or equivalently,
  the statement (3) of Theorem~\ref{thm:allbidegree} with $l=0$). 
Suppose now that it is true for some integer $l$, and pick
  $i_{l+1}\in \{ 0, \cdots , d-2\}$.  For $k(n_{l+1})< n_{l+1}$, 
define $\mathrm{Per}\,  (n_{l+1},k(n_{l+1}))$ as the set of
 $(c, a)\in\mathbb{C}^{d-1}$ such that $P_{c
    ,a}^{n_{l+1}} (c_{l+1}) = P_{c ,a}^{k(n_{l+1})} (c_{l+1})$. 
For ease of notation we write $n$ for $(n_0,\ldots, n_l)$.  By
  Theorem~\ref{thm:allbidegree}, $W_{n, n_{l+1}} = \mathrm{Per}\,
  (n_{l+1},k(n_{l+1})) \cap W_n$ has pure codimension $l+2$, and
  $\lim_{n_{l+1}\to \infty} d^{-n_{l+1}}[W_{n, n_{l+1}}] = T_{i_{l+1}}
  \wedge [W_n]$. Now, our inductive hypothesis asserts that  $
  d^{-(n_{l}+\cdots +n_0)}[W_n]\cv T_{i_l} \wedge \cdots \wedge T_{i_1}\wedge
  T_{i_0}$.  Since, $T_{i_{l+1}}$ has
  \emph{continuous potential}, we deduce (see~\cite{de}) that 
  $$T_{i_{l+1}}\wedge d^{-(n_{l}+\cdots +n_0)}[W_n] \cv
  T_{i_{l+1}}\wedge T_{i_l} \wedge \cdots\wedge T_{i_0}$$
  as
  $n_0,\ldots,n_l \cv\infty$. By the convergence Theorem~\ref{thm:cvg}
  we get that $d^{-n_{l+1}}\mathrm{Per}\, (n_{l+1},k(n_{l+1}))\cv
  T_{l+1}$ as $n_{l+1}\cv\infty$, whence $$d^{-n_{l+1}}\mathrm{Per}\,
  (n_{l+1},k(n_{l+1}))\wedge d^{-(n_{l}+\cdots +n_0)}
  [W_n]\longrightarrow T_{i_{l+1}}\wedge T_{i_l} \wedge \cdots \wedge
  T_{i_0},$$
  as $n_0, \ldots,n_{l+1}\cv\infty$. The corollary is
  proved.
\end{proof}

\begin{proof}[Proof of Theorem~\ref{thm:allbidegree}]
  Pick any set of integers $n_1 >m_1, \cdots, n_{d-1} >m_{d-1}$, and
  let $W^k \= \{ (c, a), \, P^{n_j}_{c, a}(c_{i_j}) = P^{m_j}_{c,
    a}(c_{i_j}) \text{ for all } 1\le j \le k \}$. Being defined by
  $k$ equations, the codimension of $W^k$ is greater than or equal to $k$.
  For a polynomial lying in $W^{d-1}$, all critical points are
  preperiodic, whence $P_{c , a}$ lies in the connected locus which is
  bounded by Proposition~\ref{p:compact}.  So $W^{d-1}$ is a finite
  set of points, i.e.  $\mathrm{codim}_\mathbb{C}\, W^{d-1} =d-1$. By
  induction we conclude that $\mathrm{codim}_\mathbb{C}\, W^k =k$ for
  all $k$.
  
  This proves (1) and (2). Theorem~\ref{thm:cvg} implies (3), because the
  variety $W$ in the statement of the theorem is quasi-projective. 
\end{proof}

\section{External rays}
\label{sec:external}

In this section we prove Theorems~\ref{thm:landing},~\ref{thm:supp} and~\ref{t:CE}.  
These rely on the combinatorial description
of polynomials in terms of external rays landing at critical points.
This technique was introduced and studied
in~\cite{Go1,BFH,kiwi-portrait}.  We describe the set $\mathsf{Cb}$ of
all these combinatorics in detail in Section~\ref{sec:comb} and show
that it is a compact and connected set endowed with a natural measure.
We do not claim originality here, but we hope that our presentation
will shed some light on the structure of this space.
Following~\cite{Go1} we construct in Section~\ref{sec:gold} a natural
map $\Phi_\mathrm{g}$ from $\mathsf{Cb} \times \mathbb{R}_+^*$ into
$\mathcal{P}_d$. The space $\mathsf{Cb} \times \mathbb{R}_+^*$ has a
natural structure of Riemann surface lamination, and $\Phi_\mathrm{g}$
provides an embedding into $\mathcal{P}_d\setminus \mathcal{C}$
preserving the lamination structure
(Proposition~\ref{p:def-goldberg}). We then describe the extension of
$\Phi_\mathrm{g}$ to a subset of $\mathsf{Cb}\times \{0\}$ of full
measure, and we state Kiwi's Continuity Theorem saying that
$\Phi_\mathrm{g}$ extends \emph{continuously} at Misiurewicz
combinatorics.  The restriction of $\Phi_\mathrm{g}$ to
$\mathsf{Cb}\times \{0\}$ defines a measurable 
``landing map'' from $\mathsf{Cb}$
into the boundary of the connectedness locus, which transports the natural
measure on the combinatorial space onto the bifurcation measure
(Section~\ref{sec:landing}).  In Section~\ref{sec:connected}, we
describe a connectedness property of a subset of the boundary of the
connectedness locus containing $\supp(\mu)$.  Section~\ref{sec:pf-end}
is devoted to the proof of Theorem~\ref{thm:supp}.
Finally, in Section~\ref{sec:CE} we prove Theorem \ref{t:CE}.


\subsection{The combinatorial space}\label{sec:comb}

We describe the set $\mathsf{Cb}$ of combinatorics of polynomials of
degree $d$ such that $g_P$ takes the same value on all critical
points. We define this space abstractly and study its geometry.


\subsubsection{The restricted combinatorial space}
We first look at the subset $\mathsf{Cb}_0$ of combinatorics of
polynomials of degree $d$ for which all critical points are marked and
simple.  In order to define it formally, we first need the following
\begin{defi}
  We denote by $\mathsf{S}$ the set of pairs $\{ \alpha , \alpha ' \}$
  contained in the circle $\mathbb{R}/\mathbb{Z} $, such that $d\alpha
  = d \alpha '$ and $\alpha \neq \alpha '$.
\end{defi}
Two finite and disjoint subsets $\theta_1, \theta_2 \subset
\mathbb{R}/\mathbb{Z} $ are said to be \emph{unlinked} if $\theta_2$
is included in a single connected component of $(\mathbb{R}/\mathbb{Z})
\setminus \theta_1$.
\begin{defi}
 We let $\mathsf{Cb}_0$ be the set of $(d-1)$-tuples $\Theta =
  (\theta_1,\cdots , \theta_{d-1})\in \mathsf{S}^{d-1}$ such that for
  all $i\neq j$, the two pairs $\theta_i$ and $\theta_j$ are disjoint
  and unlinked.
\end{defi}
We  spend the remaining part of this section describing the topology of
$\mathsf{S}$ and $\mathsf{Cb}_0$.

\smallskip

Recall that a smooth manifold of dimension $n$ is a \emph{translation
  manifold} if it admits an atlas for which the transition maps are
translations in $\mathbb{R}^n$.  Any translation manifold is endowed
with a natural  metric coming from the standard euclidean metric
on $\mathbb{R}^n$, and is naturally oriented. It is thus endowed with
a natural volume form, and a natural smooth measure.

\smallskip Let us explain how to define a translation structure  on
$\mathsf{Cb}$. The circle $\mathbb{R}/\mathbb{Z}$ admits a natural
structure of translation manifold, with a metric
$d_{\mathbb{R}/\mathbb{Z}}$.  Let $d_\mathrm{h}$ be the Hausdorff
metric on compact sets in $(\mathbb{R}/\mathbb{Z} ,
d_{\mathbb{R}/\mathbb{Z}} )$.  It endows $\mathsf{S}$ with a structure
of a compact metric space.  Denote by $\delta: \mathsf{S}\to \mathbb{R}_+^*$ the
map sending $\theta = \{ \alpha , \alpha '\}$ to the distance between
$\alpha $ and $\alpha '$ in $\mathbb{R}/\mathbb{Z} $.  This map is
continuous with values in the discrete set $ \{ 1/d, \cdots ,
[d/2]/d\}$ where $[d/2]$ is the integral part of $d/2$. One easily
checks that $\delta$ is surjective, and that each preimage is
connected.  Finally, for each $k\le [d/2]$, we have a surjective map
$\pi_k : \mathbb{R}/\mathbb{Z} \to \delta^{-1}\{ k \} \subset
\mathsf{S}$ sending $\alpha $ to $\{ \alpha , \alpha + k/d\}$. This
map is a $2$-to-$1$ covering map if $d$ is even and $k = d/2$, and
$1$-to-$1$ otherwise.  Moreover it commutes with any translation in
$\mathbb{R}/\mathbb{Z} $.  Hence every connected component of
$\mathsf{S}$ inherits by $\pi$ a translation structure coming from $\re/\zz$.

In short we conclude that:
\begin{lem}\label{l:sfS}
The space $\mathsf{S}$ has a natural structure of translation
manifold of dimension 1. It is compact  and has $[d/2]$
connected components, all of them homeomorphic to $\mathbb{R}/\mathbb{Z} $.
\end{lem}
We may now define a translation structure to 
 $\mathsf{Cb}_0$ by way of the following
\begin{prop}\label{p:trans}
  The set $\mathsf{Cb}_0$ is an open subset of $\mathsf{S}^{d-1}$, hence admits
  a natural structure of translation manifold.
  \end{prop}
\begin{proof}
  The subset $\mathsf{Cb}_0$ of $\mathsf{S}^{d-1}$ is defined by the
  two conditions $\theta_i \neq \theta_j$ and $\theta_i,\theta_j$ are
  unlinked for any $i\neq j$. Both conditions are clearly open for the
  Hausdorff topology on compact subsets. Whence $\mathsf{Cb}_0$ is open.
\end{proof}
\begin{rmk}
  In fact, it can be proven that $\mathsf{Cb}_0$ has finitely many
  connected components, each homeomorphic to $\mathbb{R}/\mathbb{Z} \times
  ]0,1[^{d-2}$. The set of connected components is in natural bijection
  with the set of finite simplicial trees having $d-1$ edges labelled with
  $\{ 1, \cdots, d-1\}$, and such that the set of branches at any
  vertex is oriented.
\end{rmk}
The translation manifold $\mathsf{S}^{d-1}$ is compact,  hence
endowed with a natural positive measure which is of finite
mass.  
\begin{defi} We define $\mu_{\mathsf{Cb}_0}$ to be the probability measure
proportional to the natural measure on $\mathsf{S}^{d-1}$ and whose
support is  precisely $\mathsf{Cb}_0$.
\end{defi}


\subsubsection{The full combinatorial space}
We now describe the set of combinatorics of \emph{all} polynomials of
degree $d$ (with all critical points marked).
\begin{defi}\label{def:cb}
  The set $\mathsf{Cb}$ is the collection of all $(d-1)$-tuples
  $(\theta_1, \cdots, \theta_{d-1})$ of finite sets in
  $\mathbb{R}/\mathbb{Z} $ satisfying the following four conditions:
\begin{itemize}
\itm for any fixed $i$, $\theta_i = \{ \alpha _1 , \cdots , \alpha
  _{k(i)}\}$ and  $d\, \alpha _j = d\, \alpha_1$ for all $j$;
\itm for any $i\neq j$, either $\theta_i \cap \theta_j = \emptyset$,
  or $\theta_i = \theta_j$;
\itm if $N$ is the total number of distinct $\theta_i$'s, then
$\mathrm{Card}\, \bigcup_i \theta_i 
  = d+N-1$;
\itm for any $i,j$ such that $\theta_i \cap \theta_j = \emptyset$,
  then $\theta_i$ and   $\theta_j$ are unlinked, that is 
$\theta_j$ is contained in a single connected component of
  $\mathbb{R}/\mathbb{Z} \setminus \theta_i$.
\end{itemize}
\end{defi}
A comment on the third item is in order. It may be rephrased as
follows: if $\theta_{i_1},\ldots, \theta_{i_N}$ is a maximal family of
disjoint sets in $(\theta_1, \cdots, \theta_{d-1})$, then $\sum
\mathrm{Card} (\theta_{i_j}-1)=d-1$. This  models the fact
that our polynomials have exactly $d-1$ critical points, counting
multiplicities. 

We now define a topology on $\mathsf{Cb}$ 
(see~\cite{kiwi-portrait}).  For any collection of open sets $O_1,
\cdots, O_{d-1}\subset \mathbb{R}/\mathbb{Z} $, define $U(O) = \{
\Theta =(\theta_i)\in \mathsf{Cb}, \, \text{s.t. } \theta_i \subset
O_i\}$.  Any intersection of such sets is of the same form. An open
set in $\mathsf{Cb}$ is by definition an arbitrary union of sets of
the form $U(O)$.
\begin{prop}[see also {\cite[Lemma 3.25]{kiwi-portrait}}]\label{p:connected}
  The set $\mathsf{Cb}$ is compact and path connected, and contains
  $\mathsf{Cb}_0$ as a dense and open subset.
\end{prop}
\begin{rmk}
  The set $\mathsf{Cb}$ can be stratified in a natural way such that
  each stratum is a translation manifold. 
\end{rmk}
\begin{defi}\label{d:meas}
  We let $\mu_\mathsf{Cb}$ be the unique measure which coincides with
  $\mu_{\mathsf{Cb}_0}$ on $\mathsf{Cb}_0$ and which does not charge
  $\mathsf{Cb}\setminus\mathsf{Cb}_0$.
\end{defi}
\begin{proof}
  In order to prove that $\mathsf{Cb}$ is path connected, it is
  sufficient to prove that any $\Theta= (\theta_i)$ can be joined by a
  continuous path to $\Theta_\star = (U_d , \cdots , U_d)$ where $
  U_d= \{ 0,1/d, \cdots, (d-1)/d\}$.  We prove this by induction on
  the number $N$ of distinct $\theta_i$'s. When this number equals
  $1$, all $\theta_i$'s are equal to the translate of $U_d$ by
  some element $\alpha \in \mathbb{R}/\mathbb{Z} $. Reducing the
  parameter of translation to $0$ gives us a path joining $\Theta$ to
  $\Theta_\star$. Now assume the claim has been proven for $N-1\ge 0$,
  and suppose $\Theta$ has $N$ distinct $\theta_i$'s. As before,
  translate all $\theta_i$ equal to $\theta_1$ while leaving the
  others fixed, until $\theta_1$ intersects another $\theta_j\neq
  \theta_1$. We conclude by using the inductive hypothesis.
  
  \smallskip The fact that $\mathsf{Cb}_0$ is open is clear. A
  descending induction on the number $N$ defined before shows that any
  $\Theta$ is the limit of a sequence of elements in $\mathsf{Cb}_0$.
  This shows the density of $\mathsf{Cb}_0$ in $\mathsf{Cb}$.

  \smallskip To prove that $\mathsf{Cb}$ is compact, we let
  $\overline{\mathsf{Cb}}_0$ be the closure of $\mathsf{Cb}_0$ in
  $\mathsf{S}^{d-1}$ (for the topology induced by the Hausdorff
  distance). This is a compact space. We claim that there exists a
  continuous surjective map $\pi: \overline{\mathsf{Cb}}_0 \to
  \mathsf{Cb}$, which yields the desired result.
  
  Define an equivalence relation $\sim$ on $\cup \theta_i$ as follows:
  $\alpha \sim \alpha'$ iff there exists a chain $\theta_{i_1},\ldots
  ,\theta_{i_k}$, with $\alpha\in \theta_{i_1}$, $\alpha'\in
  \theta_{i_k}$, and for every $j$, $\theta_{i_j}\cap
  \theta_{i_{j+1}}\neq \emptyset$. Define $\pi:(\theta_1, \cdots ,
  \theta_{d-1})\mapsto (\widetilde{\theta}_1, \cdots ,
  \widetilde{\theta}_{d-1})$, where $\widetilde{\theta_i}$ is the
  equivalence class of points contained in $\theta_i$.  By
  construction, $\pi$ is the identity map on $\mathsf{Cb}_0$.

 We  prove that $\pi(\Theta)$ belongs to $\mathsf{Cb}$. The
 first and second conditions of Definition~\ref{def:cb} are clearly
 satisfied. To check the third one,  notice that if $\theta_i$ and
 $\theta_j$ belong to $\overline{\mathsf{Cb}_0}$, $\theta_i\neq 
\theta_j$, and $\theta_i \cap
\theta_j\neq\emptyset$, then $\mathrm{Card}(\theta_i \cap
\theta_j )=1$, and use the comment after Definition~\ref{def:cb}.

It remains to prove that the $ \widetilde{\theta}_i$ are pairwise
unlinked. If, say,  $\widetilde{\theta}_1$ and $\widetilde{\theta}_2$
are linked, we get that in the chain ${\theta}_{i_1},\ldots
,\theta_{i_k}$ constituting $\widetilde{\theta}_1$, one element
 is linked
with $\widetilde{\theta}_2$. Applying the same reasoning to the chain
constituting $\widetilde{\theta}_2$ implies that two of the subsets 
in $\set{\theta_1, \ldots, \theta_{d-1}}$ are linked, a contradiction. 
  
  This shows that $\pi: \overline{\mathsf{Cb}}_0\to \mathsf{Cb}$ is
  well-defined. It is continuous by construction, hence its image is
  compact, and in particular closed. But this image contains
  $\mathsf{Cb}_0$ which is dense in $\mathsf{Cb}$ so $\pi$ is
  surjective.
\end{proof}

\subsubsection{The lamination structure on $\mathsf{Cb} \times \mathbb{R}_+^*$}
A locally compact topological space is said to be \emph{laminated by
  Riemann surfaces} if every point admits a neighborhood $U_i$ 
homeomorphic to a product $\mathbb{D}\times T_i$ where $\mathbb{D}$ is
the unit disk, $T_i$ is a compact topological space, and such that the
transition function from $U_i$ to $U_j$ are continuous and their
restriction to complex disks are holomorphic. A {\em plaque} is a set of the
form $\mathbb{D}\times \{ t_i\}$ in some chart. The leaf ${L}$
passing through a point $x$ is the smallest pathwise connected set containing
$x$ and such that if a plaque intersects ${L}$ then it is
completely included in it.  We refer to~\cite{ghys} for general facts
on laminations.

The abelian group $\mathbb{R}/\mathbb{Z} $ acts naturally on itself by
translation.  It hence induces a natural action on $\mathsf{Cb}$.  If
$\Theta = (\theta_i) \in \mathsf{Cb}$, with $\theta_i =\{ \alpha
_{ij}\}$, and $\alpha \in \mathbb{R}/\mathbb{Z} $, we write $ \Theta +
\alpha \= ( \{ \alpha +\alpha _{ij}\}_i )$.  Denote by $\mathbb{H}=
\{z\in \mathbb{C},\, \mathrm{Re}\, (z) >0\}$ the right half plane. We
then have a natural action $\mathbb{H}\times (\mathsf{Cb}\times
\mathbb{R}_+^*) \to (\mathsf{Cb}\times \mathbb{R}_+^*)$ defined by
$(s+it) \cdot (\Theta,r) \= (\Theta +rt, sr)$ --here the group
structure on $\mathbb{H}$ is given by $(s_1+it_1)\star(s_2+it_2)=
(s_1+it_1)s_2+ it_2$.
\begin{prop}\label{p:lamination}
  There exists a unique structure of lamination by Riemann surfaces on
  $\mathsf{Cb} \times \mathbb{R}_+^*$ such that for any fixed
  $(\Theta, r) \in \mathsf{Cb} \times \mathbb{R}_+^*$, the map
  $\mathbb{H} \ni u \mapsto u \cdot (\Theta,r)$ is holomorphic.
  Moreover, all complex leaves of the lamination are analytically
  diffeomorphic to the punctured unit disk $\mathbb{D}^*$.
\end{prop}

\begin{proof}
  Pick $(\Theta_\star,r_\star)\in\mathsf{Cb}\times \mathbb{R}_+^*$,
  with $\Theta_\star =(\theta_{i,\star})$ and define
  $\mathsf{T}_\star$ as the set 
$\{ \Theta =(\theta_i) \in \mathsf{Cb},\, \text{
    s.t. } M_d(\theta_1) =M_d( \theta_{1,\star}) \}$, where $M_d
  (\alpha) \= d\,\alpha$. This is a compact set.
  
  For each $i$, choose an open set $O_i\subset \mathbb{R}/\mathbb{Z} $
  containing $\theta_{i,\star}$ and such that any connected component
  of $O_i$ has length $<1/d^2$ and contains exactly one element of
  $\theta_{i,\star}$. 
  Define $U_\star \= U(O) = \{ \Theta = (\theta_i),\, \theta_i
  \subset O_i\}$. This an open set in $\mathsf{Cb}$. Now define the
  map $\pi_\star: U\times \mathbb{R}_+^* \to
  \mathsf{T}_\star\times\mathbb{R}/\mathbb{Z} \times \mathbb{R}_+^*$
  as follows. If $\Theta$ belongs to $U_\star$, then $\theta_1 \subset
  O_1$ so $M_d(\theta_1)$ is a point at distance $<1/d$ of
  $M_d(\theta_{1,\star})$. Thus there exists  a unique $\alpha
  \in\mathbb{R}/\mathbb{Z} $ at distance $<1/d$ of $1$ such that
  $M_d(\alpha + \theta_1) = M_d(\theta_{1,\star})$, i.e. $\alpha \=
  M_d(\theta_{1,\star}-\theta_1)/d$.  We set $\pi(\Theta,r) = ( \Theta
  +\alpha ,-\alpha,r)$. This map is clearly continuous, and injective.
  Its inverse is given by $(\Theta,\alpha ,r) \mapsto (\Theta +\alpha
  ,r)$, so the image of $\pi_\star$ is an open set in
  $\mathsf{T}_\star\times\mathbb{R}/\mathbb{Z} \times \mathbb{R}_+^*$.
  By postcomposing with the map $(\Theta, \alpha ,r ) \mapsto (\Theta
  , \exp(-2i \pi (\alpha + i r)))$, we get a map that we again denote
  by $\pi_\star$, which is defined on $U_\star$ with values in
  $\mathsf{T}_\star \times \mathbb{C}\setminus\overline{\mathbb{D}}$.
  By definition, $\pi_\star: U_\star \to \mathsf{T}_\star \times
  \mathbb{C}\setminus\overline{\mathbb{D}}$ is a chart for the
  lamination structure.

  We now check the compatibility of this collection of charts.  Choose
  two charts $\pi_\star,\pi_\bullet$ centered at
  $(\Theta_\star,r_\star),(\Theta_\bullet,r_\bullet)$, and suppose
  that their domains of definition have non trivial intersection.
  Then the composition $\pi_\bullet \circ \pi_\star^{-1}$ is of the
  form $(\Theta,\alpha ,r) \mapsto ( \Theta+\xi , \alpha +\xi, r)$,
  with $\xi\= M_d(\theta_{1,\bullet}-\theta_{1, \star})/d$.  For a
  fixed $(\Theta,r)$ in the transversal $\mathsf{T}_\star$, the
  composition of $(\alpha ,r) \mapsto \pi_\bullet \circ \pi_\star^{-1}
  (\Theta,\alpha ,r)$ with the projection onto $\mathbb{R}/\mathbb{Z}
  \times \mathbb{R}_+^*$ is a translation of angle $\xi$ in
  $\mathbb{R}/\mathbb{Z}$. Equivalently, it is a complex
  rotation of angle $\exp(2 i \pi \xi)$ in
  $\mathbb{C}\setminus\overline{\mathbb{D}}$, and is thus holomorphic.
  This proves that these charts patch together defining a structure
  of lamination by Riemann surfaces on $\mathsf{Cb}\times
  \mathbb{R}_+^*$.
  
  For fixed $(\Theta_\star,r_\star)$, the composition of $u=
  (s+it)\in\mathbb{H} \mapsto \pi_\star( u\cdot (\Theta,r))$ with the
  projection onto the last two factors induces a map $\mathbb{H} \to
  \mathbb{C}\setminus\overline{\mathbb{D}}$ which is equal to $u
  \mapsto \exp(2\pi r_\star u)$ and is clearly holomorphic. One checks
  that this map is surjective onto the  complex leaf of the
  lamination passing through $(\Theta_\star,r_\star)$.  The uniqueness
  of the structure of lamination follows from this remark.

  Finally suppose that $u\cdot (\Theta_\star,r_\star) = u' \cdot
  (\Theta_\star,r_\star)$ for $u = s+it$ and $u' =s'+it'$. Then
  $s=s'$, and ${rt}+\theta_{1,\star} = {irt'}+\theta_{1,\star}$ as
  sets. The latter condition is equivalent to $ t- t'$ being congruent
  to an integer $k$ depending only on the configuration of
  $\theta_{1,\star}$. The quotient of $\mathbb{H}$ by $ u \mapsto u+
  ik$ is the punctured unit disk, which concludes the proof.
\end{proof}


\subsection{The Goldberg map}\label{sec:gold}
We now explain how the sets $\mathsf{Cb}_0$ and $\mathsf{Cb}$ are
connected with the parameter space of polynomials.

\subsubsection{Definition}
For any $r>0$, let $\mathcal{G}(r)$ be the set of polynomials $P$ of
degree $d$ for which  all critical points $c$ of $P$ satisfy $g_P(c) = r$.
We shall see that $\mathsf{Cb}$ and $\mathcal{G}(r)$ are closely
related. We first recall some basic facts about B{\"o}ttcher
coordinates.  For any polynomial $P$ there exists a unique holomorphic
map $\varphi_P$, which is defined on $\{g_P > G(P)\}$,
tangent to the identity at infinity,  and such that $\varphi_P \circ P = \varphi_P^d$, 
We call $\varphi_P$ the B{\"o}ttcher map. It further satisfies  $g_P =
\log |\varphi_P|$ and  varies holomorphically with $P$.  For
$\alpha \in \mathbb{R}/\mathbb{Z}$, the path $r \mapsto
\varphi^{-1}_P(re^{2 i \pi \alpha})$ is called the external ray
associated to the angle $\alpha$. 
When $G(P)>0$, $\varphi^{-1}_P(re^{2
  i \pi \alpha})$ tends to a well defined limit point as $r\cv G(P)$. We
then say that $\alpha$ is an {\em external argument} of this point.
External rays coincide with gradient
lines of $g_P$. 

The basic proposition is the following.
\begin{prop}\label{p:def-goldberg}
  There exists a unique continuous map $\Phi_\mathrm{g}: \mathsf{Cb}
  \times \mathbb{R}_+^* \to \mathcal{P}_d$, 
  $\Phi_\mathrm{g}(\Theta, r) = (P(\Theta, r), c_i(\Theta,r))$ such
  that the following hold.
\begin{enumerate}
\item For each $i$, the set of external arguments the
  critical point of  $c_i$ is $\theta_i$, and $g_{P(\Theta
    ,r)}(c_i) =r$.
\item The map $\Phi_\mathrm{g}(\cdot, r)$ is a homeomorphism from
  $\mathsf{Cb}$ onto $\mathcal{G}(r)$.  Moreover, $\Phi_\mathrm{g}$
  restricts to a homeomorphism from $\mathsf{Cb}_0$ onto the subset of
  $\mathcal{G}(r)$ of polynomials, all critical points of which
   are simple.
\end{enumerate}
\end{prop}
\begin{proof}  
  Pick $(P,c) \in \mathcal{G}(r)$. Each critical point $c_i$ belongs
  to the closure of the domain of definition of $\varphi_P$, so we may
  look at the set $\theta_i \= \{ \alpha^i_1, \cdots, \alpha^i_k\}$ of
  external arguments of $c_i$.  We get a collection of finite sets
  $\Theta(P,c) = (\theta_i)$. It is not difficult to check that they
  satisfy all conditions of Definition~\ref{def:cb}. We thus obtain a
  continuous map $\Psi: \mathcal{G}(r) \mapsto \mathsf{Cb}\times
  \mathbb{R}_+^*$, $\Psi(P,c) = (\Theta(P,c), g_P(c))$.
  By~\cite[Proposition 3.8]{Go1}, $\Psi$ is surjective, and
  by~\cite[Lemma~3.22]{kiwi-portrait}, it is injective.  As
  $\mathcal{G}(r)$ is compact, $\Psi$ is a homeomorphism onto
  $\mathsf{Cb}\times \{ r \}$, and we denote by $\Phi_\mathrm{g} :
  \mathsf{Cb}\times \{ r \} \to \mathcal{G}(r)$ its inverse. By
  construction, (1) is satisfied.  Furthermore, for any $(P,c)\in
  \mathcal{G}(r)$, the polynomial $P$ has a simple critical point at
  $c_i$ iff $c_i$ has exactly two external arguments. This implies
  that the image of $\mathsf{Cb}_0$ by $\Phi_\mathrm{g}$ is the set of
  polynomials with only simple critical points, and proves (2).
\end{proof}  
\begin{prop}\label{p:lamin}
The restriction of $\Phi_\mathrm{g}$ to any  leaf of
  the underlying lamination of $\mathsf{Cb} \times \mathbb{R}_+^*$ is
  holomorphic.  In particular, $\Phi_\mathrm{g}$ induces an embedding of
  the natural lamination of  $\mathsf{Cb} \times \mathbb{R}_+^*
  $ into $\mathcal {P}_d\setminus \mathcal{C}$.
\end{prop}
\begin{rmk}
  For a fixed $(\Theta,r)$, the image under $\Phi_\mathrm{g}$ of the
  leaf $\mathbb{H}\cdot (\Theta,r)$ coincides with the wringing
  curve of $\Phi_\mathrm{g}(\Theta,r)$ as defined by Branner-Hubbard
  in~\cite{BrH1}.
\end{rmk}
The proposition is based on the following remark.  Suppose $P \in
\mathcal{G}(r)$, then $g_P(c_i) =r$, so that we have $g_P(P(c_i)) = d
r > G(P)$. Hence we may consider the holomorphic maps $\varphi_i(P,c)
\= \varphi_P \, (P(c_i))$ in a neighborhood of $P$. This defines a
holomorphic map $\varphi \= (\varphi_1,\cdots,\varphi_{d-1})$ in the
neighborhood of $\mathcal{G}(r)$ taking its values in
$\mathbb{C}\setminus \overline{\mathbb{D}}$. We have
 \begin{prop}\label{p:invert}
  The set of points where the differential of $\varphi$ is not locally
  invertible is a  complex hypersurface $H$, such that for all
  $r>0$,  $H \cap
  \mathcal{G}(r)$ has no interior points in $\mathcal{G}(r)$.
\end{prop}
\begin{proof}[Proof of Proposition~\ref{p:lamin}]
  Pick $(\Theta,r) \in \mathsf{Cb}$.  We claim that the restriction of
  $\Phi_\mathrm{g}$ to the complex leaf passing through $(\Theta,r)$
  is holomorphic. By Proposition~\ref{p:lamination}, this complex leaf
  is $\mathbb{H} \cdot (\Theta, r)$, and a direct computation shows
  that 
  $\mathbb{H} \ni u \mapsto (\varphi\circ \Phi_\mathrm{g})\, [u \cdot
  (\Theta, r)] = \exp(dr (u-1))\varphi (P)$ if $\Phi_\mathrm{g}
  (\Theta, r) = (P,c)$. This map is clearly holomorphic as a function
  of $u$. For
  any $(P,c)$ outside $H$, $\varphi$ gives local holomorphic
  coordinates, so $u\mapsto \Phi_\mathrm{g}(u \cdot
  (\Theta, r))$ is holomorphic for $(\Theta, r)$ in 
 the dense subset $\mathsf{Cb}_0\cap
  \Phi_\mathrm{g}^{-1}( \mathcal{P}_d \setminus H)$. 
 For general $(P,c)$,  we conclude by continuity.
\end{proof}

\begin{proof}[Proof of Proposition~\ref{p:invert}]
  First note that $\mathcal{G}(r) = \bigcap_i \{ \log |\varphi_i| =
  d\, r\}$ is a (possibly singular) real-analytic subset of
  $\mathcal{P}_d$.  Although we do not strictly need it, we also
  note that $\mathcal{G}(r)$ has pure real dimension $d-1$, as it is
  homeomorphic to $\mathsf{Cb}$ by $\Phi$. 
Consider  the
  continuous map $\psi$, sending $(\theta_1, \cdots,
  \theta_{d-1})\in\mathsf{Cb}$ to $\left[\exp(d\, r+ i M_d(\theta_1)),
    \cdots, \exp(d\, r + i M_d(\theta_{d-1})) \right]\in 
(e^{dr} \mathbb{S}^1)^{d-1}$ (with $M_d(\alpha )
  = d\times \alpha$).  By construction, we have $\varphi = \psi \circ
  \Phi^{-1}$. As $\psi$ is open, the restriction map
  $\varphi|_{\mathcal{G}(r)}$ is open. Sard's Theorem implies that its
  differential has generically maximal rank.
 In other words, if $(P,c)$
  belongs to some full measure subset  in  $\mathcal{G}(r)$, the differential of
  $\varphi$ restricted to the tangent space of $\mathcal{G}(r)$ has
  (real) rank equal to the (real) dimension of
  $\varphi(\mathcal{G}(r)) = (e^{dr} \mathbb{S}^1)^{d-1}$ which is $d-1$.
  Moreover, the image of $d\varphi_{P,c}$ contains the tangent space
  at $(e^{dr} \mathbb{S}^1)^{d-1}$ which is totally real.  But
  $d\varphi_{P,c}$ is a complex linear map, so its image necessarily
  contains a complex vector space of dimension $d-1$. This proves that
  $\varphi$ is locally invertible at generic points of
  $\mathcal{G}(r)$. In particular $H = \{ \det d\varphi = 0 \}$ is a
  complex hypersurface whose intersection with 
  $\mathcal{G}(r)$ has no interior points, and $d\varphi$ is invertible outside $H$.
   \end{proof}

\subsubsection{Kiwi's continuity property for Misiurewicz combinatorics}
Fix $\Theta \in \mathsf{Cb}$.  The 
{\em (stretching) ray} associated  to $\Theta$ is the set
$\set{ P_r= P(\Theta,r), r>0}$. When $r\cv0$, then $G(P_r)$ converges to $0$ so that any
polynomial in the cluster set of $\{P_r\}$ belongs to the
connectedness locus. It is a very delicate problem to describe this
cluster set in general. We say that a ray {\em lands} if this cluster
set is a single point.  
\begin{defi}
  A combinatorics $\Theta = (\theta_i) $ is said to be of
  \emph{Misiurewicz type}, if any $\alpha \in \bigcup\theta_i$ is strictly preperiodic
  under the map $z \mapsto d\, z$. We denote by $\mathsf{Cb}_\mathrm{mis}$
  the set of Misiurewicz combinatorics.
\end{defi}
\begin{prop}\label{p:dense}
  The set $\mathsf{Cb}_\mathrm{mis}$ is dense in $\mathsf{Cb}$.
\end{prop}
\begin{proof}
  The set of periodic orbits of $z\mapsto d\, z$ is dense in
  $\mathbb{R}/\mathbb{Z} $.  Pick any finite set $\theta= \{ \alpha
  _1, \cdots, \alpha _j\}\subset \mathbb{R}/\mathbb{Z} $ such that
  $\alpha \= d\, \alpha _1= \cdots =d\, \alpha _j$. One can then find
  a periodic point $\alpha _\star$ arbitrarily close to $\alpha $
  whose orbit does not intersect $\theta$. Let $\theta'= \theta+
  \unsur{d} (\alpha _\star-\alpha)$. This is a finite set very close
  to $\theta$, and strictly preperiodic.
  
  Now let $\Theta = (\theta_i)\in\mathsf{Cb}$.  For each $i$,
  consider the set $I_i$ of all indices $j$ such that $\theta_j =
  \theta_i$. The preceding argument shows that we may translate all
  $\{\theta_j\}_{j\in I_i}$ at the same time so that they become
  strictly preperiodic. Doing the same for each $i$, we get a
  combinatorics $\Theta'$ which is of Misiurewicz type and arbitrarily
  close to $\Theta$.
\end{proof}
We now state without proof the following deep continuity result, which
is a combination of the results of \cite{kiwi-portrait} and
\cite{BFH}, 
see~\cite[Corollary~5.3]{kiwi-portrait}.
\begin{thm}\label{thm:cont-landing}
  The map $\Phi_\mathrm{g}$ extends continuously to
  $\mathsf{Cb}_\mathrm{mis} \times \{ 0 \}$. The extended map
  $\Phi_\mathrm{g}$ induces a bijection from $\mathsf{Cb}_\mathrm{mis}
  \times \{ 0 \} $ onto the subset  of Misiurewicz
  polynomials.
\end{thm}
Note that the continuity statement is particularly strong. It  means that
for any $\Theta \in \mathsf{Cb}_\mathrm{mis}$, the sequence $P(\Theta,
r)$ \emph{converges} when $r\cv 0$, and that the polynomial is
Misiurewicz. It also means that for any sequence $(\Theta_n,r_n) \cv
(\Theta,0)$,  $P(\Theta_n,r_n)$ converges to $P(\Theta,0)$.


\subsubsection{Measurable landing}
We now prove a statement somewhat dual to the previous theorem.  Note
that this result was already observed in~\cite[Theorem~B]{BMS}.
\begin{prop}\label{p:landing}
  For any $\Theta$ and for almost every $t>0$, the map
  $\Phi_\mathrm{g}\left((s+it)\cdot (\Theta,1) \right)$ has a limit
  when $s\cv0$.  In particular, for $\mu_\mathsf{Cb}$-almost every
  $\Theta \in\mathsf{Cb}$, the map $ r\mapsto \Phi_\mathrm{g}( \Theta
  , r)$ admits a limit when $r\cv 0$.
\end{prop}
\begin{proof}
  Let us prove the first statement. Fix $\Theta\in\mathsf{Cb}$.  By
  Proposition~\ref{p:def-goldberg}, the map $\mathbb{H} \ni u =
  (s+it)\mapsto \Phi_\mathrm{g}\left((s+it)\cdot (\Theta,1) \right) =
  \Phi_\mathrm{g}(\Theta+t,s)\in\mathcal{P}_d$ is holomorphic. For
   simplicity, we lift this map to the ramified cover $\pi:
  \mathbb{C}^{d-1} \to \mathcal{P}_d$ that we used in Section~\ref{sec:param}.
  We get a map $\varphi: \mathbb{H} \to \mathbb{C}^{d-1}$ such that $G
  (\varphi(u)) = \mathrm{Re}\,(u) \times G
  \left(\Phi_\mathrm{g}(\Theta,1)\right)$.  In particular, the
  restriction of $\varphi$ to any square $\{0< \mathrm{Re}\,(u)<1,\,
  |\mathrm{Im}\,(u)|\le R \}$ with $R>0$ is a holomorphic and bounded
  function. By Fatou's theorem (see e.g. \cite[Lemma~15.1]{mibook})
  $\lim_{s\cv0} \varphi(s+it)$ exists for almost every $|t|\le R$. We
  conclude by letting $R\cv\infty$.
  
  For the second statement, let $\mathcal{B}$ be the set of
  $\Theta\in\mathsf{Cb}$ such that $\Phi_\mathrm{g}( \Theta , r)$
  admits a limit when $r\cv 0$.  We have proved that for any $\Theta$,
  the set of $t>0$ for which $ \Theta +t$ belongs to $\mathcal{B}$ has
  full (Lebesgue) measure in $\mathbb{R}$. Now recall that
  $\mu_\mathsf{Cb}$ puts full measure on $\mathsf{Cb}_0$, that
  $\mathsf{Cb}_0$ is an open set of $(\mathbb{R}/\mathbb{Z})^{d-1}$
  and that the measure $\mu_{\mathsf{Cb}}$ is the restriction of the
  natural Haar measure on $(\mathbb{R}/\mathbb{Z})^{d-1}$. 
   By Fubini's Theorem, we conclude that $\mu_r(\mathcal{B}) = 1$.
\end{proof}


\subsection{Landing measure}\label{sec:landing}
Our main result is
\begin{thm}\label{thm:landingr}
  For any $r>0$, the Monge-Amp{\`e}re measure associated to $\max \{ G , r
  \}$ is equal to the image of $\mu_\mathsf{Cb}$ under the map $\Theta
  \mapsto \Phi_\mathrm{g}(\Theta ,r)$.
\end{thm}
As a corollary, we can give a proof of Theorem~\ref{thm:landing}
 stated  in the introduction.
\begin{proof}[Proof of Theorem~\ref{thm:landing}]
  By Proposition~\ref{p:landing}, the limit $e(\Theta) \= \lim_{r\to0}
  \Phi_\mathrm{g}(\Theta,r)$ exists for $\mu_\mathsf{Cb}$-almost every
  $\Theta$. This yields a measurable map $e: \mathsf{Cb} \to
  \fr \mathcal{C}$. Now for any $r>0$,
  $\Phi_\mathrm{g}(\cdot,r)_*\,\mu_\mathsf{Cb} = (dd^c)^{d-1} \max \{
  G , r\}$ by Theorem~\ref{thm:landingr}. As $\max\{ G,r\}$ decreases
  to the continuous psh function $G$ when $r$ decreases to $0$, we
  have $\lim_{r\cv0} (dd^c)^{d-1} \max \{ G , r\}= \mu_\mathrm{bif}$,
  see~\cite{BT}. On the other hand $\Phi_\mathrm{g}(\cdot,r)$ converges
  measurably to $e$. Whence $e_* \mu_\mathsf{Cb} = \mu_\mathrm{bif}$.
\end{proof}

\begin{proof}[Proof of Theorem~\ref{thm:landingr}]
Look at the following commutative diagram
$$
\xymatrix{
\mathsf{Cb} \ar[d]_{\Phi} \ar[dr]^{\psi} \\
\mathcal{G}(r)\ar[r]_(.4){\varphi} 
&
(e^{dr} \mathbb{S}^1 )^{d-1}
}
$$
where $\Phi = \Phi_\mathrm{g}(\cdot,r)$, $\varphi (P,c)= (\varphi_P
\, P(c_i))$ and $\psi(\theta_1, \cdots, \theta_{d-1})\= (e^{d\, r+ i
  M_d(\theta_1)}, \cdots, e^{d\, r + i M_d(\theta_{d-1})})$ as in the
proof of Proposition~\ref{p:invert}. By construction, $\varphi \circ
\Phi = \psi$, and recall that by Proposition~\ref{p:def-goldberg}, the
map $\Phi$ is a homeomorphism.  Let $\mu_r = (dd^c)^{d-1} \max \{ G ,
r \}$. We  show that $\Phi^* \mu_r= \mu_\mathsf{Cb}$ (where
$\Phi^*$ stands for $(\Phi^{-1})_*$).  Denote by $d\lambda$ the Haar
measure on the $(d-1)$ real dimensional torus $(e^{dr} \mathbb{S}^1
)^{d-1}$.
\begin{lem}\label{l:compu-phi}
  There exists an open set $G\subset\Phi(\mathsf{Cb}_0)$, such
  that $\mu_r(G) = \mu_\mathsf{Cb}( \Phi^{-1}(G))=1$, and every point
  in $G$ admits an open neighborhood $U$ on which $\varphi_* \mu_r|_U$
  coincides with $d^{1-d}\times d\lambda|_{\varphi(U)}$.
\end{lem}
Pick $(P,c)\in G$.  Take $U$ as in the lemma so that $\varphi_*
\mu_r|_U = d^{1-d}\times d\lambda|_{\varphi(U)}$. As $G\subset
\Phi(\mathsf{Cb}_0)$, the polynomial $P$ has only simple critical
points, and we may assume that this property is satisfied all over
$U$.  Let $V$ be the inverse image of $U\cap\mathcal{G}(r)$ by $\Phi$.
By Proposition~\ref{p:def-goldberg}, it is contained in
$\mathsf{Cb}_0$.  Now note that $\psi: \mathsf{Cb}_0\to (e^{dr}
\mathbb{S}^1) ^{d-1}$ preserves the structure of translation manifolds
of both spaces, so $\psi_* \mu_\mathsf{Cb}$ is proportional to
$d\lambda$.  In particular, $\psi_* \mu_\mathsf{Cb}|_V =t \times
d\lambda$ where $t^{-1}= \lambda (\psi(\mathsf{Cb}_0))$ does not
depend on $V$. We infer that $\Phi^*\mu_r = t'\times \mu_\mathsf{Cb}$
on $V$ hence on $\Phi^{-1}(G)$ with $t' = d^{1-d}/t$. But $G$ has full
$\mu_r$-measure, and $\Phi^{-1}(G)$ has full $
\mu_\mathsf{Cb}$-measure, so $\Phi^*\mu_r = t'\times \mu_\mathsf{Cb}$.
Both $\mu_r$ and $\mu_\mathsf{Cb}$ being probability measures,
$\Phi^*\mu_r =\mu_\mathsf{Cb}$, as required.
 \end{proof}

\begin{proof}[Proof of Lemma~\ref{l:compu-phi}]
  The same proof as for Lemma~\ref{l:wedge} gives in a neighborhood
  $\mathcal{N}$ of $\mathcal{G}(r)$,
\begin{eqnarray*}
\mu_r = (dd^c)^{d-1}\max \{ G ,r \} 
& =& 
dd^c \max\{ g_P(c_1),r \} \wedge \cdots \wedge dd^c \max\{ g_P(c_{d-1}),r \} 
\\
&=&
 dd^c \max\{ d^{-1}\log|\varphi_1|,r\} \wedge \cdots 
\wedge  dd^c \max \{d^{-1}\log|\varphi_{d-1}|,r\}
 \end{eqnarray*}
 Apply Proposition~\ref{p:invert}, and set $G \= \{ (P,c) \in
 \mathcal{N}\setminus H, \, \text{s.t. } P \text{ has only simple critical
   points}\}$. For any point in $G$, there exists an open neighborhood
 such that the mapping $\varphi$ is a holomorphic diffeomorphism from
 $U$ onto its image, so
\begin{eqnarray*}
\mu_r|_U = 
& =& 
\varphi|_U^* \left( dd^c \max\{ d^{-1}\log|z_1|,r\} \wedge
 \cdots \wedge  dd^c \max \{d^{-1}\log|z_{d-1}|,r\}\right)
\\
& = & 
d^{1-d}\, \varphi|_U^*\left(  dd^c \max\{ \log|z_1|,r^d\} \wedge 
\cdots \wedge  dd^c \max \{\log|z_{d-1}|,r^d\}\right)
 \end{eqnarray*}
 Now for any real number $\rho>0$, the measure $dd^c \max\{
 \log|z_1|,\rho\} \wedge \cdots \wedge dd^c \max
 \{\log|z_{d-1}|,\rho\}$ in $\mathbb{C}^{d-1}$ is  a
 probability measure supported on $(e^\rho\mathbb{S}^1 )^{d-1}$ and 
 invariant under the subgroups of rotations
 $\simeq(\mathbb{S}^1)^{d-1}$.  So it equals the Haar measure
 $d\lambda$. 
 
 In order to conclude, we need to show that $\mu_r(G) =
 \mu_\mathsf{Cb}( \Phi^{-1}(G))=1$.  First $\mu_r$ does not charge
 complex analytic sets, so $\mu_r (G) = \mu_r (\mathcal{N}) =1$. Now
 $\Phi^{-1}(H)\subset\psi^{-1} \varphi(H)$, and $\varphi(H)$ is a complex
 analytic set in $\mathbb{C}^{d-1}$. Its intersection with the torus
 $(e^{dr}\mathbb{S}^1)^{d-1}$ is hence of zero Haar measure.
Finally $ \mu_\mathsf{Cb} \Phi^{-1}(H) \le \mu_\mathsf{Cb} \psi^{-1} \varphi(H)
= (\psi_*\mu_\mathsf{Cb}) \varphi(H) = t \times \lambda \varphi(H) =0$.
\end{proof}


\subsection{A connectedness property}\label{sec:connected}
The next proposition gives a more accurate picture of the geometry of
the boundary of the connectedness locus. When $d=2$, this is the
statement that the Mandelbrot set is connected. Let $\mathcal{L}$ be
the subset of the shift locus consisting of the polynomials such that
$g_P(c_1)=\cdots =g_P(c_{d-1})$, that is, $\mathcal{L} = \Phi_{\rm g}
(\mathsf{Cb}\times\re^*_+)$. Note that $\supp(\mu_{\rm bif}) \subset \mathcal{C}\cap
\overline{\mathcal{L}}$.  
\begin{prop}\label{p:connex} 
The set $\mathcal{C}\cap \overline{\mathcal{L}}$ is connected. 
\end{prop}  

\begin{proof} For $r>0$, recall the set $\mathcal{G}(r)$ of polynomials for which
  all critical points satisfy $g_P(c)=r$. It is homeomorphic to
  $\mathsf{Cb}$, hence connected.  Let now $\mathcal{G}(\leq
  \hspace{-.2em}r)$ be the union of $\mathcal{G}(s)$ for $0<s\leq r$.
  Being homeomorphic to $\mathsf{Cb}\times ]0,r]$ it is connected.
  Hence $\overline{\mathcal{G}(\leq\hspace{-.2em} r)}$ is compact and
  connected.  Finally
  $$\bigcap_{r>0}\overline{\mathcal{G}(\leq \hspace{-.2em}r)} = \overline{\el}\cap \mathcal{C}$$ is connected.
  \end{proof}


\subsection{Proof of Theorem~\ref{thm:supp}}\label{sec:pf-end}
We prove that any Misiurewicz parameter $(P,c)$ lies in the
support of $\mu_\mathrm{bif}$. Apply Theorem~\ref{thm:cont-landing},
and pick $\Theta_*\in \mathsf{Cb}_{\rm mis}$ such that
$\Phi_\mathrm{g}(\Theta_*,0) = (P,c)$. Apply now
Theorem~\ref{thm:landingr}: $e_*\mu_\mathsf{Cb} = \mu_\mathrm{bif}$.
In particular $e(\Theta)$ is well defined and belongs to the support of
$\mu_\mathrm{bif}$ for $\mu_\mathsf{Cb}$-almost every $\Theta$. Now,
$\mu_{\mathsf{Cb}}$ has full support in 
 $\mathsf{Cb}$. Thus we can find a sequence $\Theta_k$
converging to $\Theta_*$ and such that $e(\Theta_k)$ lies in the
support of $\mu_\mathrm{bif}$. For each $k$, pick $\e_k>0$ small
enough such that the distance between $\Phi_\mathrm{g}(\Theta_k,\e_k)$
and $e(\Theta_k)$ is less than $1/k$.  In particular, the distance
between $\Phi_\mathrm{g}(\Theta_k,\e_k)$ and $\mathrm{Supp}\,
\mu_\mathrm{bif}$ tends to $0$ when $k\cv\infty$.  By the Continuity
Theorem~\ref{thm:cont-landing}, we have
$\Phi_\mathrm{g}(\Theta_k,\e_k)\cv (P,c)$.  Whence $(P,c)$ is in the
support of the bifurcation measure.

The reverse inclusion was already proved in  Corollary~\ref{cor:mis}.
\hfill $\square$


\subsection{Topological Collet-Eckmann property}\label{sec:CE}

In this section we give the proof of Theorem \ref{t:CE} (and
Corollary~\ref{cor:CE}). This is an adaptation of the proof of 
Smirnov~\cite{smirnov}. Although Smirnov worked only with unicritical
polynomials, most of his arguments remains valid in the multi-critical
case as well, see~\cite[Remark 1]{smirnov} and~\cite[p. 348]{BMS}. The other
two ingredients in the proof are the work of Kiwi on the combinatorics
of multi-critical polynomials,~\cite{kiwi-real} and our landing
Theorem~\ref{thm:landing}.

\begin{proof}[Proof of Theorem \ref{t:CE}]
  Recall that the TCE condition reads as follows: the polynomial $P$
  satisfies the TCE condition if for some $A\geq1$ there exists
  constants $M>1$ and $r>0$ such that for every $x\in J_P$ there is an
  increasing sequence $(n_j)$ with $n_j\leq Aj$ such that for every
  $j$,
\begin{equation}\label{e:CE}\# \set{i, \ 0\leq i<n_j,\
  \mathrm{Comp}_{f^i(x)}f^{-(n_j-i)}B(f^{n_j} (x), r)\cap\mathrm{Crit}
\neq \emptyset  }\leq M,\end{equation}
 where $\mathrm{Crit}$ denotes the critical set, and $
\mathrm{Comp}_xX$ is the connected component of the set $X$ containing
  $x$. More precisely we refer to \eqref{e:CE} as the TCE$(M, A, r)$
  condition. 
When $P$ is a polynomial with marked critical points
  $(c_1, \ldots, c_{d-1})$, we say that $P$ satisfies
  the TCE$_k(M, A, r)$ condition if $\eqref{e:CE}$ holds with $c_k$
  instead of $\mathrm{Crit}$. 
It is clear that  if $P$ fails the TCE$(M,A,r)$ condition, it fails
  the TCE$_k(\frac{M}{d-1},A,r)$ condition for some $k$.

We now introduce the following subset $\mathsf{Cb}_1 \subset \mathsf{Cb}$
  of combinatorics $\Theta$ satisfying the following three conditions:
\begin{itemize}
\itm the ray $\{ \Phi_\mathrm{g}(\Theta, r)\} _{r>0}$ lands at a
  polynomial with marked critical points $e(\Theta) \= (P,c)$.
\itm $P$ has only simple critical points (i.e.  $\Theta\in
  \mathsf{Cb}_0$), and none of them is preperiodic.
\itm
$P$ has only repelling cycles; in particular $K_P = J_P$.
\end{itemize}
Note that the first condition is $\mu_\mathsf{Cb}$-generic by
Theorem~\ref{thm:landing}; the second is also generic as
$\mu_\mathrm{bif}$ does not charge hypersurfaces.
Finally~\cite[Theorem~1]{kiwi-portrait}, and~\cite[Lemma~5]{BMS},
implies the genericity of the last condition too.  So $\mathsf{Cb}_1$
has \emph{full measure} in $\mathsf{Cb}$.

We will prove that the set of combinatorics $\Theta\in \mathsf{Cb}_1$
for which the associated polynomial $P$ violates the
TCE$_k(\frac{M}{d-1},A,r)$ condition for some $A$ and every $M,r$ has
zero $\mu_{\rm bif}$ measure. Without loss of generality, assume
$k=1$.

\medskip

We use external rays in dynamical plane to understand the
recurrence property of critical points on the Julia set. As external
rays do not land in general, we are led to work with {\em fibers},
originally introduced by Schleicher \cite{sch}.

Recall that external rays with rational angles always land at
preperiodic points \cite{DH}.  By definition, two points $\xi,
\zeta\in J_P$ {\em do not belong to the same fiber} if there exist
external rays $R_s$ and $R_t$, with rational angles, landing at a
common point $z$, such that $R_s\cup\set{z}\cup R_t$ separates $\xi$
from $\zeta$. We say that $\theta$ is an {\em external argument} of
$\zeta\in J_P$ if $\overline{R_\theta}\cap J_P\subset {\rm
  Fiber}(\zeta)$. It follows from the work of Kiwi \cite[Theorem
3]{kiwi-real} that when $P$ has only repelling cycles, every $\zeta\in
J_P$ has a nonempty finite set of external arguments.

The following lemma is the classical connection between external arguments in 
dynamical and parameter spaces.

\begin{lem} \label{lem:crit_fiber} 
  Let $\Theta =(\theta_1, \cdots, \theta_{d-1}) \in \mathsf{Cb}_1$,
  with $\theta_1= \set{\alpha, \alpha'}\in\mathsf{S}$. Write
  $e(\Theta) = (P, c_1, \cdots, c_{d-1})$.  Then
  $\overline{R_\alpha}\cap J_P$ (resp.  $\overline{R_{\alpha'}}\cap
  J_P$) is contained in the fiber of $c_1$, that is, $\alpha$ and
  $\alpha'$ are external arguments of $c_1$.
\end{lem}

\begin{proof} Assume that some accumulation point $\zeta$ of $R_\alpha$ does
   not belong to the fiber of $c_1$. Then there exists a pair of
   rational rays $R_s$ and $R_t$ landing at $z\in J_P$
separating $\zeta$ from $c_1$. The point $z$ is necessarily
   prerepelling and not precritical, because $P$ has no preperiodic
   critical point and all cycles repelling.
Now, if we perturb $P$ in the parameter
   space of polynomials, $z$ admits a continuation as a prerepelling
   point and for the perturbed map, the rays $R_s$ and $R_t$  still
   land at the continuation of $z$ (see e.g. \cite[Lemma
   5.2]{kiwi-portrait}). But in the parameter ray associated to
   $\Theta$, $\theta_1=\set{\alpha, \alpha'}$ is the set of 
 external arguments of $c_1$ so $c_1$ and $R_\alpha$ are in the same
   connected component of $\cc\setminus(R_s\cup R_t\cup\set{z})$, a
   contradiction. 
\end{proof}
\begin{lem} \label{lem:disconnect}
  Let $I_0$ and $I_1$ be the two connected components of
  $\re/\zz\setminus \{ \alpha , \alpha'\}$. For $\e = 0,1$, define
  $J_\e$ to be the set of points in $J_P\setminus\mathrm{Fiber}(c_1)$
  having an external argument in $I_\e$.
  
  Then the three sets $J_0$, $J_1$ and $\mathrm{Fiber}(c_1)$ form a
  partition of $J_P$, and all iterates of $c_1$ lie in $J_0\cup J_1$.
\end{lem}

\begin{proof} By \cite[Proposition 3.15]{kiwi-real}, 
  $\mathrm{Fiber}(c_1)$ is compact, connected and full. Hence
  ${R_\alpha}\cup\mathrm{Fiber}(c_1)\cup {R_{\alpha'}}$ separates the
  plane into two connected components, $U_0$ and $U_1$. Define $J_\e
  \= U_\e \cap J_P$. By construction an external ray with argument in
  $I_\e$ is contained in $U_\e$. So any external angle of a point in
  $J_\e$ belongs to $I_\e$.
  
  Finally, by \cite[Corollary 2.15]{kiwi-real}, the image of a fiber
  is a fiber, and, by \cite[Lemma 4.4]{kiwi-real}, if
  $\mathrm{Fiber}(c_1)$ is periodic, then $c_1$ itself is periodic. By
  assumption $c_1$ is not preperiodic, so all iterates of $c_1$ lie in
  $J_0\cup J_1$.
\end{proof}
We now follow Smirnov's proof. Pick $\Theta\in\mathsf{Cb}_1$, and
write $e(\Theta) = (P,c)$. By the previous lemma, $c_1$ has a well
defined itinerary in the space $\Sigma_2 \= \set{0,1}^{\mathbb N^*}$.
Call this itinerary the {\em kneading sequence} of $c_1$.  In this way, we get a
map $\kappa : \mathsf{Cb}_1 \to \Sigma_2$, characterized by the
condition that $P^n(c_1) \in U_{\e_n}$ for all $n$ with $\kappa(\Theta)=
(\e_1, \e_2, \cdots)$.

By repeating \cite[Section 2]{smirnov}, we get that if $P$ fails the
TCE$_1(M,A,r)$ condition for some $A$ and every $M,r$, then
$\kappa(P)$ is {\em Strongly Recurrent}.  We refer to \cite{smirnov} for a
precise definition of this condition. The key  fact is that
the set ${\rm SR}$ of stronly recurrent sequences has zero Hausdorff
dimension in $\Sigma_2$. Here we endow $\Sigma_2$ with its usual
2-adic metric, that is, $d(x,y)= 2^{-n}$, where $n$ is the smallest
integer such that $x_n\neq y_n$.  Notice that the Hausdorff dimension
of $(\Sigma_2, d)$ is 1. We now rely on the following lemma which is
the analogue of \cite[Proposition 1]{smirnov}.
\begin{lem}\label{lem:HD}
  The Hausdorff dimension of the set $\kappa^{-1}(SR)\subset
  \mathsf{Cb}_1$ is not greater than $(d-2)+ \log(d-1)/\log d <d-1 = \dim
  (\mathsf{Cb}_1)$.
\end{lem} 
The lemma implies that $\kappa^{-1}(SR)\cap\mathsf{Cb}_1$ has
zero $\mu_{\mathsf{Cb}}$-measure, and the proof of the theorem is complete.
\end{proof}
\begin{proof}[Proof of Lemma~\ref{lem:HD}]
  First note that $\kappa$ is the composition of the projection $\pi:
  \mathsf{Cb}_1 \to \mathsf{S}$ onto the first factor $\Theta \mapsto
  \theta_1$, together with the map $K:\mathsf{S} \to \Sigma_2$ defined
  as follows.  If $\theta \= \{ \alpha, \alpha'\} \in \mathsf{S}$,
  then set the $n^{\rm th}$ term of $K(\theta)$ to be $\e$ iff $d^n\alpha=
  d^n\alpha'\in I_\e $. We will prove that the Hausdorff dimension of
  $K^{-1}(SR)$ is not greater  than $\log(d-1)/\log d$, which implies the
  lemma.
  
  It is enough to restrict to one connected component $\mathsf{S}_0$
  of $\mathsf{S}$.  For simplicity we assume that the distance between
  $\alpha$ and $\alpha'$ is $k/d$ with $k<d/2$, so the component
  $\mathsf{S}_0$ is parameterized by $\alpha\in \re/\zz$, and
  $\alpha'=\alpha+k/d$.  The remaining case $k=d/2$ is left to the
  reader.
  
  To get the dimension estimate, we prove that if $C_n$ is any
  cylinder of depth $n$ in $\Sigma_2$, then $K^{-1}(C_n)$ consists of
  at most $C^{st} \times n\,(d-k)^n$ intervals of length at most
  $d^{-n}$ (here $d-k\geq 2$ since $d\geq 3$).  This easily implies
  that the dimension of $K^{-1}(SR)\cap \mathsf{S}_0$ is at most
  $\log(d-k)/\log d$ (compare with \cite[Section 8]{bs}).
  
  First it is clear that the $n^{\rm th}$ digit of the kneading
  sequence cannot remain constant if $\alpha$ gets increased by an
  amount of $1/d^n$: indeed $d^n\alpha$ turns once around $\re/\zz$ so
  it hits $\alpha$ or $\alpha'$ which are almost fixed and at distance
  $k/d$.  The estimate of the number $N(n)$ of intervals goes by
  induction on $n$.  Let $M_d$ be the multiplication-by-$d$ map on $\re/\zz$.
  Recall that $\re/\zz = I_0\cup I_1\cup\set{\alpha, \alpha'}$, and
  let $\ell_k/d$ be the length of $I_k$, $\ell_k\in \set{k,d-k}$.  Let
  $I$ be any interval in $\re/\zz$ of length $<1/d$. There are two
  cases. Either $d\alpha \notin I$ and $M_d^{-1}(I)\cap I_k$
  consists of $\ell_k$ intervals, or $d\alpha \in I$ and
  $M_d^{-1}(I)\cap I_k$ consists of $\ell_k+1$ intervals. Of course
  the latter case occurs for at most one of the $N(n)$ intervals.  We
  get that $N(n+1)\leq (d-k)(N(n)-1) + (d-k+1)= (d-k)N(n) +1$, whence
  $N(n) \le C^{st}\times n\,(d-k)^n$.
\end{proof}

\begin{proof}[Proof of Corollary \ref{cor:CE}]
  We already noticed in the proof of the previous theorem that the
  first assertion is a consequence of~\cite[Theorem~1]{kiwi-portrait},
  and~\cite[Lemma~5]{BMS}. The third one is a combination
  of~\cite{gsmirnov,prz2} for the local connectivity statement,
  and~\cite{pr} for the statement on Hausdorff dimension.
 Lastly, the second item follows from the proof
  of the theorem. Indeed, if $P$ satisfies the TCE property, then $J_P$ is
  locally connected, so the landing of external rays defines a
  continuous map $\re/\zz\cv J_P$, semiconjugating $P$ to
  multiplication by $d$. If $P=e(\Theta)$, the angles in $\theta_i$
  are external arguments of $c_i$, and for generic $\Theta$, they have
  dense orbit on $\re/\zz$ under multiplication by $d$.  
\end{proof}

\section{The space of rational maps}\label{sec:rational}
In this last section we indicate how some of the  results of Section
\ref{sec:poly} can
be extended to the space of rational maps of degree $d$ with critical
points marked. Since the methods and results are quite similar, the
exposition is somewhat sketchy.  

Let $\mathrm{Rat}_d$ be the space of all rational maps of degree
$d\ge2$, with $(2d-2)$ marked critical points, modulo conjugation
by M{\"o}bius transformations.  A point in this space is a $(2d-1)$-tuple
$(R, c_i)$ such that $R$ is a rational map of degree $d$, and $\{
c_i\}$ is the set of critical points of $R$ (counted with repetition
according to their multiplicities).  This space is a finite ramified cover
over the space of rational maps of degree $d$ up to conjugacy, hence
$\mathrm{Rat}_d$ is a quasiprojective algebraic variety,
see~\cite{milnor1,silverman}. Note that in general it has
singularities. Denote by $T_k$ the positive closed $(1,1)$ current
describing the bifurcation of the marked critical point $c_k$ as
defined in Section~\ref{sec:bif}.

It is a consequence of~\cite[Theorem~2.2]{McM2}, that the set of
critically finite maps is the union of countably many points, together
with a quasi-projective curve $\mathcal{L}$. This curve appears only
in degrees which are squares of integers, and consists in the set of
Latt{\`e}s maps whose lift to some elliptic curve is the multiplication
by $\sqrt{d}$.
\begin{thm}
  For each $l<2d-2$, denote by $\Gamma_l$ an irreducible component of
  the quasiprojective subvariety of $\mathrm{Rat}_d\setminus
  \mathcal{L}$ consisting of rational maps such that the critical
  points $c_1, \cdots , c_l$ are preperiodic. 
  
  The codimension of $\Gamma_l$ is precisely $l$; the restriction of
  the current $T_{l+1} |_{\Gamma_l}$ is a non-zero positive closed
  $(1,1)$ current; and for any sequence of integers $k(n) <n$, we have
\begin{equation}\label{e:777}
\lim_{n\to\infty}  (d^{n}+d^{k(n)})^{-1}  \left[\Per_{l+1}(n,k(n)) \cap \Gamma_l
\right] = T_{l+1} |_{\Gamma_l}~,
\end{equation} 
where $\Per_{l+1}(n,k(n))$ denotes the set $\set{R_\lambda^n \, c_{l+1}(\lambda) 
= R_\lambda^{k(n)} \, c_{l+1}(\lambda)}$. 
\end{thm}
Another version of this result  has recently been obtained by Bassanelli-Berteloot, see~\cite{BB2}.
Their method is different although it bears the same analytic techniques: it is
based on the description of the bifurcation current in terms of the Lyapunov
function on the parameter space.
\begin{proof}
  For the proof we fix a decreasing sequence of components $\Gamma_1
  \supset \Gamma_2 \cdots \supset \Gamma_{2d-2}$ satisfying the
  conditions of the theorem.
  
  Being defined by $l$ equations, $\mathrm{codim}\, \Gamma_l \ge l$
  for all $l$.  In $\mathrm{Rat}_d\setminus \mathcal{L}$, the set of
  critically finite rational maps is countable, so any of its
  irreducible component is a point. This shows that $\Gamma_{2d-2}$ is
  necessarily a point. By induction, we conclude that
  $\mathrm{codim}\, \Gamma_l = l$. Moreover $\Gamma_{l+1}$ has
  codimension $1$ in $\Gamma_l$, so the family of rational maps
  restricted to $\Gamma_l$ with the marked critical point $c_{l+1}$ is
  non-trivial. Theorem~\ref{thm:cvg} implies~\eqref{e:777}. Finally,
  $T_{l+1} |_{\Gamma_l}$ is non-zero by Theorem~\ref{t:quasiproj}.
\end{proof}

Let now $1<l\leq 2d-2$; since the currents $T_i$ have continuous local
potentials, it is possible to take their wedge product
$T_1\wedge\cdots\wedge T_l$ (see \cite{de}). This defines a positive closed
current of bidegree $(l,l)$ of locally finite mass, which does not
charge analytic subsets. From the previous theorem we get the
following corollary, see also~\cite{BB2}.

\begin{cor}
For any $(n_1, \ldots, n_l)\in (\mathbb{N}^*)^l$, choose a collection
of integers $k(n_1)<n_1$,..., $k(n_l)<n_l$. Let $W_{n_1, \ldots,
  n_l}$  be the subvariety in  $\mathrm{Rat}_d\setminus
  \mathcal{L}$ defined  as $W_{n_1, \ldots,
  n_l}=\bigcap_{j=1}^l\Per_{j}(n_j,k(n_j))$. Then 
 $W_{n_1, \ldots,  n_l}$ has pure codimension $l$ and 
$$\lim_{n_l\cv\infty}\cdots  \lim_{n_1\cv\infty}
\unsur{(d^{n_l}+d^{k(n_l)}) \cdots (d^{n_1}+d^{k(n_1)})} [W_{n_1, \ldots,
  n_l}] \cv T_1\wedge\cdots \wedge T_l.$$ 
\end{cor}
 
Observe that the same holds when $\set{1,\ldots, l}$ is replaced by
any set of $l$ distinct integers $\set{i_1,\ldots, i_l}$. 
The proof is a straightforward adaptation of 
Corollary~\ref{cor:cvgall}). As a consequence we obtain that near any point of 
$\supp(T_1\wedge\cdots \wedge T_l)$ there are parameters where
$c_1,\ldots,c_l$ are periodic (resp. strictly preperiodic). For
$l=2d-2$, these parameters are critically finite (compare with
\ref{cor:meas}).   Notice that since the currents and measures
constructed here do not charge subvarieties, we need not care about
the curve $\el$ of flexible Latt{\`e}s examples. 
\medskip

Lastly, it follows from the convergence Theorem~\ref{thm:cvg}
  that for each $j$, $T_j\wedge T_j=0$. Indeed, $T_j$ is the limit of
  the sequence of hypersurfaces $\unsur{d^n+1}[\Per_j(n,0)]$, and 
  $[\Per_j(n,0)]\wedge T_j=0$. So by using the formulas for
  the Lyapounov exponent of the maximal entropy measure given in
  \cite{DeM2, bas-ber}, we can use the same argument as in
  Section~\ref{sec:poly} to get that for every $1\leq k\leq 2d-2$, the 
 support of the positive closed $(k,k)$ current  $(dd^c
  \mathrm{Lyap}(f))^k$, studied by \cite{bas-ber} is
   accumulated by the set of rational maps where $k$ critical points
   are preperiodic.

\end{document}